\documentclass[12pt]{article}
\usepackage[english]{babel}
\usepackage[letterpaper,hscale=.76,vscale=.76]{geometry}
\usepackage{amsmath,amssymb,amsthm,graphicx,float,array,xcolor,bm,tikz}
\usepackage[colorlinks=true,allcolors=blue,pdfusetitle]{hyperref}
\usepackage{longtable,caption,subcaption,placeins}
\newtheorem{theorem}{Theorem}[section]
\theoremstyle{definition}\newtheorem{definition}[theorem]{Definition}
\theoremstyle{plain}\newtheorem{proposition}[theorem]{Proposition}\newtheorem{lemma}[theorem]{Lemma}
\title{Numerical Investigation of Shortest Escape Path of Quadrilaterals In Bellman's Lost-in-a-forest Problem}
\author{Zhipeng Deng}
\date{}
\makeatletter
\renewcommand\section{\@startsection{section}{1}{\z@}%
  {-3.5ex \@plus -1ex \@minus -.2ex}%
  {2.3ex \@plus .2ex}%
  {\normalfont\Large\bfseries\raggedright}}
\makeatother
\begin{document}
\maketitle
\begin{abstract}

Bellman's lost-in-a-forest problem asks for a shortest prescribed path that reaches the boundary of a known region when the initial position and orientation are unknown. In this paper, we develop numerical investigation of convex and concave quadrilateral forests based on the developed general solution before, including rectangles, rhombi, parallelograms, isosceles trapezoids, and general simple quadrilaterals. For convex forests, translation is eliminated through support functions and an equilibrium multiplier polytope, yielding a finite family of competing containment constraints. For concave forests, support functions alone lose the reentrant geometry. We instead partition the polygon into two triangles and impose containment of every complete polyline segment. Ordering the vertices by their projections onto the internal diagonal's normal reduces the fitting calculation to a finite collection of three-variable linear programs, with residual corrected dual upper bounds.  The path search combines structured starting curves, derivative free exploration, and local constrained optimization. The results identify diagonal, strip-like, and bent-path regimes, resolve competing branches near geometric transitions, and quantify the effect of a concave notch relative to its triangular convex hull. Comparisons with isosceles and general triangles explain the behavior of truncated triangles and concave quadrilaterals. 
\end{abstract}



\setcounter{tocdepth}{2}
{\small\tableofcontents}
\medskip
\section{Introduction}
Bellman's lost-in-a-forest problem asks for a shortest prescribed path that reaches the boundary of a forest of known shape and size when the initial position and heading are unknown \cite{bellman,finch}. Equivalently, a path guarantees escape when no rigidly moved copy fits wholly inside the forest, so it also relates to universal coverage problem.  Classical studies address line search \cite{isbell,joris,finchtranslation}, infinite strips and arcs of prescribed width \cite{zalgaller,adhikaripitman}, and triangular forests \cite{besicovitch,coultonmovshovich,movshovichwetzel,movshovich2012,gibbs}. The surveys \cite{finch,fitscovers} explain the progress and connections with geometric fitting, universal covers, and Moser's worm problem. Approximation algorithms for forests of unknown shape are studied in past paper \cite{kubellangetepe}. 

Previously, we developed general boundary-based formulation and solution \cite{denggeneral}, equivalence and convergence arguments with further variants \cite{dengproof}, and a functional formulation for strip forests \cite{dengstrip}. Support function methods and numerical calculations for general triangles appear in \cite{dengtriangle}. The golden-gnomon solution by Temerev \cite{temerevdoria} and balanced support calibrations \cite{dengcalibration} provide additional triangular comparisons. Recent research also considered quadrilateral Moser worm covering \cite{Keller}. However, research gap exists that there is a significant lack of research on forest shape of quadrilaterals. Moreover, quadrilaterals include non-convex forest shapes, which has rarely been addressed in previous research.

Our method in this paper is partly based on the polygon extension in Section~2.5 of \cite{dengtriangle}. It eliminates translation through support functions and nonnegative equilibrium multipliers. Four side normals introduce competing dual branches, while a concave quadrilateral requires a different fitting calculation: global supports retain only its convex hull and erase its notch. We therefore partition the concave region into two triangles and impose containment of every complete polyline segment in their union. This extends the numerical treatment beyond the convex triangular forests considered in \cite{dengtriangle}. 

This study covers rectangles, rhombi, parallelograms, isosceles trapezoids, and general simple convex and concave quadrilaterals. Self-intersecting quadrilaterals are excluded before the path calculation. The emphasis and purose of this paper is on accurately evaluated paths on the quadrilateral geometry.

\section{Mathematical formulation and numerical methodology}
\subsection{Escape paths as a noncontainment problem}
Based on our previous paper \cite{denggeneral}, let $Q\subset\mathbb{R}^2$ be a compact polygonal forest
with nonempty interior. A prescribed path is
$
C=\gamma([0,1]),\gamma(0)=0,
$
of a rectifiable curve $\gamma:[0,1]\to\mathbb{R}^2$. Its length is denoted by $L(C)$. For $t\in[0,2\pi)$, let
\begin{equation}\label{eq:rotation}
R_t=\begin{pmatrix}\cos t&-\sin t\\\sin t&\cos t\end{pmatrix}.
\end{equation}
If the hiker starts at $z\in Q$ with relative orientation $t$, the walked set is $z+R_tC$.

Reaching $\partial Q$ is regarded as escape. Thus a path guarantees escape if no rotated copy can be
translated wholly into the interior of $Q$. For optimization of the infimum it is convenient to use the closed limiting condition. 

\begin{definition}[Containment scale]\label{def:scale}
For a fixed nontrivial path $C$ and orientation $t$, define
\begin{equation}\label{eq:scale}
\sigma_Q(C,t)=\sup\bigl\{s\geq0:\exists z\in\mathbb{R}^2\text{ such that }z+sR_tC\subseteq Q\bigr\}.
\end{equation}
The worst orientation fit scale is
\begin{equation}\label{eq:worstscale}
\sigma_Q^*(C)=\sup_{t\in[0,2\pi)}\sigma_Q(C,t).
\end{equation}
\end{definition}
Because $L(sC)=sL(C)$, the physical length of the scaled path that reaches the noncontainment
threshold is
\begin{equation}\label{eq:physical}
J_Q(C)=L(C)\,\sigma_Q^*(C).
\end{equation}
The quantity $J_Q(C)$ is invariant under positive rescaling of the raw candidate $C$. 

\subsection{Convex support formulation and polyline optimization}
The convex formulation below is partly based on Section~2.5 of \cite{dengtriangle}. That section extends the triangular support constraint to a polygon represented as an intersection of supporting half-planes, with nonnegative equilibrium multipliers, and gives the corresponding polygonal chain optimization. We use that construction for convex quadrilaterals, reduce it to extreme multiplier branches, and evaluate the fitting scale throughout the orientation circle. A concave polygon is not directly represented by that half-plane intersection. Its whole segment fitting formulation is developed in Section~\ref{sec:concave-fitting}.

Assume first that $Q$ is convex
and represented by its supporting half-planes \cite{schneider},
\begin{equation}\label{eq:halfplanes}
Q=\bigl\{x\in\mathbb{R}^2:u_j^{\mathsf T} x\leq b_j,\ j=1,\ldots,m\bigr\},\qquad \|u_j\|_2=1.
\end{equation}
Write
$
u_j=\mathbf{e}(\phi_j),\mathbf{e}(\theta)=(\cos\theta,\sin\theta)^{\mathsf T}.
$
The support function of the path is
\begin{equation}\label{eq:support}
h_C(\theta)=\max_{x\in C}x^{\mathsf T} \mathbf{e}(\theta).
\end{equation}
Since $0\in C$, $h_C(\theta)\geq0$. Rotation satisfies
\begin{equation}\label{eq:supportrotation}
h_{R_tC}(\phi)=h_C(\phi-t).
\end{equation}
For a fixed $t$ and scale $s$, the translated copy $z+sR_tC$ is contained in $Q$ exactly when
\begin{equation}\label{eq:translation}
u_j^{\mathsf T} z+s\,h_C(\phi_j-t)\leq b_j,\qquad j=1,\ldots,m.
\end{equation}
The only variables in this feasibility problem are the two translation coordinates. Define the
normalized equilibrium multiplier polytope
\begin{equation}\label{eq:lambda}
\Lambda=\left\{\lambda\in\mathbb{R}_+^m:\mathbf{1}^{\mathsf T}\lambda=1,\ \sum_{j=1}^m\lambda_ju_j=0\right\}.
\end{equation}
\begin{proposition}[Farkas containment criterion]\label{prop:farkas}
For a bounded convex polygon $Q$, the system in Eq.~\eqref{eq:translation} is
feasible if and only if
\begin{equation}\label{eq:farkas}
s\sum_{j=1}^m\lambda_jh_C(\phi_j-t)\leq\sum_{j=1}^m\lambda_jb_j\qquad\text{for every }\lambda\in\Lambda.
\end{equation}
The obstruction to fitting a translated copy wholly in the interior is
\begin{equation}\label{eq:obstruction}
\max_{\lambda\in\Lambda}\sum_{j=1}^m\lambda_j\bigl[s\,h_C(\phi_j-t)-b_j\bigr]\geq0.
\end{equation}
\end{proposition}
\begin{proof}
Put $U\in\mathbb{R}^{m\times2}$ with rows $u_j^{\mathsf T}$ and $c=b-sh(t)$, where $h(t)_j=h_C(\phi_j-t)$. The system
$Uz\leq c$ is feasible if and only if $\lambda^{\mathsf T} c\geq0$ for every $\lambda\geq0$ satisfying $U^{\mathsf T}\lambda=0$. Positive scaling of $\lambda$ is
immaterial, so nonzero multipliers may be normalized by $\mathbf{1}^{\mathsf T}\lambda=1$, giving Eq.~\eqref{eq:lambda} and Eq.~\eqref{eq:farkas}. Negating
the strict interior feasibility criterion gives Eq.~\eqref{eq:obstruction}.
\end{proof}
Let
$\operatorname{ext}(\Lambda)=\{\lambda^{(1)},\ldots,\lambda^{(N_\Lambda)}\}.$
As the expression in Eq.~\eqref{eq:obstruction} is linear in $\lambda$, only extreme multipliers are needed. Define
\begin{equation}\label{eq:branches}
B_\nu=\sum_{j=1}^m\lambda_j^{(\nu)}b_j,\qquad
A_\nu(t;C)=\sum_{j=1}^m\lambda_j^{(\nu)}h_C(\phi_j-t),
\end{equation}
and, for branches with $B_\nu>0$,
\begin{equation}\label{eq:margin}
G_C(t)=\max_{1\leq\nu\leq N_\Lambda}\frac{A_\nu(t;C)}{B_\nu}.
\end{equation}
Then the largest scale that can be fitted at orientation $t$ is
\begin{equation}\label{eq:reciprocal}
\sigma_Q(C,t)=\frac{1}{G_C(t)},
\end{equation}
with the usual extended-value interpretation if $G_C(t)=0$. Consequently,
\begin{equation}\label{eq:objective}
J_Q(C)=\frac{L(C)}{\displaystyle\min_{t\in[0,2\pi)}G_C(t)}.
\end{equation}

The continuous curve is
approximated by an open polyline
\begin{equation}\label{eq:polyline}
P=(P_1,\ldots,P_K),\qquad P_1=0,
\end{equation}
with increments
$
e_i=P_{i+1}-P_i=(d_{x,i},d_{y,i})^{\mathsf T}, i=1,\ldots,K-1.
$
The vertices are recovered cumulatively and the length is
\begin{equation}\label{eq:length}
L(P)=\sum_{i=1}^{K-1}\|e_i\|_2=\sum_{i=1}^{K-1}\sqrt{d_{x,i}^2+d_{y,i}^2}.
\end{equation}
For any direction $\theta$, the support function is evaluated \emph{exactly} from the vertices,
\begin{equation}\label{eq:polylinesupport}
h_P(\theta)=\max_{1\leq i\leq K}P_i^{\mathsf T} \mathbf{e}(\theta).
\end{equation}
No sampling along individual segments is necessary because a linear functional attains its maximum
over a segment at an endpoint.

For convex $Q$, substitution of Eq.~\eqref{eq:polylinesupport} into Eq.~\eqref{eq:objective} gives a nonsmooth, nonconvex objective in the edge
variables. The raw scale is redundant. The numerical normalization is
\begin{equation}\label{eq:penalty}
\mu(P)=\min_tG_P(t),\qquad\widehat P=\frac{P}{\mu(P)},
\qquad\min_tG_{\widehat P}(t)=1.
\end{equation}
The final physical path therefore has length $L(P)/\mu(P)$.

The outer search uses multiple path families and exploratory differential
evolution \cite{stornprice}, followed by local refinement. Convex candidates
are refined by sequential constrained minimization: at each orientation,
the supporting vertices and an active extreme dual branch are fixed from
the current path, and the Euclidean polyline length is minimized subject
to the resulting linear support inequalities. The numerical sequence is
\begin{equation}\label{eq:search}
\text{multistart exploration}\longrightarrow\text{local shape refinement}\longrightarrow\text{continuous angle fit check}.
\end{equation}
For the convex refinement, sequential least squares quadratic programming uses analytic length gradients. 

For concave candidates, changing segment modes makes the outer objective less suitable for the convex support contact refinement. Multiple geometric starts are compared using the whole segment fitting routine in Section~\ref{sec:concave-fitting}. Derivative free shape exploration includes the Nelder--Mead method; promising proposals are evaluated on the complete orientation grid and compared using their final guarded lengths. The additional refinement of the selected general quadrilaterals is described with those results in Section~3. 

\subsection{Quadrilateral geometry: convex and concave cases}\label{sec:quadrilateral-parameters}
A nondegenerate bounded triangle formed by three straight lines is necessarily convex: it is the intersection of the three inward half-planes. A simple quadrilateral can instead be convex or nonconvex (concave), with one interior angle greater than $\pi$. Thus the fourth side changes both the convex dual formulation and the possible geometry of the forest. In particular, the three-line support method cannot be applied unchanged to a concave quadrilateral, whose notch must remain part of the containment test.

As shown in Figure 1, we assume rigid motion and scale are removed by
fixing the first side to unit length. We use the ordered vertices
\begin{equation}\label{eq:basevertices}
A=(0,0),\qquad B=(1,0),
\end{equation}
\begin{equation}\label{eq:quadvertices}
C=\bigl(1-q\cos\beta,q\sin\beta\bigr),\qquad
D=\bigl(p\cos\alpha,p\sin\alpha\bigr),
\end{equation}
where $p,q>0$ and $0<\alpha,\beta<\pi$. Thus $|AD|=p$, $|BC|=q$, the ray $AD$ makes angle $\alpha$ with $AB$,
and the ray $BC$ makes angle $\pi-\beta$ with the positive $x$-axis. Let $\operatorname{orient}(X,Y,Z)=\det(Y-X,Z-X)$. Before any convexity test, the nonadjacent edge
pairs $AB$--$CD$ and $BC$--$DA$ are tested for proper intersection; a proper intersection produces a
self-intersecting quadrilateral and is discarded. Zero-area and self-touching degeneracies are likewise excluded.

For a simple counterclockwise quadrilateral, the four consecutive signed turns are
\begin{align}
\kappa_A&=p\sin\alpha,\label{eq:turna}\\
\kappa_B&=q\sin\beta,\label{eq:turnb}\\
\kappa_C&=q\bigl[\sin\beta-p\sin(\alpha+\beta)\bigr],\label{eq:turnc}\\
\kappa_D&=p\bigl[\sin\alpha-q\sin(\alpha+\beta)\bigr].\label{eq:turnd}
\end{align}
Since $\kappa_A,\kappa_B>0$ under the stated parameter ranges, strict convexity is equivalent to
\begin{equation}\label{eq:convexity}
\sin\beta>p\sin(\alpha+\beta),\qquad\sin\alpha>q\sin(\alpha+\beta).
\end{equation}
If exactly one of $\kappa_C,\kappa_D$ is negative, the corresponding vertex is the unique reflex vertex. The main structured families are summarized in Table~\ref{tab:families}.

\begin{table}[tbp]
\centering
\caption{Special quadrilateral families in the normalized parameterization.}\label{tab:families}
\small
\begin{tabular}{lll}
Family&Parameter restrictions&Additional geometry\\
\hline
Rectangle&$p=q=r,\ \alpha=\beta=\pi/2$&size $1\times r$\\
Rhombus&$p=q=1,\ \beta=\pi-\alpha$&all sides $1$\\
Parallelogram&$p=q=s,\ \beta=\pi-\alpha$&side lengths $1,s$\\
Isosceles trapezoid&$p=q=s,\ \beta=\alpha$&top base $1-2s\cos\alpha$
\end{tabular}
\end{table}

\begin{figure}[tbp]
\centering
\scalebox{1.5}{
\begin{tikzpicture}[scale=1.35,>=stealth]
\definecolor{forestred}{HTML}{D8202A}
\definecolor{forestlight}{HTML}{FCE8E9}
\coordinate (A) at (0,0);
\coordinate (B) at (2,0);
\coordinate (D) at (.76,1.43);
\coordinate (C) at (1.44,1.25);
\filldraw[draw=forestred,fill=forestlight,line width=1.35pt] (A)--(B)--(C)--(D)--cycle;
\draw[dashed] (-.12,0)--(2.28,0);
\draw[->,thick] (A)--(.43,.809);
\draw[->,thick] (B)--(1.68,.714);
\node[below left] at (A) {$A$};
\node[below right] at (B) {$B$};
\node[above right] at (C) {$C$};
\node[above left] at (D) {$D$};
\node[left] at (.33,.63) {$p$};
\node[right] at (1.78,.63) {$q$};
\node at (.27,.12) {$\alpha$};
\node at (1.73,.12) {$\beta$};
\end{tikzpicture}}
\caption{Normalized quadrilateral parameterization. The base $AB$ has unit length,
$AD$ has length $p$, $BC$ has length $q$, and the two side directions are controlled by $\alpha$
and $\beta$.}\label{fig:parameterization}
\end{figure}
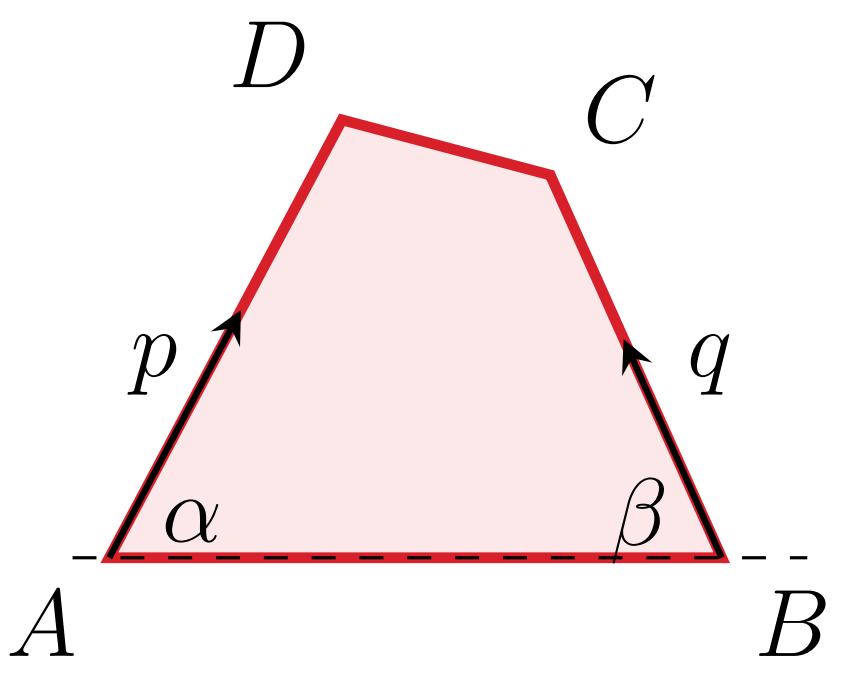

For a generic convex quadrilateral, $m=4$. The matrix in Eq.~\eqref{eq:lambda} has three independent equality
constraints, so $\Lambda_4$ is generically one-dimensional: a line segment. Its two endpoints typically
correspond to different positive circuits of three side normals. Hence
\begin{equation}\label{eq:extreme}
\max_{\lambda\in\Lambda_4}(\cdots)=\max_{\lambda\in\operatorname{ext}(\Lambda_4)}(\cdots)
\end{equation}
is a genuine finite maximum. Which branch is active can change with the orientation of the path.
Imposing all branch inequalities simultaneously would incorrectly replace this maximum by a
minimum-type condition and would solve a more restrictive problem.

Every nondegenerate extreme point of $\Lambda_4$ has at most three positive components. For a triple
$\{i,j,k\}$, a null vector is obtained from the planar determinants
\begin{equation}\label{eq:circuits}
\widetilde\lambda_i=\det(u_j,u_k),\qquad
\widetilde\lambda_j=\det(u_k,u_i),\qquad
\widetilde\lambda_k=\det(u_i,u_j).
\end{equation}
When these coefficients have a common nonzero sign, they can be sign-adjusted and normalized
to sum to one, producing an extreme positive circuit. Antipodal pairs of normals produce the
additional two-element circuit $(1/2,1/2)$ on that pair. This enumeration is faster than resolving an
inner linear program at every orientation and every outer iteration.

For the parameterization Eq.~\eqref{eq:quadvertices}, define
\begin{equation}\label{eq:delta}
\Delta_x=p\cos\alpha-1+q\cos\beta,\quad
\Delta_y=p\sin\alpha-q\sin\beta,\quad
\ell_3=\sqrt{\Delta_x^2+\Delta_y^2}.
\end{equation}
When the polygon is convex and counterclockwise, outward unit normals and offsets can be chosen
as
\begin{align}
u_1&=(0,-1),& b_1&=0,\label{eq:normal1}\\
u_2&=(\sin\beta,\cos\beta),& b_2&=\sin\beta,\label{eq:normal2}\\
u_3&=\frac{1}{\ell_3}(\Delta_y,-\Delta_x),&
b_3&=\frac{p[\sin\alpha-q\sin(\alpha+\beta)]}{\ell_3},\label{eq:normal3}\\
u_4&=(-\sin\alpha,\cos\alpha),& b_4&=0.\label{eq:normal4}
\end{align}
Thus
\begin{equation}\label{eq:quadintersection}
Q=\bigcap_{j=1}^4\{x:u_j^{\mathsf T} x\leq b_j\}.
\end{equation}
The equilibrium set is
\begin{equation}\label{eq:lambda4}
\Lambda_4=\left\{\lambda\in\mathbb{R}_+^4:\mathbf{1}^{\mathsf T}\lambda=1,\ \sum_{j=1}^4\lambda_ju_j=0\right\}.
\end{equation}
For each extreme point $\lambda^{(\nu)}$, the normalized branch is
\begin{equation}\label{eq:normalizedbranch}
E_\nu(t;P)=\frac{\sum_{j=1}^4\lambda_j^{(\nu)}h_P(\phi_j-t)}{\sum_{j=1}^4\lambda_j^{(\nu)}b_j}.
\end{equation}
The fit margin is $G_P(t)=\max_\nu E_\nu(t;P)$ and the physical length follows from Eq.~\eqref{eq:objective}.

Computationally, $\operatorname{ext}(\Lambda_4)$ is constructed once for each quadrilateral. This removes the translation
variables and avoids repeatedly solving a linear program inside the outer path search. 

\subsection{Whole segment fitting in concave quadrilaterals}\label{sec:concave-fitting}
The support formulation above is exact only because a convex polygon is an intersection of half-planes. A strictly
concave quadrilateral cannot satisfy
\[
Q=\bigcap_{j=1}^4\{x:u_j^{\mathsf T} x\leq b_j\},
\]
since any intersection of half-planes is convex. There is a second, more fundamental obstruction: for
every compact set $S$,
\begin{equation}\label{eq:hullsupport}
h_S=h_{\operatorname{conv}S}.
\end{equation}
Thus a support-only representation cannot tell whether a segment of the path cuts through a
reentrant notch of a concave forest.

Every simple concave quadrilateral has exactly one reflex vertex and admits an exact triangulation
through that vertex. In the present vertex order,
\begin{equation}\label{eq:triangulationc}
\kappa_C<0\quad\Longrightarrow\quad Q=\triangle ABC\cup\triangle ACD,
\end{equation}
where the common internal diagonal is $AC$, while
\begin{equation}\label{eq:triangulationd}
\kappa_D<0\quad\Longrightarrow\quad Q=\triangle ABD\cup\triangle BCD,
\end{equation}
with common diagonal $BD$.

Let the two counterclockwise triangles be
\begin{equation}\label{eq:trianglehalfplanes}
T_r=\{x:A_rx\leq c_r\},\qquad r=1,2,
\end{equation}
where each $A_r$ has three rows. Point membership is the disjunction
\begin{equation}\label{eq:disjunction}
x\in Q\quad\Longleftrightarrow\quad x\in T_1\lor x\in T_2.
\end{equation}
Checking only polyline vertices is not sufficient, because two vertices can lie in $Q$ while their joining
segment passes through the exterior notch.

\begin{lemma}[Segment containment in a two triangle quadrilateral]\label{lem:segment}
Let $T_1$ and $T_2$ be the two
triangles in Eq.~\eqref{eq:triangulationc} or Eq.~\eqref{eq:triangulationd}, sharing the internal diagonal $D_{\mathrm{int}}$. A segment $[X,Y]$ lies in $T_1\cup T_2$ if and only if one
of the following holds: (i) $X,Y\in T_1$; (ii) $X,Y\in T_2$; (iii) $X\in T_1$, $Y\in T_2$, and $[X,Y]$ intersects
$D_{\mathrm{int}}$; or (iv) the symmetric alternative with $1$ and $2$ interchanged.
\end{lemma}
\begin{proof}
The same triangle cases follow from convexity of each $T_r$. In a cross triangle case, if
$Z\in[X,Y]\cap D_{\mathrm{int}}$, then $[X,Z]\subseteq T_1$ and $[Z,Y]\subseteq T_2$, again by convexity. Conversely, a connected
segment passing from the relative interior of one triangle to the relative interior of the other must
pass through their intersection, which is the common diagonal (up to shared endpoints).
\end{proof}

Fix a raw polyline $P=(P_1,\ldots,P_K)$ and orientation $t$. The placed
vertices at scale $s$ are
\begin{equation}\label{eq:placedvertices}
X_i=z+sR_tP_i,\qquad i=1,\ldots,K.
\end{equation}
Let $\ell_i\in\{1,2\}$ indicate a triangle containing $X_i$. The segment
mode is the ordered pair
\begin{equation}\label{eq:modes}
m_i=(\ell_i,\ell_{i+1})\in\{(1,1),(2,2),(1,2),(2,1)\}.
\end{equation}
We write these four modes as $11,22,12,21$. Once the labels are fixed,
endpoint containment is expressed by the linear inequalities
$A_{\ell_i}(z+sR_tP_i)\leq c_{\ell_i}$. Modes $11$ and $22$ require
only these endpoint inequalities. The cross triangle modes also require
intersection with the finite shared diagonal.

Write that diagonal as $[d_0,d_1]$, set $r_i=R_tP_i$, and let
$w_i=r_{i+1}-r_i$. For zero-length segments, endpoint conditions suffice;
the following crossing constraints apply to nonzero segments. Define
\begin{equation}\label{eq:crosspoint}
F_i(d;z,s)=\frac{\det(w_i,d-z)+s\det(r_i,r_{i+1})}{\|w_i\|_2}.
\end{equation}
It equals $\det(w_i,d-X_i)/\|w_i\|_2$ and is affine in
$(z_x,z_y,s)$. Let
\begin{equation}\label{eq:diagonalpoint}
\varepsilon=\operatorname{sign}\det(d_1-d_0,\bar c_1-d_0),
\end{equation}
where $\bar c_1$ is the centroid of $T_1$. In mode $12$, impose
\begin{equation}\label{eq:crossingequality}
\varepsilon F_i(d_0;z,s)\leq0,\qquad
\varepsilon F_i(d_1;z,s)\geq0;
\end{equation}
in mode $21$, reverse both signs. The endpoint constraints already
force a crossing of the line through the diagonal. Eq.~\eqref{eq:crossingequality}
places that crossing on the finite segment $[d_0,d_1]$. The normalization
by $\|w_i\|_2$ keeps the crossing residual meaningful for short path
segments. Thus complete segment containment is expressed using only
linear constraints in three placement variables.

The possible vertex labels have a simple ordering. Choose a unit normal
$n_d$ to the internal diagonal pointing toward $T_1$. The triangles lie
on opposite sides of its supporting line. For $s>0$, define
\begin{equation}\label{eq:diagonal-threshold}
\xi_i=n_d^{\mathsf T}R_tP_i,\qquad
\eta=\frac{n_d^{\mathsf T}d_0-n_d^{\mathsf T}z}{s}.
\end{equation}
Then
\[
n_d^{\mathsf T}(X_i-d_0)=s(\xi_i-\eta).
\]
A vertex with $\xi_i>\eta$ must belong to $T_1$, while one with
$\xi_i<\eta$ must belong to $T_2$. If $\xi_i=\eta$, a contained vertex
lies on the shared diagonal and belongs to both triangles. Sorting the
fixed projections as $\xi_{\pi(1)}\leq\cdots\leq\xi_{\pi(K)}$ therefore
leaves only $K+1$ label patterns:
\begin{equation}\label{eq:ordered-cuts}
\ell_{\pi(i)}^{(k)}=
\begin{cases}
2,&i\leq k,\\
1,&i>k,
\end{cases}
\qquad k=0,\ldots,K.
\end{equation}
Every positive scale placement determines one of these cuts. Equal
projections cause no loss of feasible placements because vertices on
the diagonal may be assigned to either triangle. Conversely, the
endpoint and crossing constraints for a cut ensure containment of every
complete segment by Lemma~\ref{lem:segment}. 

Finite bounds on the placement variables are available directly from the
geometry. After choosing one path vertex as the origin, its placed point
is $z$, so $z$ lies in the coordinate bounding box of $Q$. With
$
D_Q=\operatorname{diam}(Q),D_P=\max_{i,j}\|P_i-P_j\|_2,
$
every feasible placement satisfies
\begin{equation}\label{eq:scalebound}
0\leq s\leq\frac{D_Q}{D_P}\qquad(D_P>0).
\end{equation}
The path anchor is chosen to reduce the largest distance to another path
vertex. Write $y=(z_x,z_y,s)^{\mathsf T}$ and let
$\underline y\leq y\leq\overline y$ denote these finite box bounds.
For each cut $k$, collect the endpoint and whole segment crossing
constraints in $A^{(k)}y\leq b^{(k)}$. The fixed orientation fit is
\begin{equation}\label{eq:milp}
\begin{aligned}
\sigma_Q(P,t)&=\max_{0\leq k\leq K}s_k(P,t),\\
s_k(P,t)&=\max_y\left\{e_s^{\mathsf T}y:
A^{(k)}y\leq b^{(k)},\quad
\underline y\leq y\leq\overline y\right\},
\end{aligned}
\end{equation}
where $e_s=(0,0,1)^{\mathsf T}$. The calculation is the maximum over
$K+1$ three variable linear programs. The zero scale case is included by
closure. This finite reduction is exact in real arithmetic for a specified polyline and orientation with linear programming.

We also evaluate an explicit upper bound for each linear program. For
one cut, abbreviate its constraints as $Ay\leq b$ and choose any
nonnegative multiplier vector $\lambda$. Define
\begin{equation}\label{eq:residual-dual-upper}
\begin{aligned}
r&=e_s-A^{\mathsf T}\lambda,\qquad \lambda\geq0,\\
U&=\lambda^{\mathsf T}b+
\sum_{j=1}^{3}\max\{r_j\underline y_j,r_j\overline y_j\}.
\end{aligned}
\end{equation}
For every feasible $y$,
$e_s^{\mathsf T}y=\lambda^{\mathsf T}Ay+r^{\mathsf T}y\leq U$.
Consequently, Eq.~\eqref{eq:residual-dual-upper} remains an upper bound when the numerical multipliers have a nonzero stationarity residual. It follows from nonnegative multipliers and the finite box bounds, without relying on an asserted solver optimum. The implementation uses the linear program dual output to choose $\lambda$, enforces nonnegativity, and evaluates the residual correction with extended precision and a floating point allowance. 

Whole segment containment is checked independently by intersecting each segment with the two triangles and verifying that their parameter intervals cover the complete segment. An additional Boolean test uses the two non-reflex edge half-planes and the union of the two reflex-edge half-planes; its maximum along a segment occurs at an endpoint or a switch of the two affine reflex edge expressions. Residual contact errors are removed by a small contraction toward a point in the kernel of $Q$. The feasible placement and the dual upper bound are retained separately throughout the orientation calculation.

For the free coordinate outer search, finite paths also admit a contact enumeration calculation. A path vertex on a polygon edge, or a path segment through the reflex vertex, defines a plane in the three translation--scale variables. Intersections of three such planes are tested by complete segment membership, and the largest feasible scale supplies an active fitting branch. Derivatives of these branches guide changes in all free path vertex coordinates. 

The two-vertex case admits a faster special treatment. For a fixed direction, each triangle gives an interval on a parallel line, and the longest connected component of their union determines the longest chord lying wholly in $Q$. Between consecutive projected polygon vertices, the chord length is piecewise linear, so the directional maximum is obtained by checking the corresponding breakpoints. 

Triangulation and the convex hull also give a comparison for the
true escape length $L^*(\cdot)$. Monotonicity under set inclusion gives
\begin{equation}\label{eq:inclusionbound}
\max\{L^*(T_1),L^*(T_2)\}\leq L^*(Q)
\leq L^*(\operatorname{conv}Q).
\end{equation}
Every path that guarantees escape from the convex hull also guarantees
escape from $Q$. The converse can fail substantially because a path
can encounter the notch sooner. The two triangles in this
triangulation are computational pieces of one forest: their fitting
scales cannot simply be maximized or added to recover the fit in
their union.

\subsection{Orientation resolution and continuous angle checks}
\label{sec:angular-checks}
Every final candidate is evaluated using the common orientation grid
\begin{equation}\label{eq:anglegrid}
t_k=\frac{2\pi k}{M},\qquad k=0,\ldots,M-1,\qquad
M=3\,600,\quad\Delta t=\frac{2\pi}{3\,600}=0.1^\circ.
\end{equation}
The corresponding
settings and results are reported in Section~3. 

For convex polygons, a support function is Lipschitz in angle. With
$R_P=\max_i\|P_i\|_2$,
\begin{equation}\label{eq:lipschitzsupport}
|h_P(\theta_1)-h_P(\theta_2)|\leq R_P|\theta_1-\theta_2|.
\end{equation}
For branch coefficients $\mu_j^{(\nu)}=\lambda_j^{(\nu)}/B_\nu$,
the margin $G_P$ consequently has a Lipschitz constant bounded by
\begin{equation}\label{eq:lipschitzconstant}
L_G=R_P\max_\nu\sum_j\mu_j^{(\nu)}.
\end{equation}
Writing $G_{P,M}=\min_kG_P(t_k)$ gives the coarse diagnostic bound
\begin{equation}\label{eq:continuumcertificate}
\min_tG_P(t)\geq G_{P,M}-\frac{L_G\Delta t}{2}.
\end{equation}

For concave quadrilaterals, the fitting modes can change with angle
and the convex support formula does not determine the answer.
Nevertheless, $H=\operatorname{conv}Q$ gives the directional bound
\begin{equation}\label{eq:hull-direction-bound}
\sigma_Q(P,t)\leq\sigma_H(P,t).
\end{equation}

The remaining angular intervals are controlled by a kernel ball
bound. Let $B(c,\rho)$ lie in the kernel of $Q$, so every point of
this ball sees the entire polygon. The largest such numerical ball
is obtained from the directed boundary half-planes, with a small
inward safety allowance. After choosing a path anchor $a$, set
\begin{equation}\label{eq:inverse-scale-lipschitz}
R=\max_i\|P_i-a\|_2,\qquad q(t)=\frac{1}{\sigma_Q(P,t)},\qquad
|q(t)-q(t')|\leq\frac{R}{\rho}|t-t'|.
\end{equation}
The inequality follows by contracting a feasible rotated path
toward the kernel ball to absorb the displacement caused by an
angular change. It controls the inverse fitting scale even when
the selected triangle modes change.

At each endpoint $a,b$ of an angular interval, let $u_a,u_b$ be
numerical upper bounds on the fit scale, obtained either from the
convex hull or from the cut-wise dual bounds. With $K_q=R/\rho$,
\begin{equation}\label{eq:inverse-cones}
q(t)\geq
\max\left\{\frac1{u_a}-K_q(t-a),\,
             \frac1{u_b}-K_q(b-t)\right\},\qquad a\leq t\leq b.
\end{equation}
The minimum of the right hand side is obtained from the endpoints and the intersection of its two affine pieces. Its positive reciprocal bounds the fit throughout that interval; the continuous convex hull upper bound supplies a second bound. Intervals whose upper bounds can still exceed the best feasible scale by the target tolerance are subdivided. The calculation therefore reports an observed fit scale $s_-$ and a guarded scale $s_+$, with lengths $L(P)s_-$ and $L(P)s_+$. The latter is used for the published escape-path length, and their relative difference measures the unresolved orientation error for that fixed candidate.

\FloatBarrier
\section{Numerical Results}
\label{sec:numerical}

The numerical study follows five families: square and rectangle, diamond,
parallelogram, isosceles trapezoid, and arbitrary convex and concave
quadrilaterals. The results compare
numerical candidate lengths and path geometries.

\subsection{Square and Rectangle}
\label{sec:rectangles}

The first family has $p=q=r$ and $\alpha=\beta=\pi/2$, giving a
$1\times r$ rectangle. The ten cases in Figure 2 are
$r=0.1,0.2,0.3,0.4,0.45,0.48,0.5,0.6,0.8,1$.
The square returns the diagonal length $\sqrt2$; the moderate-aspect-ratio
examples likewise return their diagonals as escape path. Thin rectangles select a strip-shaped escape path.

\FloatBarrier
\begingroup
\def\LTcaptype{figure}
\makeatletter\def\@captype{figure}\let\caption@stepcounter\@gobble\makeatother
\setlength{\LTpre}{5pt}
\setlength{\LTpost}{5pt}
\setlength{\LTcapwidth}{\linewidth}
\setlength{\tabcolsep}{0pt}
\captionsetup{position=bottom,skip=4pt}
\begin{longtable}{@{}c@{\hspace{0.009333\linewidth}}c@{\hspace{0.009333\linewidth}}c@{\hspace{0.009333\linewidth}}c@{}}
\endfirsthead
\endhead
\caption[]{Square and rectangle paths: 10 cases in parameter order. The aspect ratio $r$ increases from left to right and then from top to bottom. (continued on following page).}\\
\endfoot
\caption{Square and rectangle escape paths in 10 cases. }\label{fig:rectangles-1}\\
\endlastfoot
\begin{subfigure}[b]{0.243\linewidth}\centering
\includegraphics[width=\linewidth,height=100pt,keepaspectratio]{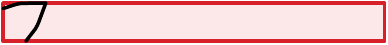}
\par{\footnotesize $r=0.1$\\$L=0.227829$\par}
\end{subfigure}
&
\begin{subfigure}[b]{0.243\linewidth}\centering
\includegraphics[width=\linewidth,height=100pt,keepaspectratio]{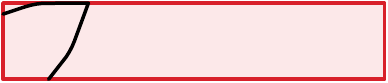}
\par{\footnotesize $r=0.2$\\$L=0.455658$\par}
\end{subfigure}
&
\begin{subfigure}[b]{0.243\linewidth}\centering
\includegraphics[width=\linewidth,height=100pt,keepaspectratio]{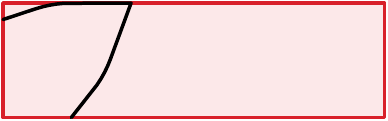}
\par{\footnotesize $r=0.3$\\$L=0.683488$\par}
\end{subfigure}
&
\begin{subfigure}[b]{0.243\linewidth}\centering
\includegraphics[width=\linewidth,height=100pt,keepaspectratio]{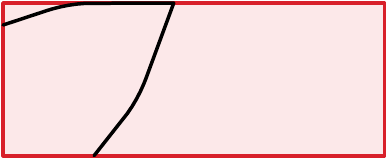}
\par{\footnotesize $r=0.4$\\$L=0.911317$\par}
\end{subfigure}
\\
\begin{subfigure}[b]{0.243\linewidth}\centering
\includegraphics[width=\linewidth,height=100pt,keepaspectratio]{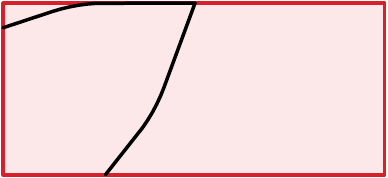}
\par{\footnotesize $r=0.45$\\$L=1.025231$\par}
\end{subfigure}
&
\begin{subfigure}[b]{0.243\linewidth}\centering
\includegraphics[width=\linewidth,height=100pt,keepaspectratio]{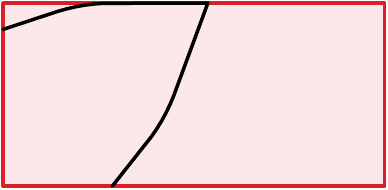}
\par{\footnotesize $r=0.48$\\$L=1.093580$\par}
\end{subfigure}
&
\begin{subfigure}[b]{0.243\linewidth}\centering
\includegraphics[width=\linewidth,height=100pt,keepaspectratio]{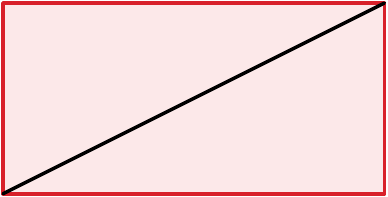}
\par{\footnotesize $r=0.5$\\$L=1.118034$\par}
\end{subfigure}
&
\begin{subfigure}[b]{0.243\linewidth}\centering
\includegraphics[width=\linewidth,height=100pt,keepaspectratio]{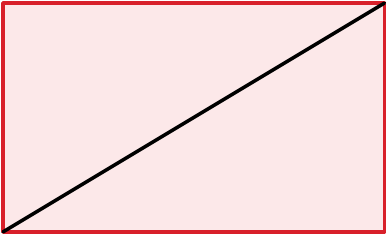}
\par{\footnotesize $r=0.6$\\$L=1.166190$\par}
\end{subfigure}
\\
\begin{subfigure}[b]{0.243\linewidth}\centering
\includegraphics[width=\linewidth,height=100pt,keepaspectratio]{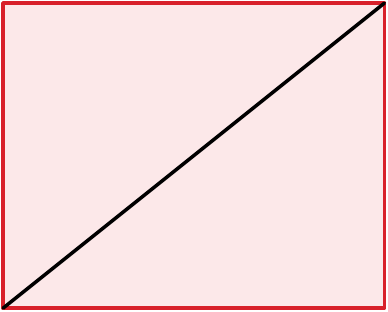}
\par{\footnotesize $r=0.8$\\$L=1.280625$\par}
\end{subfigure}
&
\begin{subfigure}[b]{0.243\linewidth}\centering
\includegraphics[angle=90,width=\linewidth,height=100pt,keepaspectratio]{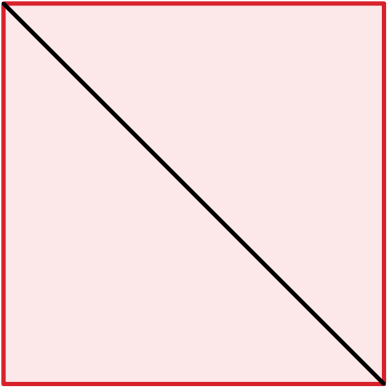}
\par{\footnotesize $r=1$\\$L=1.414214$\par}
\end{subfigure}
&&\\
\end{longtable}
\endgroup
\begin{figure}[htbp]
\centering
\includegraphics[width=.98\linewidth]{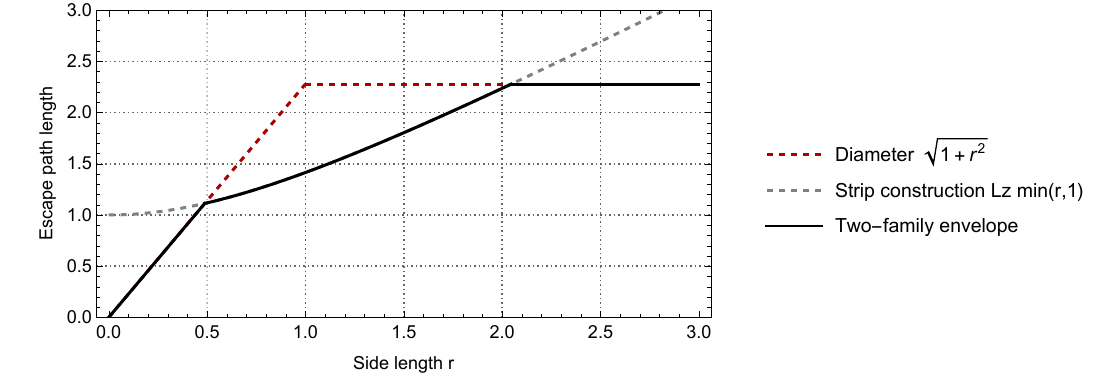}
\caption{Square and rectangle relationship plot over the full
range $0\le r\le3$. The black curve is the comparison envelope
$\min\{\ell_Z\min(r,1),\sqrt{1+r^2}\}$. The two transition locations are $r_c$ and $1/r_c$. Results for $r>1$ follow directly from those for $1/r<1$: a
$90^\circ$ rotation and scaling by $r$ give the corresponding path
and length $L(r)=rL(1/r)$.}
\label{fig:rectangle-relationship}
\end{figure}

The transition value $0.4885$ follows from the equality between the diagonal
and strip constructions. Strip constant
from the analytic strip expression gives
\begin{align}
\phi&=\arcsin\!\left[\frac16+\frac43\sin\!\left(\frac13\arcsin\frac{17}{64}\right)\right],\qquad
\psi=\arctan\!\left(\frac{\sec\phi}{2}\right),\\
\ell_Z&=2\left(\frac\pi2-\phi-2\psi+\tan\phi+\tan\psi\right)
       =2.278291641440\ldots,\\
r_c&=\sqrt{\frac{1}{\ell_Z^2-1}}
       =0.488496247834\ldots .
\end{align}

The analytic strip construction is developed by Zalgaller, and rectangle relationship is already discussed in~\cite{finch, dengstrip} as shown in Figure 3.


\subsection{Diamond}
\label{sec:diamonds}

Here ``diamond'' denotes the unit-side rhombus,
$p=q=1$ and $\beta=\pi-\alpha$. The gallery in Figure 4 retains the numerical results at each selected angle. And Figure 5 shows the numerical lengths for unit-side rhombi. The long diagonal has length $2\cos(\alpha/2)$, while the strip construction has length $\ell_Z\sin\alpha$.
\FloatBarrier
\begingroup
\def\LTcaptype{figure}
\makeatletter\def\@captype{figure}\let\caption@stepcounter\@gobble\makeatother
\setlength{\LTpre}{5pt}
\setlength{\LTpost}{5pt}
\setlength{\LTcapwidth}{\linewidth}
\setlength{\tabcolsep}{0pt}
\captionsetup{position=bottom,skip=4pt}
\begin{longtable}{@{}c@{\hspace{0.011000\linewidth}}c@{\hspace{0.011000\linewidth}}c@{}}
\endfirsthead
\endhead
\caption[]{Diamond (rhombus) with escape paths.}\\
\endfoot
\caption{Diamond (rhombus) with escape paths.}\label{fig:rhombi-1}\\
\endlastfoot
\begin{subfigure}[b]{0.326\linewidth}\centering
\includegraphics[width=\linewidth,height=119pt,keepaspectratio]{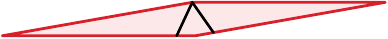}
\par{\footnotesize $\alpha=10^\circ$\\$L=0.383199$\par}
\end{subfigure}
&
\begin{subfigure}[b]{0.326\linewidth}\centering
\includegraphics[width=\linewidth,height=119pt,keepaspectratio]{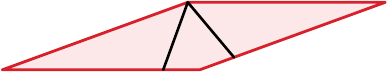}
\par{\footnotesize $\alpha=20^\circ$\\$L=0.727940$\par}
\end{subfigure}
&
\begin{subfigure}[b]{0.326\linewidth}\centering
\includegraphics[width=\linewidth,height=119pt,keepaspectratio]{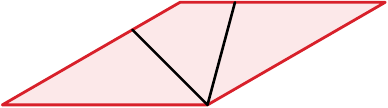}
\par{\footnotesize $\alpha=30^\circ$\\$L=1.035276$\par}
\end{subfigure}
\\
\begin{subfigure}[b]{0.326\linewidth}\centering
\includegraphics[width=\linewidth,height=119pt,keepaspectratio]{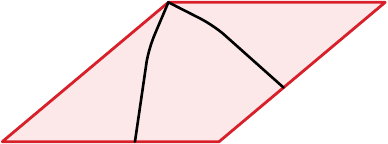}
\par{\footnotesize $\alpha=40^\circ$\\$L=1.330369$\par}
\end{subfigure}
&
\begin{subfigure}[b]{0.326\linewidth}\centering
\includegraphics[width=\linewidth,height=119pt,keepaspectratio]{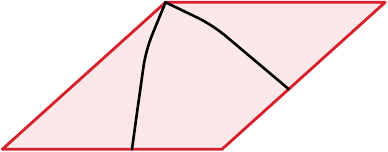}
\par{\footnotesize $\alpha=42^\circ$\\$L=1.380651$\par}
\end{subfigure}
&
\begin{subfigure}[b]{0.326\linewidth}\centering
\includegraphics[width=\linewidth,height=119pt,keepaspectratio]{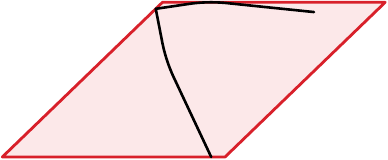}
\par{\footnotesize $\alpha=44^\circ$\\$L=1.429137$\par}
\end{subfigure}
\\
\begin{subfigure}[b]{0.326\linewidth}\centering
\includegraphics[width=\linewidth,height=119pt,keepaspectratio]{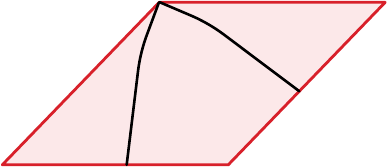}
\par{\footnotesize $\alpha=46^\circ$\\$L=1.475833$\par}
\end{subfigure}
&
\begin{subfigure}[b]{0.326\linewidth}\centering
\includegraphics[width=\linewidth,height=119pt,keepaspectratio]{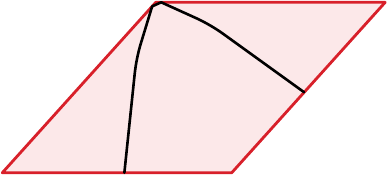}
\par{\footnotesize $\alpha=48^\circ$\\$L=1.520387$\par}
\end{subfigure}
&
\begin{subfigure}[b]{0.326\linewidth}\centering
\includegraphics[width=\linewidth,height=119pt,keepaspectratio]{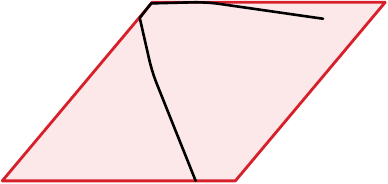}
\par{\footnotesize $\alpha=50^\circ$\\$L=1.562427$\par}
\end{subfigure}
\\
\begin{subfigure}[b]{0.326\linewidth}\centering
\includegraphics[width=\linewidth,height=119pt,keepaspectratio]{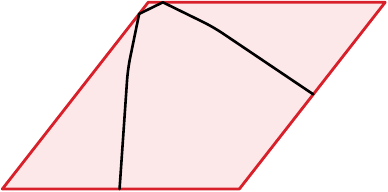}
\par{\footnotesize $\alpha=52^\circ$\\$L=1.601832$\par}
\end{subfigure}
&
\begin{subfigure}[b]{0.326\linewidth}\centering
\includegraphics[width=\linewidth,height=119pt,keepaspectratio]{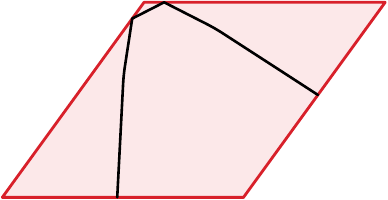}
\par{\footnotesize $\alpha=54^\circ$\\$L=1.638593$\par}
\end{subfigure}
&
\begin{subfigure}[b]{0.326\linewidth}\centering
\includegraphics[width=\linewidth,height=119pt,keepaspectratio]{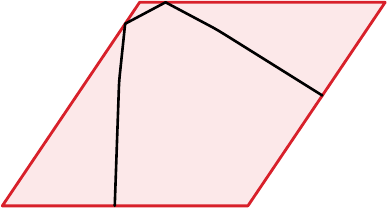}
\par{\footnotesize $\alpha=56^\circ$\\$L=1.672601$\par}
\end{subfigure}
\\
\begin{subfigure}[b]{0.326\linewidth}\centering
\includegraphics[width=\linewidth,height=119pt,keepaspectratio]{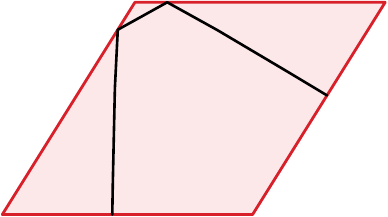}
\par{\footnotesize $\alpha=58^\circ$\\$L=1.703786$\par}
\end{subfigure}
&
\begin{subfigure}[b]{0.326\linewidth}\centering
\includegraphics[width=\linewidth,height=119pt,keepaspectratio]{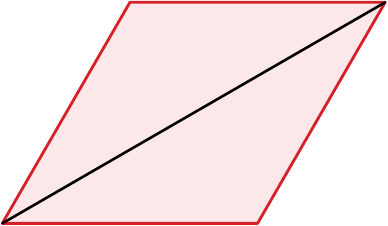}
\par{\footnotesize $\alpha=60^\circ$\\$L=1.732051$\par}
\end{subfigure}
&
\begin{subfigure}[b]{0.326\linewidth}\centering
\includegraphics[width=\linewidth,height=119pt,keepaspectratio]{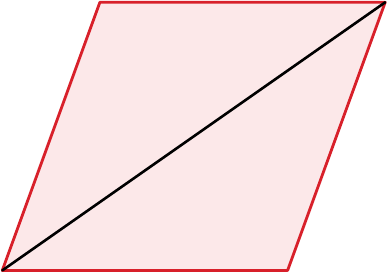}
\par{\footnotesize $\alpha=70^\circ$\\$L=1.638304$\par}
\end{subfigure}
\\
\begin{subfigure}[b]{0.326\linewidth}\centering
\includegraphics[width=\linewidth,height=119pt,keepaspectratio]{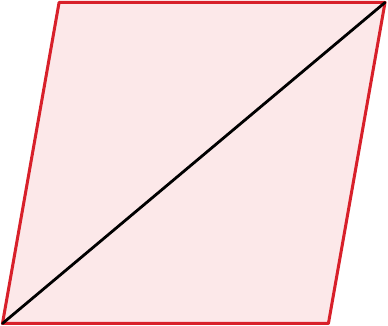}
\par{\footnotesize $\alpha=80^\circ$\\$L=1.532089$\par}
\end{subfigure}
\\
\end{longtable}
\endgroup
\begin{figure}[htbp]
\centering
\includegraphics[width=.95\linewidth]{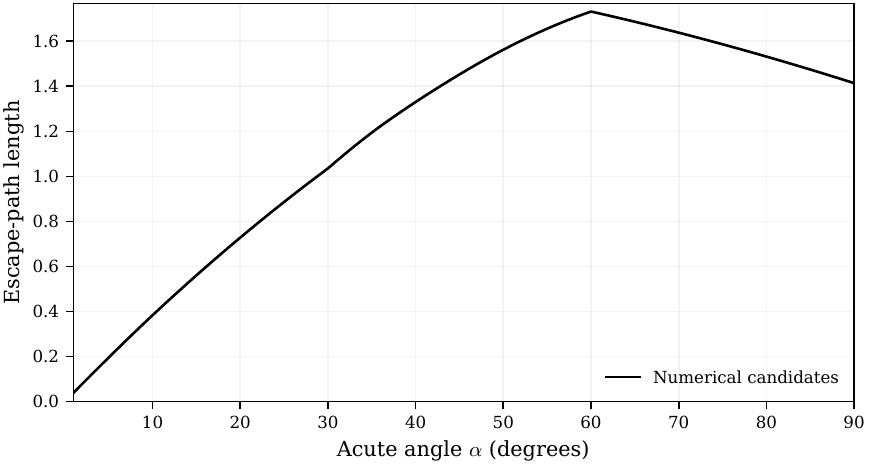}
\caption{Numerical escape path lengths for unit-side rhombi. The bent-path branch persists up to the transition to the long diagonal near $60^\circ$.}
\label{fig:rhombus-relationship}
\end{figure}

The two simple constructions provide a useful comparison, but their minimum does not describe the entire numerical curve. Bent paths are shorter over a substantial range below $60^\circ$. In particular, independent U-shaped starts recover the near-transition bent branch. At $60^\circ$ and above the computation returns the long diagonal.

\subsection{Parallelogram}
\label{sec:parallelograms}

The parallelogram has $p=q=s$ and $\beta=\pi-\alpha$.
All side-length and angle combinations are calculated together on 
$s\in\{0.4,0.6,0.8,1.2,1.4,1.5,1.6,1.8,2\}$ and
$\alpha\in\{10^\circ,20^\circ,\ldots,80^\circ\}$ as shown in Figure 6.
Unit-side rhombi are discussed in
Section~\ref{sec:diamonds} and are omitted from this gallery.
\FloatBarrier
\begingroup
\def\LTcaptype{figure}
\makeatletter\def\@captype{figure}\let\caption@stepcounter\@gobble\makeatother
\setlength{\LTpre}{5pt}
\setlength{\LTpost}{5pt}
\setlength{\LTcapwidth}{\linewidth}
\setlength{\tabcolsep}{0pt}
\captionsetup{position=bottom,skip=4pt}
\begin{longtable}{@{}c@{\hspace{0.009333\linewidth}}c@{\hspace{0.009333\linewidth}}c@{\hspace{0.009333\linewidth}}c@{}}
\endfirsthead
\endhead
\caption[]{Parallelogram paths: 72 cases in parameter order. The unit-side diamond cases are already presented in Fig.~\ref{fig:rhombi-1} and Section~\ref{sec:diamonds}. (continued on following page).}\\
\endfoot
\caption{Parallelogram paths: 72 cases in parameter order. The unit-side diamond cases are already presented in Fig.~\ref{fig:rhombi-1} and Section~\ref{sec:diamonds}.}\label{fig:parallelograms-1}\\
\endlastfoot
\begin{subfigure}[b]{0.243\linewidth}\centering
\includegraphics[width=\linewidth,height=100pt,keepaspectratio]{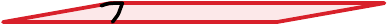}
\par{\footnotesize $s=0.4,\ \alpha=10^\circ$\\$L=0.158248$\par}
\end{subfigure}
&
\begin{subfigure}[b]{0.243\linewidth}\centering
\includegraphics[width=\linewidth,height=100pt,keepaspectratio]{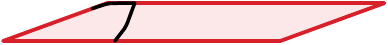}
\par{\footnotesize $s=0.4,\ \alpha=20^\circ$\\$L=0.311689$\par}
\end{subfigure}
&
\begin{subfigure}[b]{0.243\linewidth}\centering
\includegraphics[width=\linewidth,height=100pt,keepaspectratio]{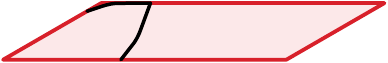}
\par{\footnotesize $s=0.4,\ \alpha=30^\circ$\\$L=0.455658$\par}
\end{subfigure}
&
\begin{subfigure}[b]{0.243\linewidth}\centering
\includegraphics[width=\linewidth,height=100pt,keepaspectratio]{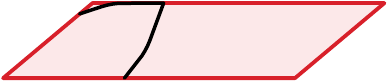}
\par{\footnotesize $s=0.4,\ \alpha=40^\circ$\\$L=0.585783$\par}
\end{subfigure}
\\
\begin{subfigure}[b]{0.243\linewidth}\centering
\includegraphics[width=\linewidth,height=100pt,keepaspectratio]{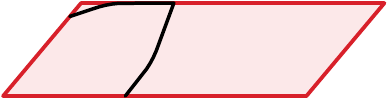}
\par{\footnotesize $s=0.4,\ \alpha=50^\circ$\\$L=0.698109$\par}
\end{subfigure}
&
\begin{subfigure}[b]{0.243\linewidth}\centering
\includegraphics[width=\linewidth,height=100pt,keepaspectratio]{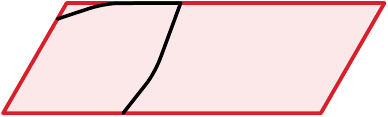}
\par{\footnotesize $s=0.4,\ \alpha=60^\circ$\\$L=0.789223$\par}
\end{subfigure}
&
\begin{subfigure}[b]{0.243\linewidth}\centering
\includegraphics[width=\linewidth,height=100pt,keepaspectratio]{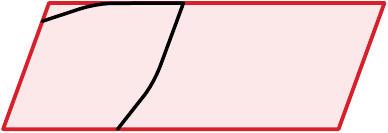}
\par{\footnotesize $s=0.4,\ \alpha=70^\circ$\\$L=0.856358$\par}
\end{subfigure}
&
\begin{subfigure}[b]{0.243\linewidth}\centering
\includegraphics[width=\linewidth,height=100pt,keepaspectratio]{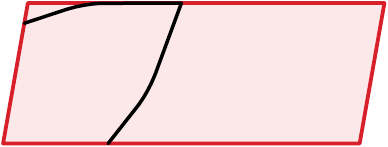}
\par{\footnotesize $s=0.4,\ \alpha=80^\circ$\\$L=0.897472$\par}
\end{subfigure}
\\
\begin{subfigure}[b]{0.243\linewidth}\centering
\includegraphics[width=\linewidth,height=100pt,keepaspectratio]{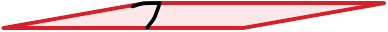}
\par{\footnotesize $s=0.6,\ \alpha=10^\circ$\\$L=0.237373$\par}
\end{subfigure}
&
\begin{subfigure}[b]{0.243\linewidth}\centering
\includegraphics[width=\linewidth,height=100pt,keepaspectratio]{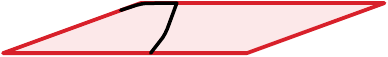}
\par{\footnotesize $s=0.6,\ \alpha=20^\circ$\\$L=0.467533$\par}
\end{subfigure}
&
\begin{subfigure}[b]{0.243\linewidth}\centering
\includegraphics[width=\linewidth,height=100pt,keepaspectratio]{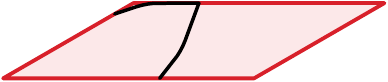}
\par{\footnotesize $s=0.6,\ \alpha=30^\circ$\\$L=0.683488$\par}
\end{subfigure}
&
\begin{subfigure}[b]{0.243\linewidth}\centering
\includegraphics[width=\linewidth,height=100pt,keepaspectratio]{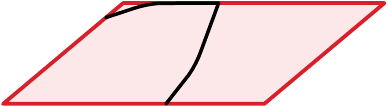}
\par{\footnotesize $s=0.6,\ \alpha=40^\circ$\\$L=0.878675$\par}
\end{subfigure}
\\
\begin{subfigure}[b]{0.243\linewidth}\centering
\includegraphics[width=\linewidth,height=100pt,keepaspectratio]{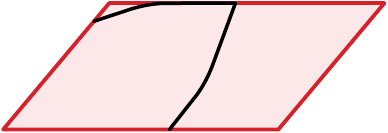}
\par{\footnotesize $s=0.6,\ \alpha=50^\circ$\\$L=1.047164$\par}
\end{subfigure}
&
\begin{subfigure}[b]{0.243\linewidth}\centering
\includegraphics[width=\linewidth,height=100pt,keepaspectratio]{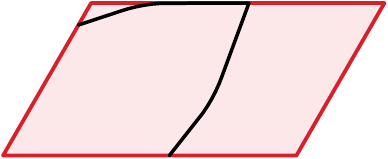}
\par{\footnotesize $s=0.6,\ \alpha=60^\circ$\\$L=1.183835$\par}
\end{subfigure}
&
\begin{subfigure}[b]{0.243\linewidth}\centering
\includegraphics[width=\linewidth,height=100pt,keepaspectratio]{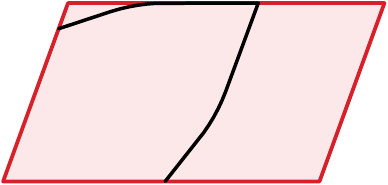}
\par{\footnotesize $s=0.6,\ \alpha=70^\circ$\\$L=1.284536$\par}
\end{subfigure}
&
\begin{subfigure}[b]{0.243\linewidth}\centering
\includegraphics[width=\linewidth,height=100pt,keepaspectratio]{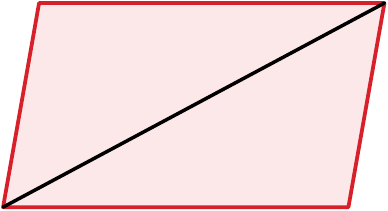}
\par{\footnotesize $s=0.6,\ \alpha=80^\circ$\\$L=1.252349$\par}
\end{subfigure}
\\
\begin{subfigure}[b]{0.243\linewidth}\centering
\includegraphics[width=\linewidth,height=100pt,keepaspectratio]{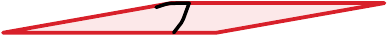}
\par{\footnotesize $s=0.8,\ \alpha=10^\circ$\\$L=0.316497$\par}
\end{subfigure}
&
\begin{subfigure}[b]{0.243\linewidth}\centering
\includegraphics[width=\linewidth,height=100pt,keepaspectratio]{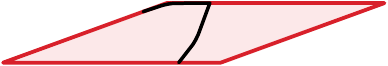}
\par{\footnotesize $s=0.8,\ \alpha=20^\circ$\\$L=0.623377$\par}
\end{subfigure}
&
\begin{subfigure}[b]{0.243\linewidth}\centering
\includegraphics[width=\linewidth,height=100pt,keepaspectratio]{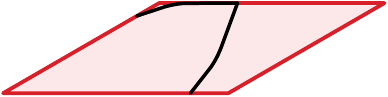}
\par{\footnotesize $s=0.8,\ \alpha=30^\circ$\\$L=0.911317$\par}
\end{subfigure}
&
\begin{subfigure}[b]{0.243\linewidth}\centering
\includegraphics[width=\linewidth,height=100pt,keepaspectratio]{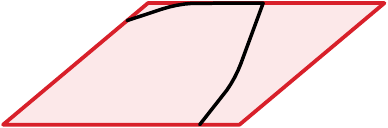}
\par{\footnotesize $s=0.8,\ \alpha=40^\circ$\\$L=1.171566$\par}
\end{subfigure}
\\
\begin{subfigure}[b]{0.243\linewidth}\centering
\includegraphics[width=\linewidth,height=100pt,keepaspectratio]{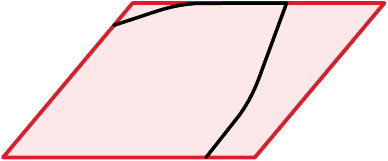}
\par{\footnotesize $s=0.8,\ \alpha=50^\circ$\\$L=1.396218$\par}
\end{subfigure}
&
\begin{subfigure}[b]{0.243\linewidth}\centering
\includegraphics[width=\linewidth,height=100pt,keepaspectratio]{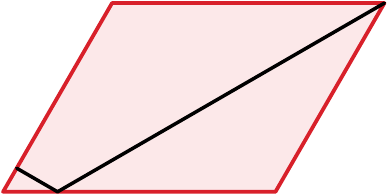}
\par{\footnotesize $s=0.8,\ \alpha=60^\circ$\\$L=1.558846$\par}
\end{subfigure}
&
\begin{subfigure}[b]{0.243\linewidth}\centering
\includegraphics[width=\linewidth,height=100pt,keepaspectratio]{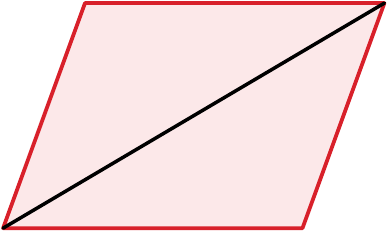}
\par{\footnotesize $s=0.8,\ \alpha=70^\circ$\\$L=1.478929$\par}
\end{subfigure}
&
\begin{subfigure}[b]{0.243\linewidth}\centering
\includegraphics[width=\linewidth,height=100pt,keepaspectratio]{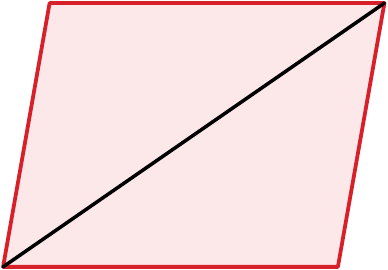}
\par{\footnotesize $s=0.8,\ \alpha=80^\circ$\\$L=1.384860$\par}
\end{subfigure}
\\
\begin{subfigure}[b]{0.243\linewidth}\centering
\includegraphics[width=\linewidth,height=100pt,keepaspectratio]{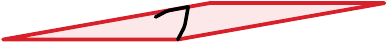}
\par{\footnotesize $s=1.2,\ \alpha=10^\circ$\\$L=0.395621$\par}
\end{subfigure}
&
\begin{subfigure}[b]{0.243\linewidth}\centering
\includegraphics[width=\linewidth,height=100pt,keepaspectratio]{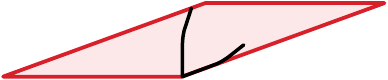}
\par{\footnotesize $s=1.2,\ \alpha=20^\circ$\\$L=0.779222$\par}
\end{subfigure}
&
\begin{subfigure}[b]{0.243\linewidth}\centering
\includegraphics[width=\linewidth,height=100pt,keepaspectratio]{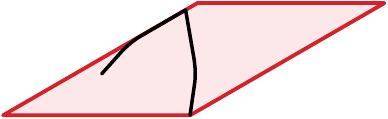}
\par{\footnotesize $s=1.2,\ \alpha=30^\circ$\\$L=1.139146$\par}
\end{subfigure}
&
\begin{subfigure}[b]{0.243\linewidth}\centering
\includegraphics[width=\linewidth,height=100pt,keepaspectratio]{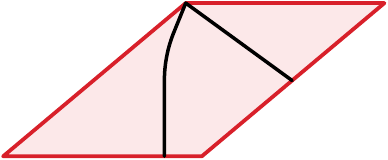}
\par{\footnotesize $s=1.2,\ \alpha=40^\circ$\\$L=1.451510$\par}
\end{subfigure}
\\
\begin{subfigure}[b]{0.243\linewidth}\centering
\includegraphics[width=\linewidth,height=100pt,keepaspectratio]{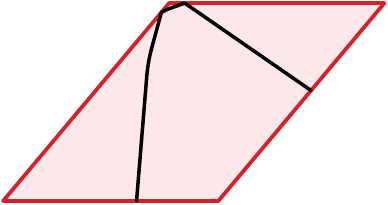}
\par{\footnotesize $s=1.2,\ \alpha=50^\circ$\\$L=1.716782$\par}
\end{subfigure}
&
\begin{subfigure}[b]{0.243\linewidth}\centering
\includegraphics[width=\linewidth,height=100pt,keepaspectratio]{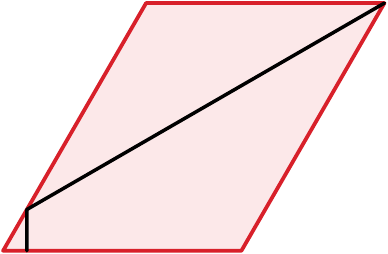}
\par{\footnotesize $s=1.2,\ \alpha=60^\circ$\\$L=1.905256$\par}
\end{subfigure}
&
\begin{subfigure}[b]{0.243\linewidth}\centering
\includegraphics[width=\linewidth,height=100pt,keepaspectratio]{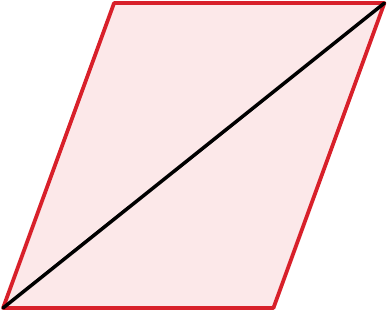}
\par{\footnotesize $s=1.2,\ \alpha=70^\circ$\\$L=1.805782$\par}
\end{subfigure}
&
\begin{subfigure}[b]{0.243\linewidth}\centering
\includegraphics[width=\linewidth,height=100pt,keepaspectratio]{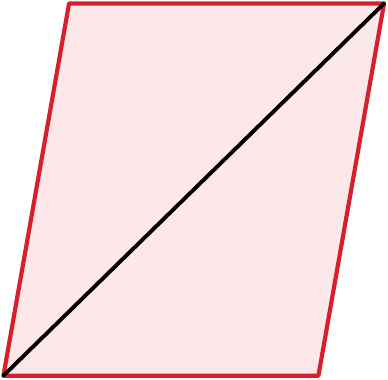}
\par{\footnotesize $s=1.2,\ \alpha=80^\circ$\\$L=1.690194$\par}
\end{subfigure}
\\
\begin{subfigure}[b]{0.243\linewidth}\centering
\includegraphics[width=\linewidth,height=100pt,keepaspectratio]{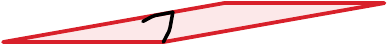}
\par{\footnotesize $s=1.4,\ \alpha=10^\circ$\\$L=0.395621$\par}
\end{subfigure}
&
\begin{subfigure}[b]{0.243\linewidth}\centering
\includegraphics[width=\linewidth,height=100pt,keepaspectratio]{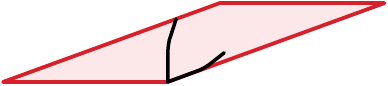}
\par{\footnotesize $s=1.4,\ \alpha=20^\circ$\\$L=0.779222$\par}
\end{subfigure}
&
\begin{subfigure}[b]{0.243\linewidth}\centering
\includegraphics[width=\linewidth,height=100pt,keepaspectratio]{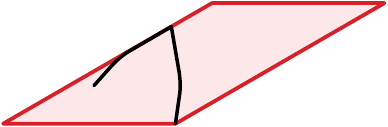}
\par{\footnotesize $s=1.4,\ \alpha=30^\circ$\\$L=1.139146$\par}
\end{subfigure}
&
\begin{subfigure}[b]{0.243\linewidth}\centering
\includegraphics[width=\linewidth,height=100pt,keepaspectratio]{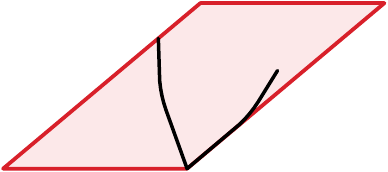}
\par{\footnotesize $s=1.4,\ \alpha=40^\circ$\\$L=1.464458$\par}
\end{subfigure}
\\
\begin{subfigure}[b]{0.243\linewidth}\centering
\includegraphics[width=\linewidth,height=100pt,keepaspectratio]{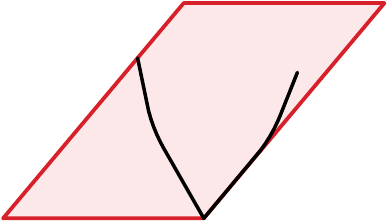}
\par{\footnotesize $s=1.4,\ \alpha=50^\circ$\\$L=1.745273$\par}
\end{subfigure}
&
\begin{subfigure}[b]{0.243\linewidth}\centering
\includegraphics[width=\linewidth,height=100pt,keepaspectratio]{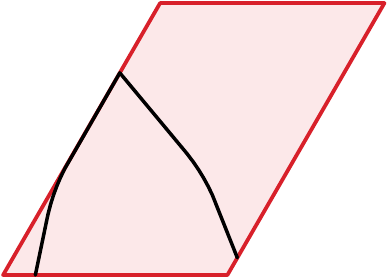}
\par{\footnotesize $s=1.4,\ \alpha=60^\circ$\\$L=1.973059$\par}
\end{subfigure}
&
\begin{subfigure}[b]{0.243\linewidth}\centering
\includegraphics[width=\linewidth,height=100pt,keepaspectratio]{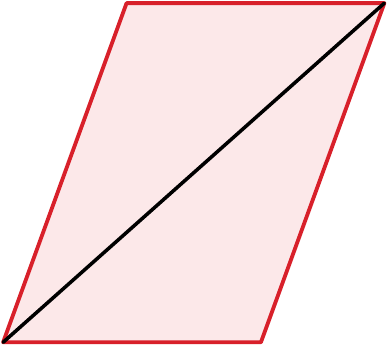}
\par{\footnotesize $s=1.4,\ \alpha=70^\circ$\\$L=1.979307$\par}
\end{subfigure}
&
\begin{subfigure}[b]{0.243\linewidth}\centering
\includegraphics[width=\linewidth,height=100pt,keepaspectratio]{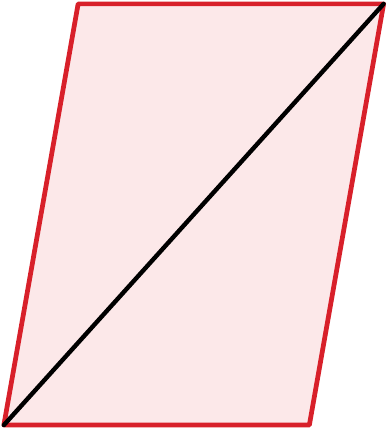}
\par{\footnotesize $s=1.4,\ \alpha=80^\circ$\\$L=1.856398$\par}
\end{subfigure}
\\
\begin{subfigure}[b]{0.243\linewidth}\centering
\includegraphics[width=\linewidth,height=100pt,keepaspectratio]{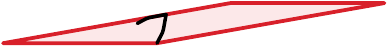}
\par{\footnotesize $s=1.5,\ \alpha=10^\circ$\\$L=0.395621$\par}
\end{subfigure}
&
\begin{subfigure}[b]{0.243\linewidth}\centering
\includegraphics[width=\linewidth,height=100pt,keepaspectratio]{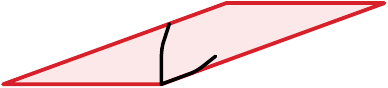}
\par{\footnotesize $s=1.5,\ \alpha=20^\circ$\\$L=0.779222$\par}
\end{subfigure}
&
\begin{subfigure}[b]{0.243\linewidth}\centering
\includegraphics[width=\linewidth,height=100pt,keepaspectratio]{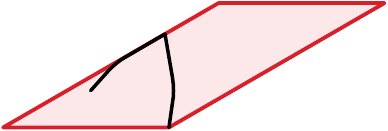}
\par{\footnotesize $s=1.5,\ \alpha=30^\circ$\\$L=1.139146$\par}
\end{subfigure}
&
\begin{subfigure}[b]{0.243\linewidth}\centering
\includegraphics[width=\linewidth,height=100pt,keepaspectratio]{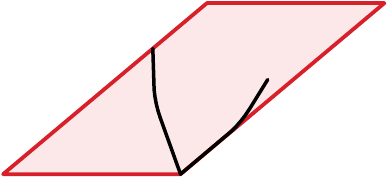}
\par{\footnotesize $s=1.5,\ \alpha=40^\circ$\\$L=1.464458$\par}
\end{subfigure}
\\
\begin{subfigure}[b]{0.243\linewidth}\centering
\includegraphics[width=\linewidth,height=100pt,keepaspectratio]{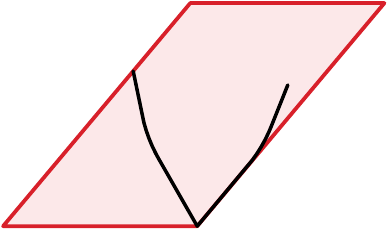}
\par{\footnotesize $s=1.5,\ \alpha=50^\circ$\\$L=1.745273$\par}
\end{subfigure}
&
\begin{subfigure}[b]{0.243\linewidth}\centering
\includegraphics[width=\linewidth,height=100pt,keepaspectratio]{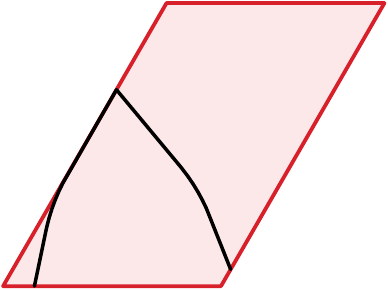}
\par{\footnotesize $s=1.5,\ \alpha=60^\circ$\\$L=1.973059$\par}
\end{subfigure}
&
\begin{subfigure}[b]{0.243\linewidth}\centering
\includegraphics[width=\linewidth,height=100pt,keepaspectratio]{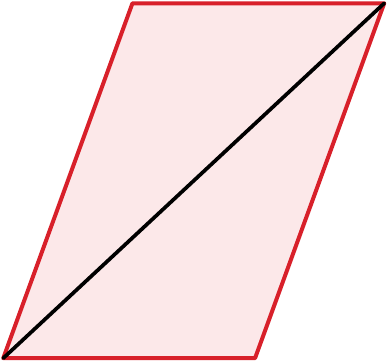}
\par{\footnotesize $s=1.5,\ \alpha=70^\circ$\\$L=2.067864$\par}
\end{subfigure}
&
\begin{subfigure}[b]{0.243\linewidth}\centering
\includegraphics[width=\linewidth,height=100pt,keepaspectratio]{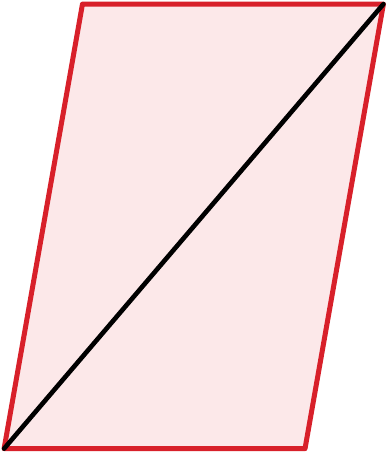}
\par{\footnotesize $s=1.5,\ \alpha=80^\circ$\\$L=1.941892$\par}
\end{subfigure}
\\
\begin{subfigure}[b]{0.243\linewidth}\centering
\includegraphics[width=\linewidth,height=100pt,keepaspectratio]{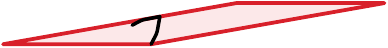}
\par{\footnotesize $s=1.6,\ \alpha=10^\circ$\\$L=0.395621$\par}
\end{subfigure}
&
\begin{subfigure}[b]{0.243\linewidth}\centering
\includegraphics[width=\linewidth,height=100pt,keepaspectratio]{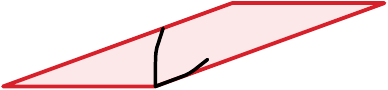}
\par{\footnotesize $s=1.6,\ \alpha=20^\circ$\\$L=0.779222$\par}
\end{subfigure}
&
\begin{subfigure}[b]{0.243\linewidth}\centering
\includegraphics[width=\linewidth,height=100pt,keepaspectratio]{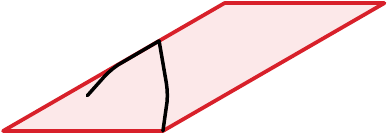}
\par{\footnotesize $s=1.6,\ \alpha=30^\circ$\\$L=1.139146$\par}
\end{subfigure}
&
\begin{subfigure}[b]{0.243\linewidth}\centering
\includegraphics[width=\linewidth,height=100pt,keepaspectratio]{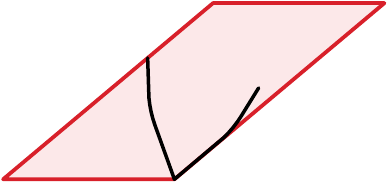}
\par{\footnotesize $s=1.6,\ \alpha=40^\circ$\\$L=1.464458$\par}
\end{subfigure}
\\
\begin{subfigure}[b]{0.243\linewidth}\centering
\includegraphics[width=\linewidth,height=100pt,keepaspectratio]{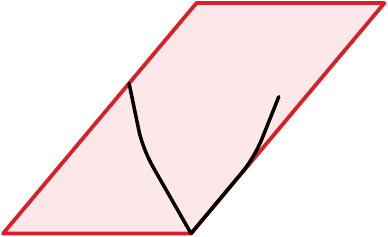}
\par{\footnotesize $s=1.6,\ \alpha=50^\circ$\\$L=1.745273$\par}
\end{subfigure}
&
\begin{subfigure}[b]{0.243\linewidth}\centering
\includegraphics[width=\linewidth,height=100pt,keepaspectratio]{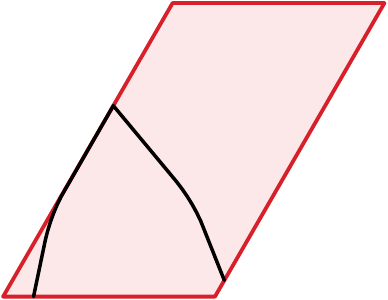}
\par{\footnotesize $s=1.6,\ \alpha=60^\circ$\\$L=1.973059$\par}
\end{subfigure}
&
\begin{subfigure}[b]{0.243\linewidth}\centering
\includegraphics[width=\linewidth,height=100pt,keepaspectratio]{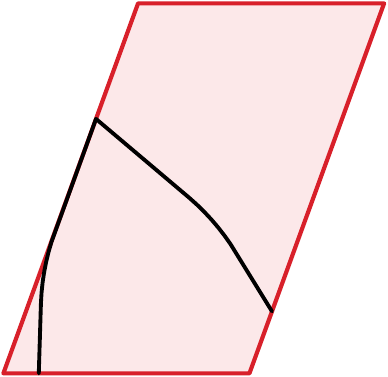}
\par{\footnotesize $s=1.6,\ \alpha=70^\circ$\\$L=2.140894$\par}
\end{subfigure}
&
\begin{subfigure}[b]{0.243\linewidth}\centering
\includegraphics[width=\linewidth,height=100pt,keepaspectratio]{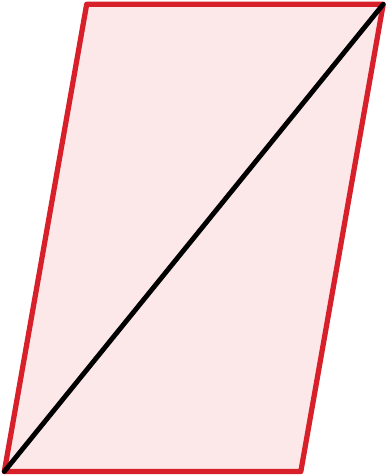}
\par{\footnotesize $s=1.6,\ \alpha=80^\circ$\\$L=2.028712$\par}
\end{subfigure}
\\
\begin{subfigure}[b]{0.243\linewidth}\centering
\includegraphics[width=\linewidth,height=100pt,keepaspectratio]{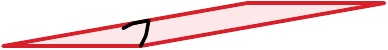}
\par{\footnotesize $s=1.8,\ \alpha=10^\circ$\\$L=0.395621$\par}
\end{subfigure}
&
\begin{subfigure}[b]{0.243\linewidth}\centering
\includegraphics[width=\linewidth,height=100pt,keepaspectratio]{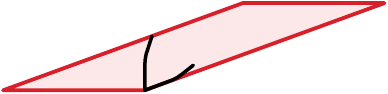}
\par{\footnotesize $s=1.8,\ \alpha=20^\circ$\\$L=0.779222$\par}
\end{subfigure}
&
\begin{subfigure}[b]{0.243\linewidth}\centering
\includegraphics[width=\linewidth,height=100pt,keepaspectratio]{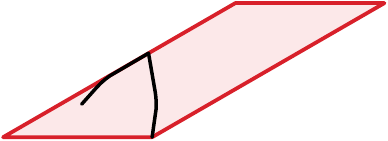}
\par{\footnotesize $s=1.8,\ \alpha=30^\circ$\\$L=1.139146$\par}
\end{subfigure}
&
\begin{subfigure}[b]{0.243\linewidth}\centering
\includegraphics[width=\linewidth,height=100pt,keepaspectratio]{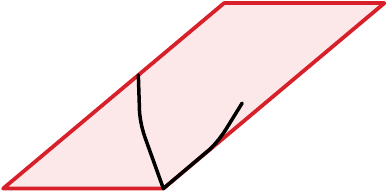}
\par{\footnotesize $s=1.8,\ \alpha=40^\circ$\\$L=1.464458$\par}
\end{subfigure}
\\
\begin{subfigure}[b]{0.243\linewidth}\centering
\includegraphics[width=\linewidth,height=100pt,keepaspectratio]{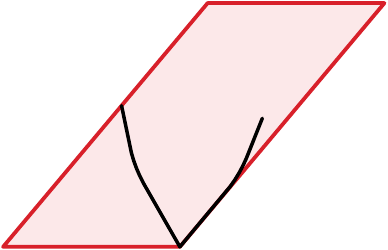}
\par{\footnotesize $s=1.8,\ \alpha=50^\circ$\\$L=1.745273$\par}
\end{subfigure}
&
\begin{subfigure}[b]{0.243\linewidth}\centering
\includegraphics[width=\linewidth,height=100pt,keepaspectratio]{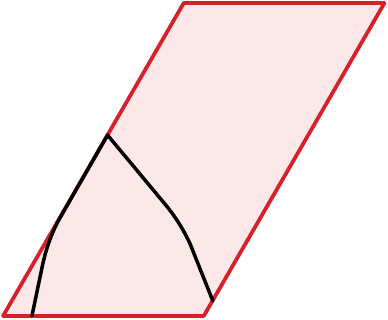}
\par{\footnotesize $s=1.8,\ \alpha=60^\circ$\\$L=1.973059$\par}
\end{subfigure}
&
\begin{subfigure}[b]{0.243\linewidth}\centering
\includegraphics[width=\linewidth,height=100pt,keepaspectratio]{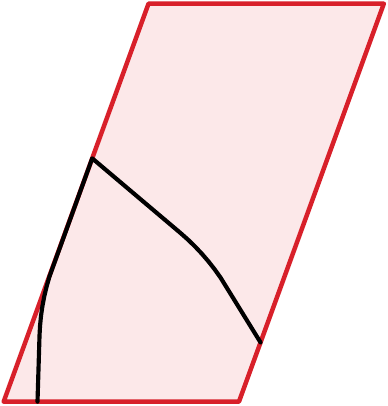}
\par{\footnotesize $s=1.8,\ \alpha=70^\circ$\\$L=2.140894$\par}
\end{subfigure}
&
\begin{subfigure}[b]{0.243\linewidth}\centering
\includegraphics[width=\linewidth,height=100pt,keepaspectratio]{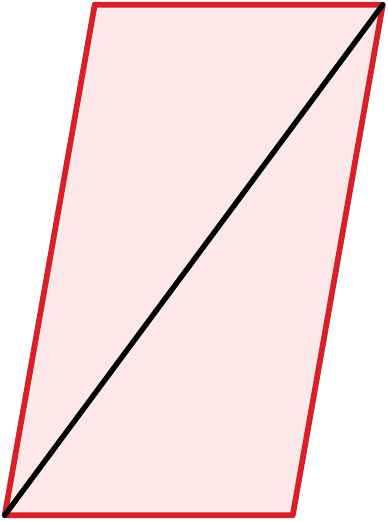}
\par{\footnotesize $s=1.8,\ \alpha=80^\circ$\\$L=2.205705$\par}
\end{subfigure}
\\
\begin{subfigure}[b]{0.243\linewidth}\centering
\includegraphics[width=\linewidth,height=100pt,keepaspectratio]{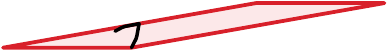}
\par{\footnotesize $s=2,\ \alpha=10^\circ$\\$L=0.395621$\par}
\end{subfigure}
&
\begin{subfigure}[b]{0.243\linewidth}\centering
\includegraphics[width=\linewidth,height=100pt,keepaspectratio]{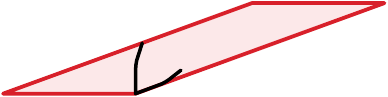}
\par{\footnotesize $s=2,\ \alpha=20^\circ$\\$L=0.779222$\par}
\end{subfigure}
&
\begin{subfigure}[b]{0.243\linewidth}\centering
\includegraphics[width=\linewidth,height=100pt,keepaspectratio]{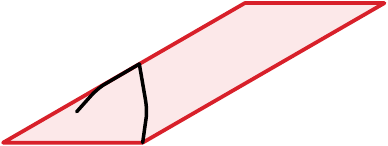}
\par{\footnotesize $s=2,\ \alpha=30^\circ$\\$L=1.139146$\par}
\end{subfigure}
&
\begin{subfigure}[b]{0.243\linewidth}\centering
\includegraphics[width=\linewidth,height=100pt,keepaspectratio]{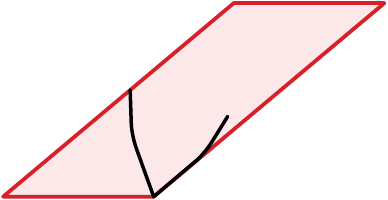}
\par{\footnotesize $s=2,\ \alpha=40^\circ$\\$L=1.464458$\par}
\end{subfigure}
\\
\begin{subfigure}[b]{0.243\linewidth}\centering
\includegraphics[width=\linewidth,height=100pt,keepaspectratio]{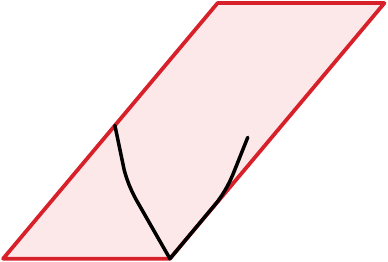}
\par{\footnotesize $s=2,\ \alpha=50^\circ$\\$L=1.745273$\par}
\end{subfigure}
&
\begin{subfigure}[b]{0.243\linewidth}\centering
\includegraphics[width=\linewidth,height=100pt,keepaspectratio]{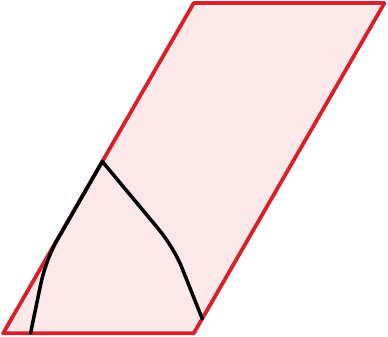}
\par{\footnotesize $s=2,\ \alpha=60^\circ$\\$L=1.973059$\par}
\end{subfigure}
&
\begin{subfigure}[b]{0.243\linewidth}\centering
\includegraphics[width=\linewidth,height=100pt,keepaspectratio]{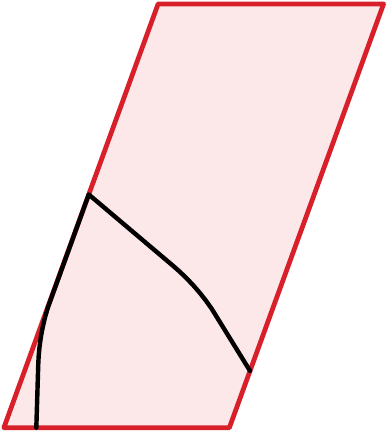}
\par{\footnotesize $s=2,\ \alpha=70^\circ$\\$L=2.140894$\par}
\end{subfigure}
&
\begin{subfigure}[b]{0.243\linewidth}\centering
\includegraphics[width=\linewidth,height=100pt,keepaspectratio]{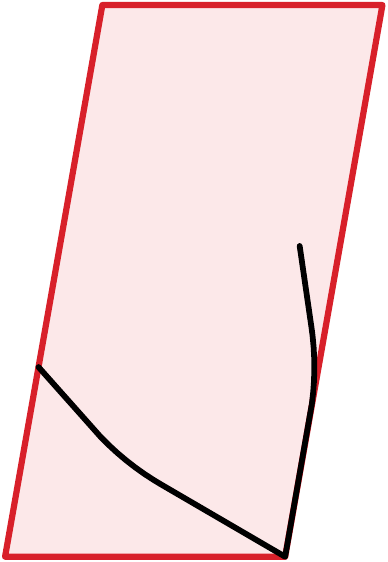}
\par{\footnotesize $s=2,\ \alpha=80^\circ$\\$L=2.243679$\par}
\end{subfigure}
\\
\end{longtable}
\endgroup
\begin{figure}[htbp]
\centering
\includegraphics[width=.83\linewidth]{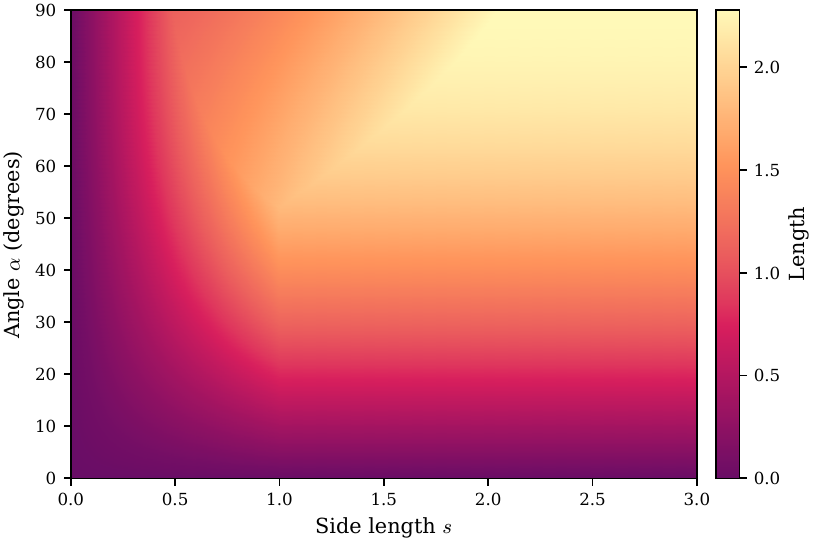}
\caption{Parallelogram comparison relationship from
Eq.~\eqref{eq:comparison-envelope}, with the original draft's color scheme.
Color indicates the comparison length and fills continuously across the branch boundaries.}
\label{fig:parallelogram-density}
\end{figure}

Figure 7 shows the summary of the parallelogram case.The summary uses side length on the horizontal axis, angle on the
vertical axis, and path length by color. The minimum of the strip and
diagonal constructions is retained as the explicit comparison
\begin{equation}\label{eq:comparison-envelope}
L_{\rm env}(s,\alpha)=\min\!\left\{
 \ell_Z\min(1,s)\sin\alpha,\sqrt{1+s^2+2s\cos\alpha}\right\}.
\end{equation}
The two terms are independently valid constructions. The numerical search
also finds bent paths below this envelope. The density plot therefore
represents the explicit comparison in Eq.~\eqref{eq:comparison-envelope}.

\subsection{Isosceles Trapezoid}
\label{sec:trapezoids}

For this family $p=q=s$ and $\alpha=\beta$. Its height is
$h=s\sin\alpha$ and upper base is
\begin{equation}
b_{\rm top}=1-2s\cos\alpha.
\end{equation}
A proper convex trapezoid requires $b_{\rm top}>0$.
The gallery in Figure 8 retains all 32 parameter slots of the parameter grid, at
$s\in\{0.2,0.4,0.6,0.8,1,1.5,2,2.5\}$ and
$\alpha\in\{20^\circ,40^\circ,60^\circ,80^\circ\}$.
\FloatBarrier
\begingroup
\def\LTcaptype{figure}
\makeatletter\def\@captype{figure}\let\caption@stepcounter\@gobble\makeatother
\setlength{\LTpre}{5pt}
\setlength{\LTpost}{5pt}
\setlength{\LTcapwidth}{\linewidth}
\setlength{\tabcolsep}{0pt}
\captionsetup{position=bottom,skip=4pt}
\begin{longtable}{@{}c@{\hspace{0.009333\linewidth}}c@{\hspace{0.009333\linewidth}}c@{\hspace{0.009333\linewidth}}c@{}}
\endfirsthead
\endhead
\caption[]{Isosceles trapezoid paths and excluded slots: 32 cases in parameter order. (continued on following page).}\\
\endfoot
\caption{Isosceles trapezoid paths and excluded slots: 32 cases in parameter order.}\label{fig:trapezoids-1}\\
\endlastfoot
\begin{subfigure}[b]{0.243\linewidth}\centering
\includegraphics[width=\linewidth,height=100pt,keepaspectratio]{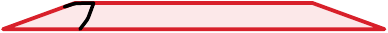}
\par{\footnotesize $s=0.2,\ \alpha=20^\circ$\\$L=0.155844$\par}
\end{subfigure}
&
\begin{subfigure}[b]{0.243\linewidth}\centering
\includegraphics[width=\linewidth,height=100pt,keepaspectratio]{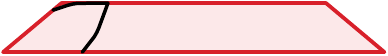}
\par{\footnotesize $s=0.2,\ \alpha=40^\circ$\\$L=0.292892$\par}
\end{subfigure}
&
\begin{subfigure}[b]{0.243\linewidth}\centering
\includegraphics[width=\linewidth,height=100pt,keepaspectratio]{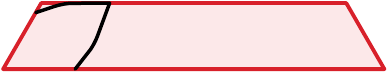}
\par{\footnotesize $s=0.2,\ \alpha=60^\circ$\\$L=0.394612$\par}
\end{subfigure}
&
\begin{subfigure}[b]{0.243\linewidth}\centering
\includegraphics[width=\linewidth,height=100pt,keepaspectratio]{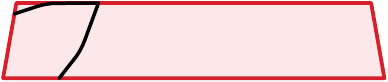}
\par{\footnotesize $s=0.2,\ \alpha=80^\circ$\\$L=0.448736$\par}
\end{subfigure}
\\
\begin{subfigure}[b]{0.243\linewidth}\centering
\includegraphics[width=\linewidth,height=100pt,keepaspectratio]{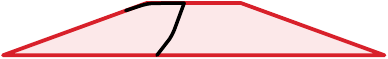}
\par{\footnotesize $s=0.4,\ \alpha=20^\circ$\\$L=0.311689$\par}
\end{subfigure}
&
\begin{subfigure}[b]{0.243\linewidth}\centering
\includegraphics[width=\linewidth,height=100pt,keepaspectratio]{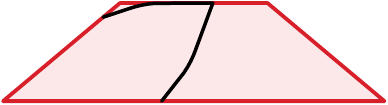}
\par{\footnotesize $s=0.4,\ \alpha=40^\circ$\\$L=0.585783$\par}
\end{subfigure}
&
\begin{subfigure}[b]{0.243\linewidth}\centering
\includegraphics[width=\linewidth,height=100pt,keepaspectratio]{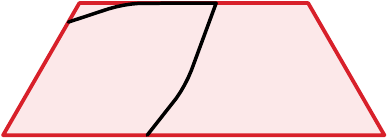}
\par{\footnotesize $s=0.4,\ \alpha=60^\circ$\\$L=0.789223$\par}
\end{subfigure}
&
\begin{subfigure}[b]{0.243\linewidth}\centering
\includegraphics[width=\linewidth,height=100pt,keepaspectratio]{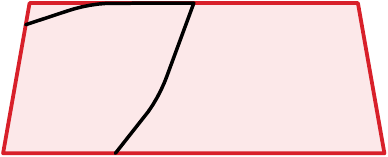}
\par{\footnotesize $s=0.4,\ \alpha=80^\circ$\\$L=0.897472$\par}
\end{subfigure}
\\
\begin{subfigure}[b]{0.243\linewidth}\centering
\includegraphics[width=\linewidth,height=100pt,keepaspectratio]{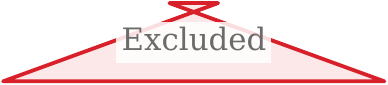}
\par{\footnotesize $s=0.6,\ \alpha=20^\circ$\par}
\end{subfigure}
&
\begin{subfigure}[b]{0.243\linewidth}\centering
\includegraphics[width=\linewidth,height=100pt,keepaspectratio]{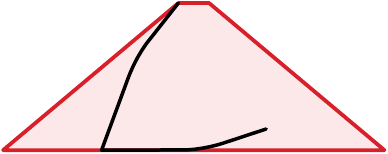}
\par{\footnotesize $s=0.6,\ \alpha=40^\circ$\\$L=0.878675$\par}
\end{subfigure}
&
\begin{subfigure}[b]{0.243\linewidth}\centering
\includegraphics[width=\linewidth,height=100pt,keepaspectratio]{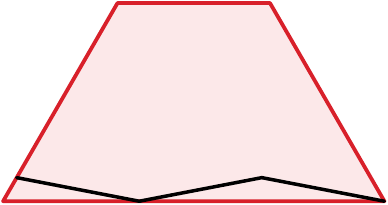}
\par{\footnotesize $s=0.6,\ \alpha=60^\circ$\\$L=0.981981$\par}
\end{subfigure}
&
\begin{subfigure}[b]{0.243\linewidth}\centering
\includegraphics[width=\linewidth,height=100pt,keepaspectratio]{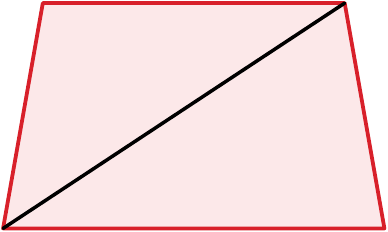}
\par{\footnotesize $s=0.6,\ \alpha=80^\circ$\\$L=1.073137$\par}
\end{subfigure}
\\
\begin{subfigure}[b]{0.243\linewidth}\centering
\includegraphics[width=\linewidth,height=100pt,keepaspectratio]{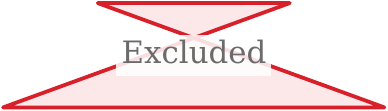}
\par{\footnotesize $s=0.8,\ \alpha=20^\circ$\par}
\end{subfigure}
&
\begin{subfigure}[b]{0.243\linewidth}\centering
\includegraphics[width=\linewidth,height=100pt,keepaspectratio]{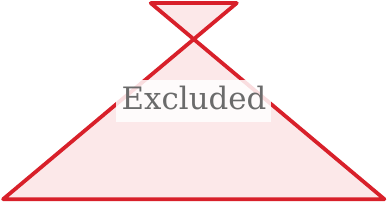}
\par{\footnotesize $s=0.8,\ \alpha=40^\circ$\par}
\end{subfigure}
&
\begin{subfigure}[b]{0.243\linewidth}\centering
\includegraphics[width=\linewidth,height=100pt,keepaspectratio]{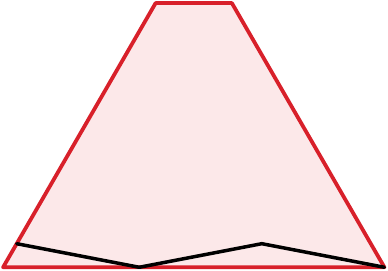}
\par{\footnotesize $s=0.8,\ \alpha=60^\circ$\\$L=0.981981$\par}
\end{subfigure}
&
\begin{subfigure}[b]{0.243\linewidth}\centering
\includegraphics[width=\linewidth,height=100pt,keepaspectratio]{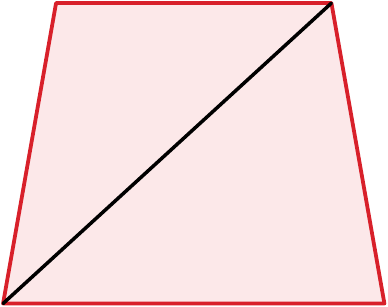}
\par{\footnotesize $s=0.8,\ \alpha=80^\circ$\\$L=1.167117$\par}
\end{subfigure}
\\
\begin{subfigure}[b]{0.243\linewidth}\centering
\includegraphics[width=\linewidth,height=100pt,keepaspectratio]{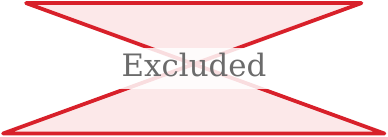}
\par{\footnotesize $s=1,\ \alpha=20^\circ$\par}
\end{subfigure}
&
\begin{subfigure}[b]{0.243\linewidth}\centering
\includegraphics[width=\linewidth,height=100pt,keepaspectratio]{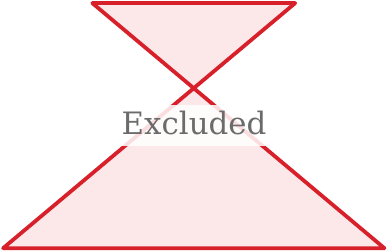}
\par{\footnotesize $s=1,\ \alpha=40^\circ$\par}
\end{subfigure}
&
\begin{subfigure}[b]{0.243\linewidth}\centering
\includegraphics[width=\linewidth,height=100pt,keepaspectratio]{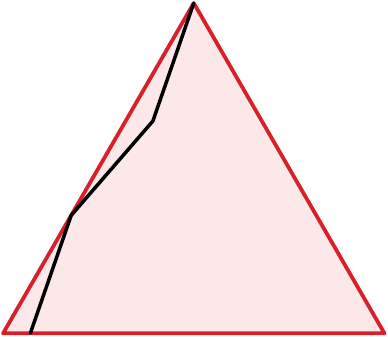}
\par{\footnotesize $s=1,\ \alpha=60^\circ$\\$L=0.981981$\\triangle limit\par}
\end{subfigure}
&
\begin{subfigure}[b]{0.243\linewidth}\centering
\includegraphics[width=\linewidth,height=100pt,keepaspectratio]{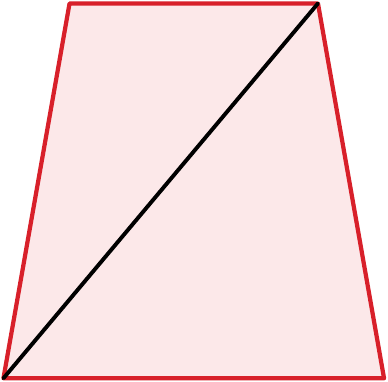}
\par{\footnotesize $s=1,\ \alpha=80^\circ$\\$L=1.285575$\par}
\end{subfigure}
\\
\begin{subfigure}[b]{0.243\linewidth}\centering
\includegraphics[width=\linewidth,height=100pt,keepaspectratio]{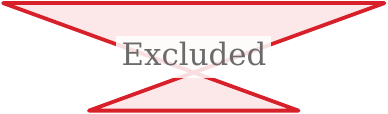}
\par{\footnotesize $s=1.5,\ \alpha=20^\circ$\par}
\end{subfigure}
&
\begin{subfigure}[b]{0.243\linewidth}\centering
\includegraphics[width=\linewidth,height=100pt,keepaspectratio]{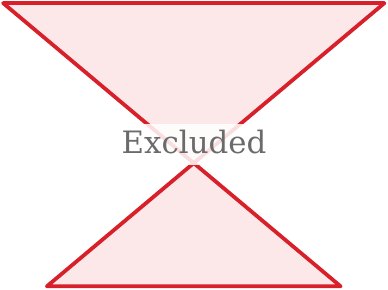}
\par{\footnotesize $s=1.5,\ \alpha=40^\circ$\par}
\end{subfigure}
&
\begin{subfigure}[b]{0.243\linewidth}\centering
\includegraphics[width=\linewidth,height=100pt,keepaspectratio]{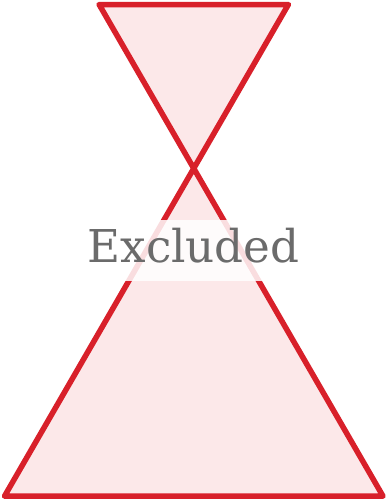}
\par{\footnotesize $s=1.5,\ \alpha=60^\circ$\par}
\end{subfigure}
&
\begin{subfigure}[b]{0.243\linewidth}\centering
\includegraphics[width=\linewidth,height=100pt,keepaspectratio]{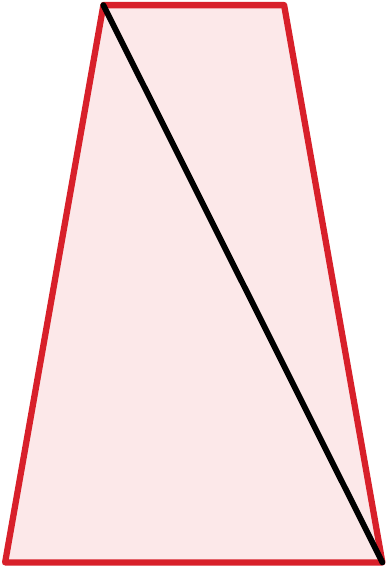}
\par{\footnotesize $s=1.5,\ \alpha=80^\circ$\\$L=1.651985$\par}
\end{subfigure}
\\
\begin{subfigure}[b]{0.243\linewidth}\centering
\includegraphics[width=\linewidth,height=100pt,keepaspectratio]{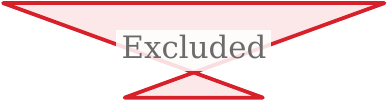}
\par{\footnotesize $s=2,\ \alpha=20^\circ$\par}
\end{subfigure}
&
\begin{subfigure}[b]{0.243\linewidth}\centering
\includegraphics[width=\linewidth,height=100pt,keepaspectratio]{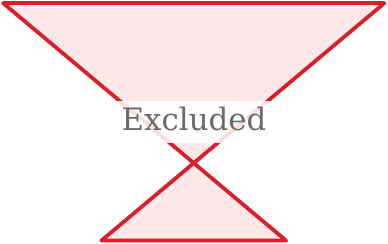}
\par{\footnotesize $s=2,\ \alpha=40^\circ$\par}
\end{subfigure}
&
\begin{subfigure}[b]{0.243\linewidth}\centering
\includegraphics[width=\linewidth,height=100pt,keepaspectratio]{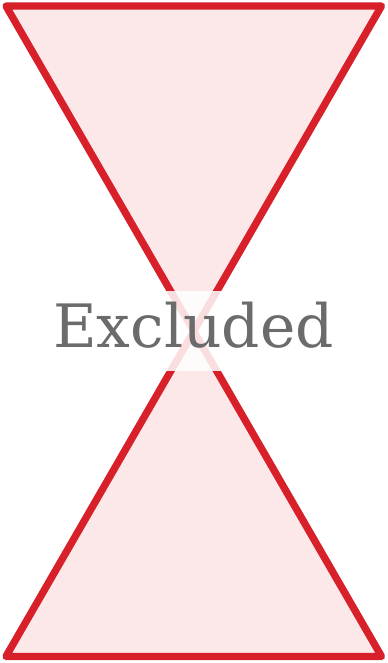}
\par{\footnotesize $s=2,\ \alpha=60^\circ$\par}
\end{subfigure}
&
\begin{subfigure}[b]{0.243\linewidth}\centering
\includegraphics[width=\linewidth,height=100pt,keepaspectratio]{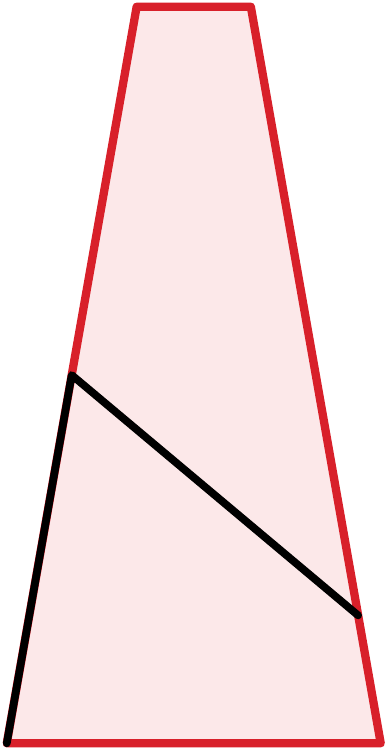}
\par{\footnotesize $s=2,\ \alpha=80^\circ$\\$L=2.000000$\par}
\end{subfigure}
\\
\begin{subfigure}[b]{0.243\linewidth}\centering
\includegraphics[width=\linewidth,height=100pt,keepaspectratio]{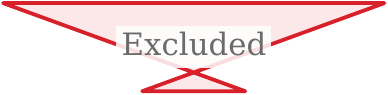}
\par{\footnotesize $s=2.5,\ \alpha=20^\circ$\par}
\end{subfigure}
&
\begin{subfigure}[b]{0.243\linewidth}\centering
\includegraphics[width=\linewidth,height=100pt,keepaspectratio]{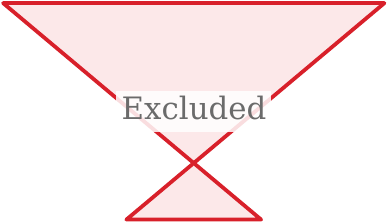}
\par{\footnotesize $s=2.5,\ \alpha=40^\circ$\par}
\end{subfigure}
&
\begin{subfigure}[b]{0.243\linewidth}\centering
\includegraphics[width=\linewidth,height=100pt,keepaspectratio]{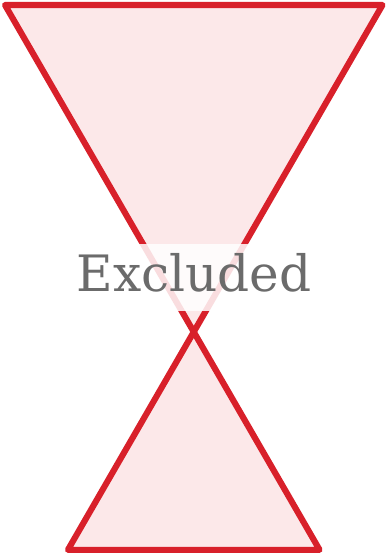}
\par{\footnotesize $s=2.5,\ \alpha=60^\circ$\par}
\end{subfigure}
&
\begin{subfigure}[b]{0.243\linewidth}\centering
\includegraphics[width=\linewidth,height=100pt,keepaspectratio]{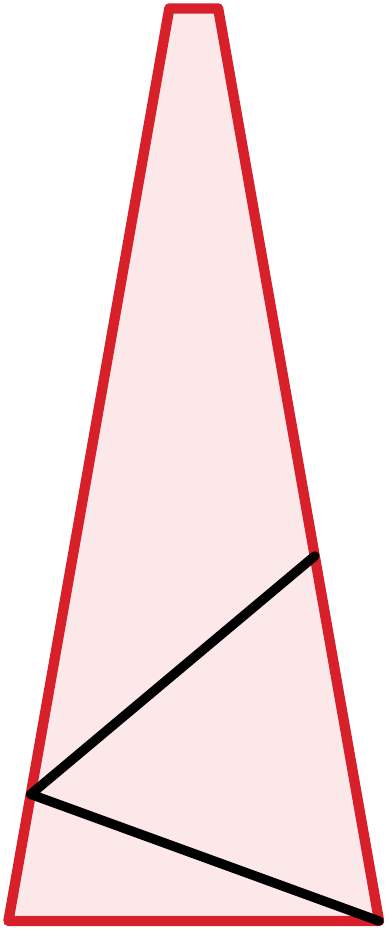}
\par{\footnotesize $s=2.5,\ \alpha=80^\circ$\\$L=2.000000$\par}
\end{subfigure}
\\
\end{longtable}
\endgroup

For the summary of the isosceles trapezoid case, the three recurring path categories are described below. 
\paragraph{Parallel lines.} When the two bases control fitting, a strip-shaped path provides the comparison $\ell_Z h$. The computation also checks all other supporting circuits before assigning a final scale. 
\paragraph{Diagonals.} A straight path of length equal to the polygon diameter supplies the second comparison. Its endpoints and orientation are computed from the actual vertices, so tall trapezoids are displayed without clipping or rescaling errors. 
\paragraph{Isosceles triangle with the top section removed.} As the upper base contracts, bent paths inherited from the triangular geometry become competitive. The truncation is retained explicitly in the quadrilateral oracle: these paths are not drawn by simply cutting an illustration of a triangular solution. The numerical changes between the three categories are evaluated using the full polygon and path coordinates.

The comparison with isosceles triangles has a precise geometric basis.
For $0<\alpha<\pi/2$, extend the two nonparallel sides until they meet at
$V=(1/2,\tfrac12\tan\alpha)$. The resulting isosceles triangle
$T_\alpha=\operatorname{conv}\{(0,0),(1,0),V\}$ contains the trapezoid,
and its base angles are both $\alpha$. As $s$ increases to
$1/(2\cos\alpha)$, the upper base contracts to zero and the trapezoid
becomes exactly $T_\alpha$. Thus the triangular limit can be compared with the
isosceles-triangle investigations in~\cite{besicovitch,coultonmovshovich,movshovich2012,gibbs,dengtriangle}
after rescaling to the same unit base.
At $\alpha=60^\circ$ this is the equilateral limit. Away from that limit,
the upper base supplies an additional supporting constraint: a triangular
path may fit differently after truncation, so its scale is determined
using all four trapezoid sides. The observed bent and curved portions
are consistent with the range of triangle path morphologies.

The relation to parallelograms is different. At the same $(s,\alpha)$,
the parallelogram has vertices $(0,0)$, $(1,0)$,
$(1+s\cos\alpha,s\sin\alpha)$, and
$(s\cos\alpha,s\sin\alpha)$ and contains the proper trapezoid.
Both forests share the height $h=s\sin\alpha$ and therefore the
horizontal-strip comparison, but a parallelogram has two pairs of
parallel sides whereas a trapezoid has only one. The converging sides
and upper cap change which contacts determine fitting.
In the right-angle limit both families reduce to a rectangle.

\subsection{Arbitrary Convex and Concave Quadrilateral}
\label{sec:arbitrary}

The last gallery in Figure 9 uses the four-parameter representation of
Section~\ref{sec:quadrilateral-parameters}. Self-intersecting parameter slots are identified and excluded from all path optimization and length summaries.

We compare
\begin{equation}\label{eq:mesh189}
p,q\in\{0.3,0.6,0.9\},\qquad
\alpha,\beta\in\{15^\circ,30^\circ,45^\circ,60^\circ,75^\circ,90^\circ\},
\qquad\alpha\le\beta.
\end{equation}
Both side-length parameters have increment $0.3$. The geometry check gives
91 convex, 76 concave, and 22 self-intersecting quadrilaterals. 

For each convex case, all feasible extreme equilibrium multipliers are enumerated; an assumed pair of normal triples is not used indiscriminately. For each concave case, the reflex vertex determines the two-triangle partition, and every segment is constrained either to one triangle or to cross the common internal diagonal. Candidate ranking uses a refined worst-orientation fit, while the coordinate files retain both the observed fit scale and the conservative scale used for the displayed escape path. This separation makes the effect of the angular guard directly inspectable. For the additional angular bound, let $B(c,r)$ lie in the polygon kernel,
and let $R$ bound the distance of every point of the raw path from a chosen
path vertex. If $g$ is an upper bound on every fit scale on an $n$-angle
uniform grid, then the kernel-contraction estimate gives
\begin{equation}\label{eq:kernel-numerical}
\sigma_Q^*(P)\le
\min\!\left\{\sigma_{\operatorname{conv}Q}^*(P),
\frac{g}{1-gR\,2\sin(\pi/(2n))/r}\right\},
\end{equation}
when the displayed denominator is positive. Otherwise the independently computed convex-hull bound is retained. Angular interval refinement can further reduce this guard. This estimate is evaluated with numerical padding; the complete derivation and floating-point qualifications are provided in the supplementary verification report. The principal reported quantities remain the computed path lengths, coordinates, containment residuals, and resolution comparisons.

The connection to general triangles is particularly direct for a simple concave quadrilateral: its convex hull $T=\operatorname{conv}Q$ is a triangle. The support-function formulation for arbitrary triangles in~\cite{dengtriangle} therefore provides a natural hull comparison and an efficient directional upper bound for fitting. Since $Q\subset T$, $L^*(Q)\le L^*(T)$, and any verified escape path for $T$ is also a valid comparison path for $Q$. A triangle has one normalized positive equilibrium relation among its three side normals, so its containment condition reduces to a single support branch. For the concave polygon, the same support data describe only $T$ and omit the reentrant notch. The actual computation must retain the two constituent triangles and the order in which each complete path segment crosses their shared diagonal.

This difference explains both the similarities and the departures in the gallery. If the critical triangular placement already avoids the removed notch, the hull and concave fitting scales can agree for that candidate. If a segment passes through the notch, the triangular placement is inadmissible in $Q$, even when its segment endpoints lie inside $Q$. A smaller fitting scale and a differently bent path may then be favorable. The asymmetry of a general triangular hull and the location of the reflex vertex both influence the result, so the isosceles-triangle symmetry cannot be imposed on a general concave case. As the reflex vertex approaches its adjacent hull edge, the notch closes and the triangular containment problem is recovered in the geometric limit. 
\FloatBarrier
\begingroup
\def\LTcaptype{figure}
\makeatletter\def\@captype{figure}\let\caption@stepcounter\@gobble\makeatother
\setlength{\LTpre}{5pt}
\setlength{\LTpost}{5pt}
\setlength{\LTcapwidth}{\linewidth}
\setlength{\tabcolsep}{0pt}
\captionsetup{position=bottom,skip=4pt}
\begin{longtable}{@{}c@{\hspace{0.011000\linewidth}}c@{\hspace{0.011000\linewidth}}c@{}}
\endfirsthead
\endhead
\caption[]{Arbitrary quadrilateral paths and excluded: 189 cases. (continued on following page).}\\
\endfoot
\caption{Arbitrary quadrilateral paths and excluded slots: 189 cases in parameter order.}\label{fig:arbitrary-1}\\
\endlastfoot
\begin{subfigure}[b]{0.326\linewidth}\centering
\includegraphics[width=\linewidth,height=119pt,keepaspectratio]{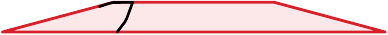}
\par{\scriptsize \mbox{$p=0.3,\ q=0.3,\ \alpha=15^\circ,\ \beta=15^\circ$}\par}
\end{subfigure}
&
\begin{subfigure}[b]{0.326\linewidth}\centering
\includegraphics[width=\linewidth,height=119pt,keepaspectratio]{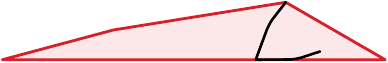}
\par{\scriptsize \mbox{$p=0.3,\ q=0.3,\ \alpha=15^\circ,\ \beta=30^\circ$}\par}
\end{subfigure}
&
\begin{subfigure}[b]{0.326\linewidth}\centering
\includegraphics[width=\linewidth,height=119pt,keepaspectratio]{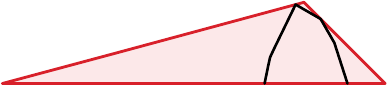}
\par{\scriptsize \mbox{$p=0.3,\ q=0.3,\ \alpha=15^\circ,\ \beta=45^\circ$}\par}
\end{subfigure}
\\
\begin{subfigure}[b]{0.326\linewidth}\centering
\includegraphics[width=\linewidth,height=119pt,keepaspectratio]{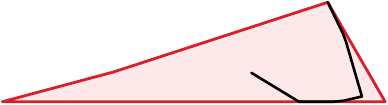}
\par{\scriptsize \mbox{$p=0.3,\ q=0.3,\ \alpha=15^\circ,\ \beta=60^\circ$}\par}
\end{subfigure}
&
\begin{subfigure}[b]{0.326\linewidth}\centering
\includegraphics[width=\linewidth,height=119pt,keepaspectratio]{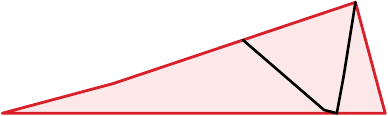}
\par{\scriptsize \mbox{$p=0.3,\ q=0.3,\ \alpha=15^\circ,\ \beta=75^\circ$}\par}
\end{subfigure}
&
\begin{subfigure}[b]{0.326\linewidth}\centering
\includegraphics[width=\linewidth,height=119pt,keepaspectratio]{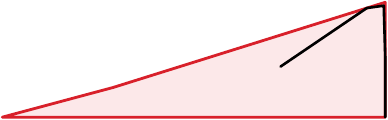}
\par{\scriptsize \mbox{$p=0.3,\ q=0.3,\ \alpha=15^\circ,\ \beta=90^\circ$}\par}
\end{subfigure}
\\
\begin{subfigure}[b]{0.326\linewidth}\centering
\includegraphics[width=\linewidth,height=119pt,keepaspectratio]{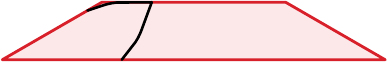}
\par{\scriptsize \mbox{$p=0.3,\ q=0.3,\ \alpha=30^\circ,\ \beta=30^\circ$}\par}
\end{subfigure}
&
\begin{subfigure}[b]{0.326\linewidth}\centering
\includegraphics[width=\linewidth,height=119pt,keepaspectratio]{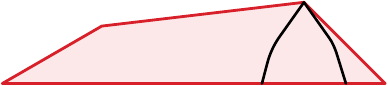}
\par{\scriptsize \mbox{$p=0.3,\ q=0.3,\ \alpha=30^\circ,\ \beta=45^\circ$}\par}
\end{subfigure}
&
\begin{subfigure}[b]{0.326\linewidth}\centering
\includegraphics[width=\linewidth,height=119pt,keepaspectratio]{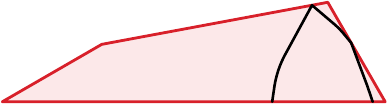}
\par{\scriptsize \mbox{$p=0.3,\ q=0.3,\ \alpha=30^\circ,\ \beta=60^\circ$}\par}
\end{subfigure}
\\
\begin{subfigure}[b]{0.326\linewidth}\centering
\includegraphics[width=\linewidth,height=119pt,keepaspectratio]{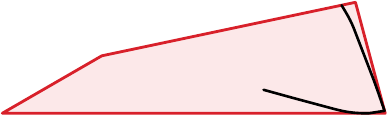}
\par{\scriptsize \mbox{$p=0.3,\ q=0.3,\ \alpha=30^\circ,\ \beta=75^\circ$}\par}
\end{subfigure}
&
\begin{subfigure}[b]{0.326\linewidth}\centering
\includegraphics[width=\linewidth,height=119pt,keepaspectratio]{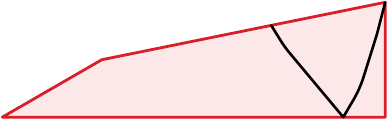}
\par{\scriptsize \mbox{$p=0.3,\ q=0.3,\ \alpha=30^\circ,\ \beta=90^\circ$}\par}
\end{subfigure}
&
\begin{subfigure}[b]{0.326\linewidth}\centering
\includegraphics[width=\linewidth,height=119pt,keepaspectratio]{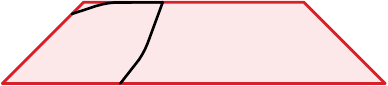}
\par{\scriptsize \mbox{$p=0.3,\ q=0.3,\ \alpha=45^\circ,\ \beta=45^\circ$}\par}
\end{subfigure}
\\
\begin{subfigure}[b]{0.326\linewidth}\centering
\includegraphics[width=\linewidth,height=119pt,keepaspectratio]{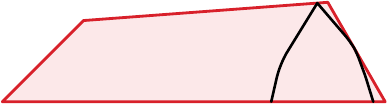}
\par{\scriptsize \mbox{$p=0.3,\ q=0.3,\ \alpha=45^\circ,\ \beta=60^\circ$}\par}
\end{subfigure}
&
\begin{subfigure}[b]{0.326\linewidth}\centering
\includegraphics[width=\linewidth,height=119pt,keepaspectratio]{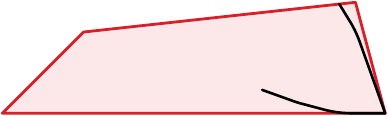}
\par{\scriptsize \mbox{$p=0.3,\ q=0.3,\ \alpha=45^\circ,\ \beta=75^\circ$}\par}
\end{subfigure}
&
\begin{subfigure}[b]{0.326\linewidth}\centering
\includegraphics[width=\linewidth,height=119pt,keepaspectratio]{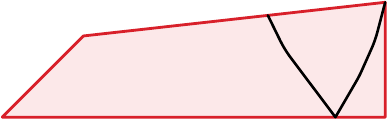}
\par{\scriptsize \mbox{$p=0.3,\ q=0.3,\ \alpha=45^\circ,\ \beta=90^\circ$}\par}
\end{subfigure}
\\
\begin{subfigure}[b]{0.326\linewidth}\centering
\includegraphics[width=\linewidth,height=119pt,keepaspectratio]{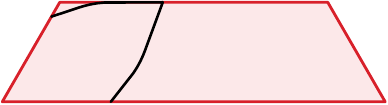}
\par{\scriptsize \mbox{$p=0.3,\ q=0.3,\ \alpha=60^\circ,\ \beta=60^\circ$}\par}
\end{subfigure}
&
\begin{subfigure}[b]{0.326\linewidth}\centering
\includegraphics[width=\linewidth,height=119pt,keepaspectratio]{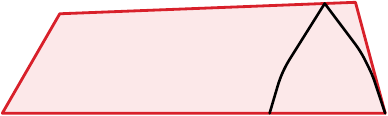}
\par{\scriptsize \mbox{$p=0.3,\ q=0.3,\ \alpha=60^\circ,\ \beta=75^\circ$}\par}
\end{subfigure}
&
\begin{subfigure}[b]{0.326\linewidth}\centering
\includegraphics[width=\linewidth,height=119pt,keepaspectratio]{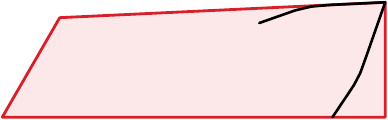}
\par{\scriptsize \mbox{$p=0.3,\ q=0.3,\ \alpha=60^\circ,\ \beta=90^\circ$}\par}
\end{subfigure}
\\
\begin{subfigure}[b]{0.326\linewidth}\centering
\includegraphics[width=\linewidth,height=119pt,keepaspectratio]{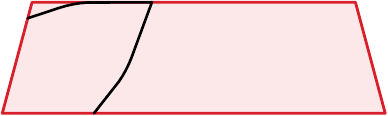}
\par{\scriptsize \mbox{$p=0.3,\ q=0.3,\ \alpha=75^\circ,\ \beta=75^\circ$}\par}
\end{subfigure}
&
\begin{subfigure}[b]{0.326\linewidth}\centering
\includegraphics[width=\linewidth,height=119pt,keepaspectratio]{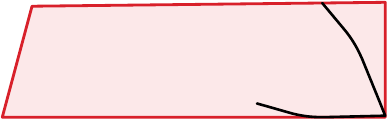}
\par{\scriptsize \mbox{$p=0.3,\ q=0.3,\ \alpha=75^\circ,\ \beta=90^\circ$}\par}
\end{subfigure}
&
\begin{subfigure}[b]{0.326\linewidth}\centering
\includegraphics[width=\linewidth,height=119pt,keepaspectratio]{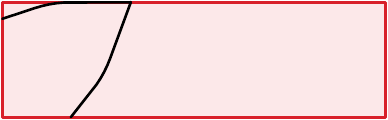}
\par{\scriptsize \mbox{$p=0.3,\ q=0.3,\ \alpha=90^\circ,\ \beta=90^\circ$}\par}
\end{subfigure}
\\
\begin{subfigure}[b]{0.326\linewidth}\centering
\includegraphics[width=\linewidth,height=119pt,keepaspectratio]{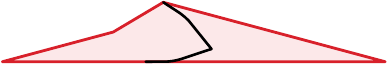}
\par{\scriptsize \mbox{$p=0.3,\ q=0.6,\ \alpha=15^\circ,\ \beta=15^\circ$}\par}
\end{subfigure}
&
\begin{subfigure}[b]{0.326\linewidth}\centering
\includegraphics[width=\linewidth,height=119pt,keepaspectratio]{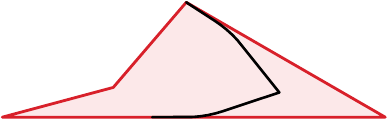}
\par{\scriptsize \mbox{$p=0.3,\ q=0.6,\ \alpha=15^\circ,\ \beta=30^\circ$}\par}
\end{subfigure}
&
\begin{subfigure}[b]{0.326\linewidth}\centering
\includegraphics[width=\linewidth,height=119pt,keepaspectratio]{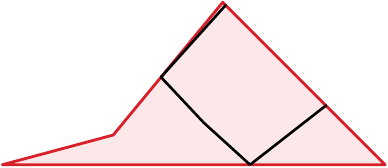}
\par{\scriptsize \mbox{$p=0.3,\ q=0.6,\ \alpha=15^\circ,\ \beta=45^\circ$}\par}
\end{subfigure}
\\
\begin{subfigure}[b]{0.326\linewidth}\centering
\includegraphics[width=\linewidth,height=119pt,keepaspectratio]{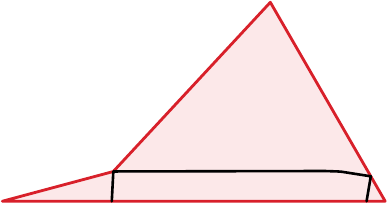}
\par{\scriptsize \mbox{$p=0.3,\ q=0.6,\ \alpha=15^\circ,\ \beta=60^\circ$}\par}
\end{subfigure}
&
\begin{subfigure}[b]{0.326\linewidth}\centering
\includegraphics[width=\linewidth,height=119pt,keepaspectratio]{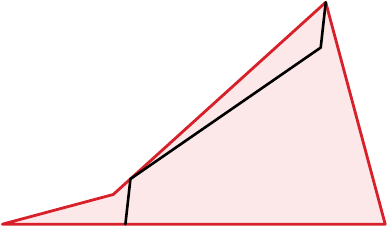}
\par{\scriptsize \mbox{$p=0.3,\ q=0.6,\ \alpha=15^\circ,\ \beta=75^\circ$}\par}
\end{subfigure}
&
\begin{subfigure}[b]{0.326\linewidth}\centering
\includegraphics[width=\linewidth,height=119pt,keepaspectratio]{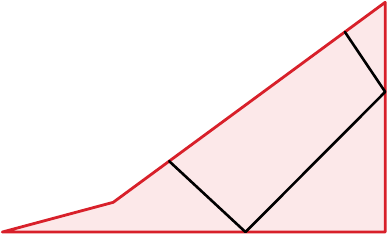}
\par{\scriptsize \mbox{$p=0.3,\ q=0.6,\ \alpha=15^\circ,\ \beta=90^\circ$}\par}
\end{subfigure}
\\
\begin{subfigure}[b]{0.326\linewidth}\centering
\includegraphics[width=\linewidth,height=119pt,keepaspectratio]{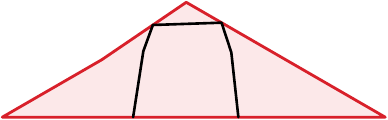}
\par{\scriptsize \mbox{$p=0.3,\ q=0.6,\ \alpha=30^\circ,\ \beta=30^\circ$}\par}
\end{subfigure}
&
\begin{subfigure}[b]{0.326\linewidth}\centering
\includegraphics[width=\linewidth,height=119pt,keepaspectratio]{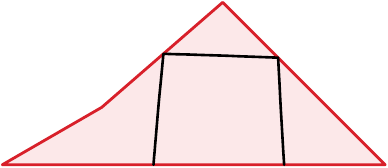}
\par{\scriptsize \mbox{$p=0.3,\ q=0.6,\ \alpha=30^\circ,\ \beta=45^\circ$}\par}
\end{subfigure}
&
\begin{subfigure}[b]{0.326\linewidth}\centering
\includegraphics[width=\linewidth,height=119pt,keepaspectratio]{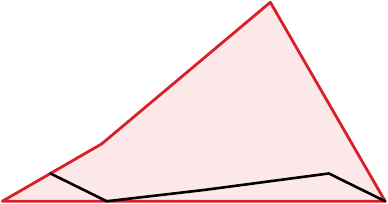}
\par{\scriptsize \mbox{$p=0.3,\ q=0.6,\ \alpha=30^\circ,\ \beta=60^\circ$}\par}
\end{subfigure}
\\
\begin{subfigure}[b]{0.326\linewidth}\centering
\includegraphics[width=\linewidth,height=119pt,keepaspectratio]{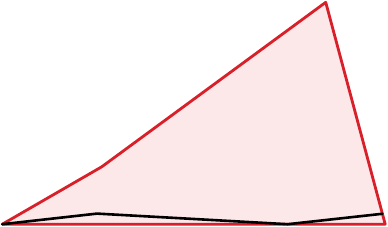}
\par{\scriptsize \mbox{$p=0.3,\ q=0.6,\ \alpha=30^\circ,\ \beta=75^\circ$}\par}
\end{subfigure}
&
\begin{subfigure}[b]{0.326\linewidth}\centering
\includegraphics[width=\linewidth,height=119pt,keepaspectratio]{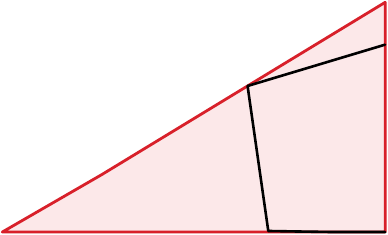}
\par{\scriptsize \mbox{$p=0.3,\ q=0.6,\ \alpha=30^\circ,\ \beta=90^\circ$}\par}
\end{subfigure}
&
\begin{subfigure}[b]{0.326\linewidth}\centering
\includegraphics[width=\linewidth,height=119pt,keepaspectratio]{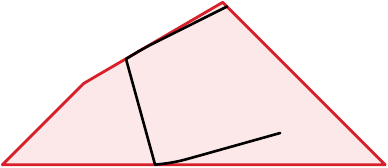}
\par{\scriptsize \mbox{$p=0.3,\ q=0.6,\ \alpha=45^\circ,\ \beta=45^\circ$}\par}
\end{subfigure}
\\
\begin{subfigure}[b]{0.326\linewidth}\centering
\includegraphics[width=\linewidth,height=119pt,keepaspectratio]{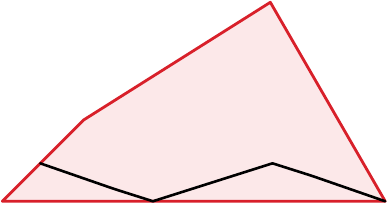}
\par{\scriptsize \mbox{$p=0.3,\ q=0.6,\ \alpha=45^\circ,\ \beta=60^\circ$}\par}
\end{subfigure}
&
\begin{subfigure}[b]{0.326\linewidth}\centering
\includegraphics[width=\linewidth,height=119pt,keepaspectratio]{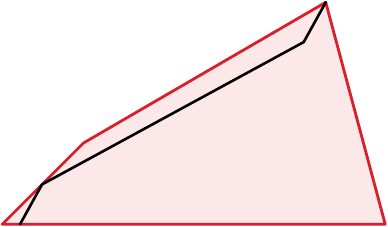}
\par{\scriptsize \mbox{$p=0.3,\ q=0.6,\ \alpha=45^\circ,\ \beta=75^\circ$}\par}
\end{subfigure}
&
\begin{subfigure}[b]{0.326\linewidth}\centering
\includegraphics[width=\linewidth,height=119pt,keepaspectratio]{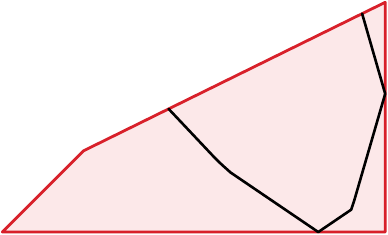}
\par{\scriptsize \mbox{$p=0.3,\ q=0.6,\ \alpha=45^\circ,\ \beta=90^\circ$}\par}
\end{subfigure}
\\
\begin{subfigure}[b]{0.326\linewidth}\centering
\includegraphics[width=\linewidth,height=119pt,keepaspectratio]{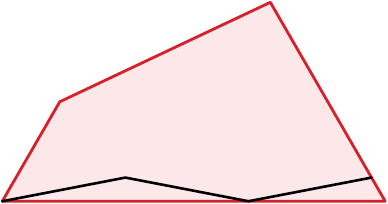}
\par{\scriptsize \mbox{$p=0.3,\ q=0.6,\ \alpha=60^\circ,\ \beta=60^\circ$}\par}
\end{subfigure}
&
\begin{subfigure}[b]{0.326\linewidth}\centering
\includegraphics[width=\linewidth,height=119pt,keepaspectratio]{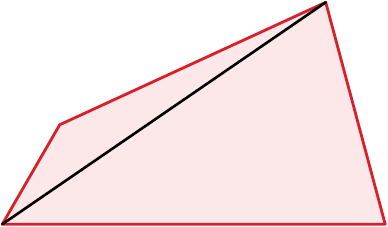}
\par{\scriptsize \mbox{$p=0.3,\ q=0.6,\ \alpha=60^\circ,\ \beta=75^\circ$}\par}
\end{subfigure}
&
\begin{subfigure}[b]{0.326\linewidth}\centering
\includegraphics[width=\linewidth,height=119pt,keepaspectratio]{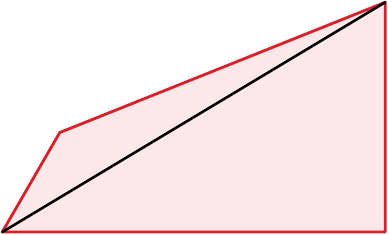}
\par{\scriptsize \mbox{$p=0.3,\ q=0.6,\ \alpha=60^\circ,\ \beta=90^\circ$}\par}
\end{subfigure}
\\
\begin{subfigure}[b]{0.326\linewidth}\centering
\includegraphics[width=\linewidth,height=119pt,keepaspectratio]{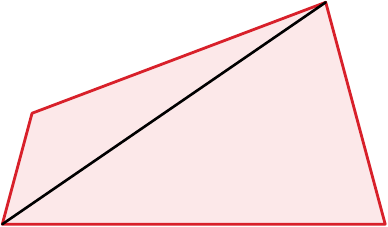}
\par{\scriptsize \mbox{$p=0.3,\ q=0.6,\ \alpha=75^\circ,\ \beta=75^\circ$}\par}
\end{subfigure}
&
\begin{subfigure}[b]{0.326\linewidth}\centering
\includegraphics[width=\linewidth,height=119pt,keepaspectratio]{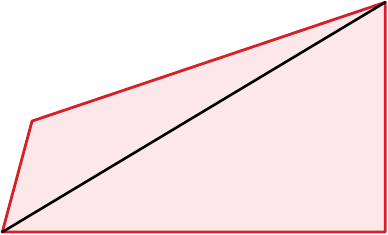}
\par{\scriptsize \mbox{$p=0.3,\ q=0.6,\ \alpha=75^\circ,\ \beta=90^\circ$}\par}
\end{subfigure}
&
\begin{subfigure}[b]{0.326\linewidth}\centering
\includegraphics[width=\linewidth,height=119pt,keepaspectratio]{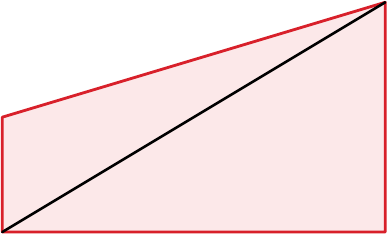}
\par{\scriptsize \mbox{$p=0.3,\ q=0.6,\ \alpha=90^\circ,\ \beta=90^\circ$}\par}
\end{subfigure}
\\
\begin{subfigure}[b]{0.326\linewidth}\centering
\includegraphics[width=\linewidth,height=119pt,keepaspectratio]{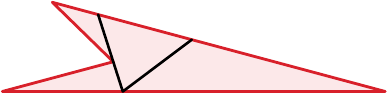}
\par{\scriptsize \mbox{$p=0.3,\ q=0.9,\ \alpha=15^\circ,\ \beta=15^\circ$}\par}
\end{subfigure}
&
\begin{subfigure}[b]{0.326\linewidth}\centering
\includegraphics[width=\linewidth,height=119pt,keepaspectratio]{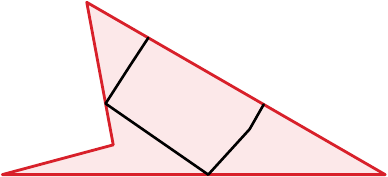}
\par{\scriptsize \mbox{$p=0.3,\ q=0.9,\ \alpha=15^\circ,\ \beta=30^\circ$}\par}
\end{subfigure}
&
\begin{subfigure}[b]{0.326\linewidth}\centering
\includegraphics[width=\linewidth,height=119pt,keepaspectratio]{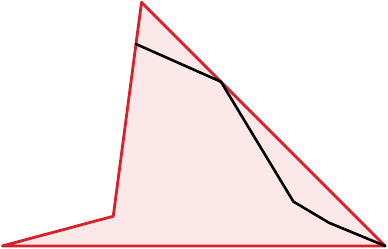}
\par{\scriptsize \mbox{$p=0.3,\ q=0.9,\ \alpha=15^\circ,\ \beta=45^\circ$}\par}
\end{subfigure}
\\
\begin{subfigure}[b]{0.326\linewidth}\centering
\includegraphics[width=\linewidth,height=119pt,keepaspectratio]{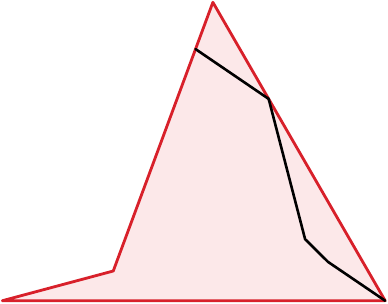}
\par{\scriptsize \mbox{$p=0.3,\ q=0.9,\ \alpha=15^\circ,\ \beta=60^\circ$}\par}
\end{subfigure}
&
\begin{subfigure}[b]{0.326\linewidth}\centering
\includegraphics[width=\linewidth,height=119pt,keepaspectratio]{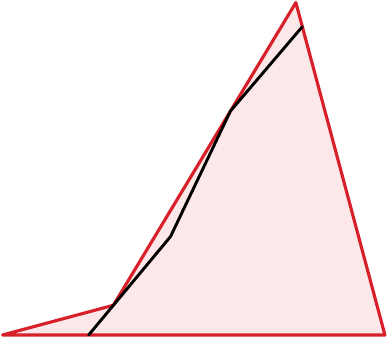}
\par{\scriptsize \mbox{$p=0.3,\ q=0.9,\ \alpha=15^\circ,\ \beta=75^\circ$}\par}
\end{subfigure}
&
\begin{subfigure}[b]{0.326\linewidth}\centering
\includegraphics[width=\linewidth,height=119pt,keepaspectratio]{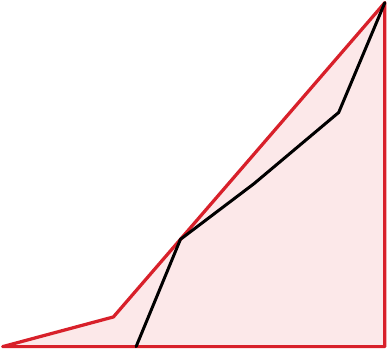}
\par{\scriptsize \mbox{$p=0.3,\ q=0.9,\ \alpha=15^\circ,\ \beta=90^\circ$}\par}
\end{subfigure}
\\
\begin{subfigure}[b]{0.326\linewidth}\centering
\includegraphics[width=\linewidth,height=119pt,keepaspectratio]{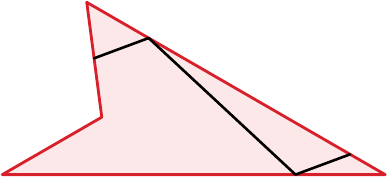}
\par{\scriptsize \mbox{$p=0.3,\ q=0.9,\ \alpha=30^\circ,\ \beta=30^\circ$}\par}
\end{subfigure}
&
\begin{subfigure}[b]{0.326\linewidth}\centering
\includegraphics[width=\linewidth,height=119pt,keepaspectratio]{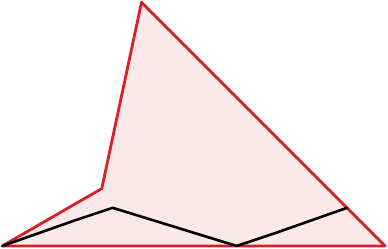}
\par{\scriptsize \mbox{$p=0.3,\ q=0.9,\ \alpha=30^\circ,\ \beta=45^\circ$}\par}
\end{subfigure}
&
\begin{subfigure}[b]{0.326\linewidth}\centering
\includegraphics[width=\linewidth,height=119pt,keepaspectratio]{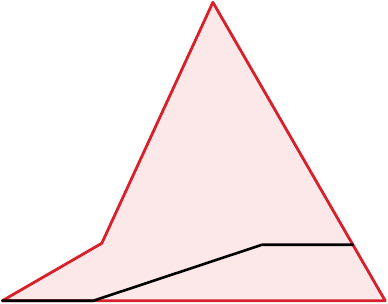}
\par{\scriptsize \mbox{$p=0.3,\ q=0.9,\ \alpha=30^\circ,\ \beta=60^\circ$}\par}
\end{subfigure}
\\
\begin{subfigure}[b]{0.326\linewidth}\centering
\includegraphics[width=\linewidth,height=119pt,keepaspectratio]{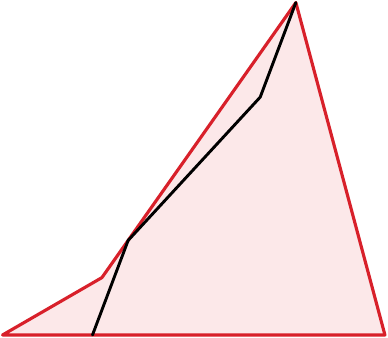}
\par{\scriptsize \mbox{$p=0.3,\ q=0.9,\ \alpha=30^\circ,\ \beta=75^\circ$}\par}
\end{subfigure}
&
\begin{subfigure}[b]{0.326\linewidth}\centering
\includegraphics[width=\linewidth,height=119pt,keepaspectratio]{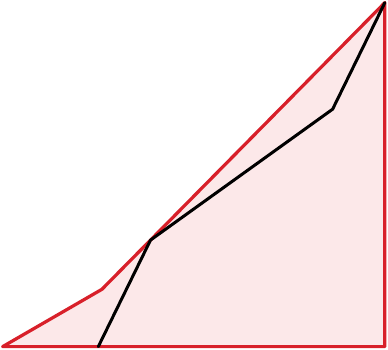}
\par{\scriptsize \mbox{$p=0.3,\ q=0.9,\ \alpha=30^\circ,\ \beta=90^\circ$}\par}
\end{subfigure}
&
\begin{subfigure}[b]{0.326\linewidth}\centering
\includegraphics[width=\linewidth,height=119pt,keepaspectratio]{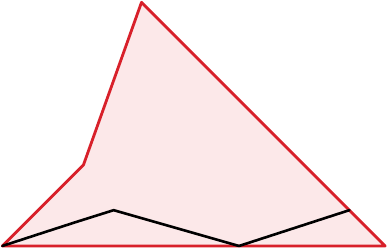}
\par{\scriptsize \mbox{$p=0.3,\ q=0.9,\ \alpha=45^\circ,\ \beta=45^\circ$}\par}
\end{subfigure}
\\
\begin{subfigure}[b]{0.326\linewidth}\centering
\includegraphics[width=\linewidth,height=119pt,keepaspectratio]{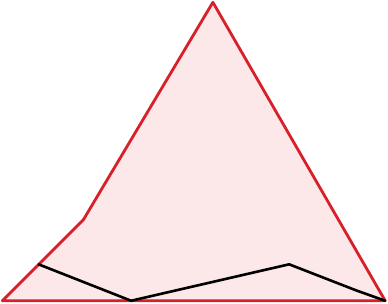}
\par{\scriptsize \mbox{$p=0.3,\ q=0.9,\ \alpha=45^\circ,\ \beta=60^\circ$}\par}
\end{subfigure}
&
\begin{subfigure}[b]{0.326\linewidth}\centering
\includegraphics[width=\linewidth,height=119pt,keepaspectratio]{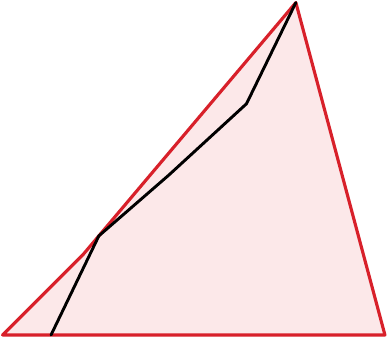}
\par{\scriptsize \mbox{$p=0.3,\ q=0.9,\ \alpha=45^\circ,\ \beta=75^\circ$}\par}
\end{subfigure}
&
\begin{subfigure}[b]{0.326\linewidth}\centering
\includegraphics[width=\linewidth,height=119pt,keepaspectratio]{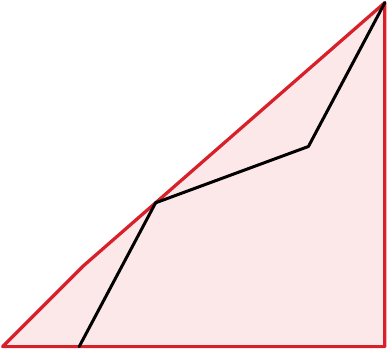}
\par{\scriptsize \mbox{$p=0.3,\ q=0.9,\ \alpha=45^\circ,\ \beta=90^\circ$}\par}
\end{subfigure}
\\
\begin{subfigure}[b]{0.326\linewidth}\centering
\includegraphics[width=\linewidth,height=119pt,keepaspectratio]{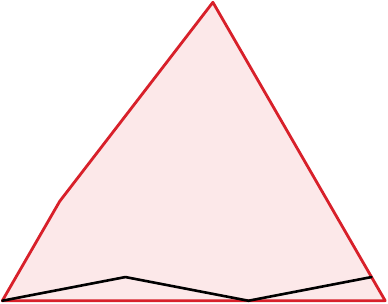}
\par{\scriptsize \mbox{$p=0.3,\ q=0.9,\ \alpha=60^\circ,\ \beta=60^\circ$}\par}
\end{subfigure}
&
\begin{subfigure}[b]{0.326\linewidth}\centering
\includegraphics[width=\linewidth,height=119pt,keepaspectratio]{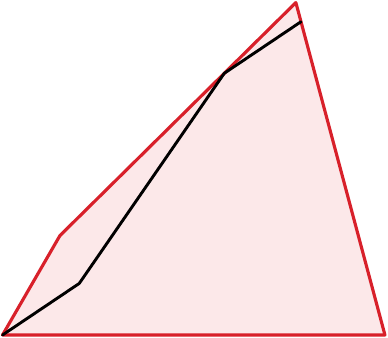}
\par{\scriptsize \mbox{$p=0.3,\ q=0.9,\ \alpha=60^\circ,\ \beta=75^\circ$}\par}
\end{subfigure}
&
\begin{subfigure}[b]{0.326\linewidth}\centering
\includegraphics[width=\linewidth,height=119pt,keepaspectratio]{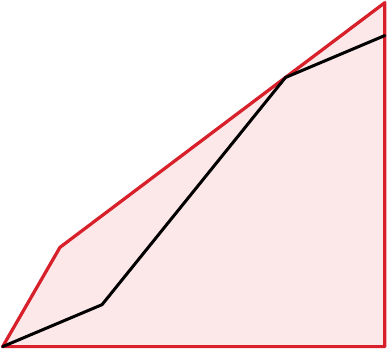}
\par{\scriptsize \mbox{$p=0.3,\ q=0.9,\ \alpha=60^\circ,\ \beta=90^\circ$}\par}
\end{subfigure}
\\
\begin{subfigure}[b]{0.326\linewidth}\centering
\includegraphics[width=\linewidth,height=119pt,keepaspectratio]{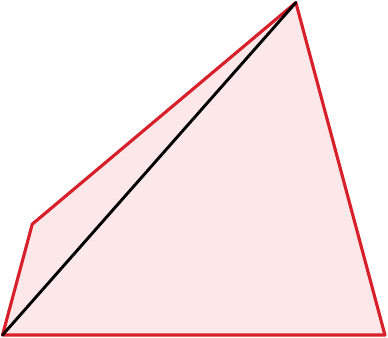}
\par{\scriptsize \mbox{$p=0.3,\ q=0.9,\ \alpha=75^\circ,\ \beta=75^\circ$}\par}
\end{subfigure}
&
\begin{subfigure}[b]{0.326\linewidth}\centering
\includegraphics[width=\linewidth,height=119pt,keepaspectratio]{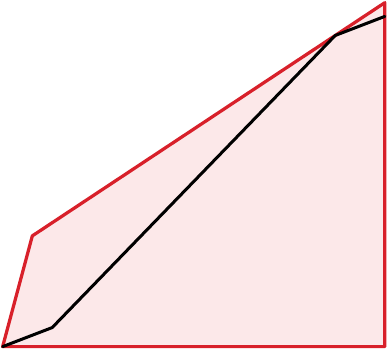}
\par{\scriptsize \mbox{$p=0.3,\ q=0.9,\ \alpha=75^\circ,\ \beta=90^\circ$}\par}
\end{subfigure}
&
\begin{subfigure}[b]{0.326\linewidth}\centering
\includegraphics[width=\linewidth,height=119pt,keepaspectratio]{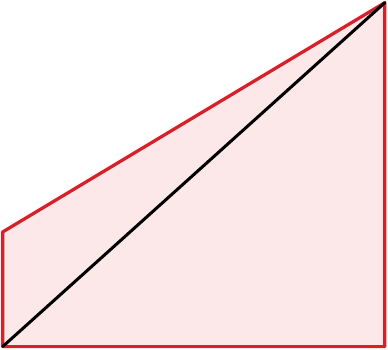}
\par{\scriptsize \mbox{$p=0.3,\ q=0.9,\ \alpha=90^\circ,\ \beta=90^\circ$}\par}
\end{subfigure}
\\
\begin{subfigure}[b]{0.326\linewidth}\centering
\includegraphics[width=\linewidth,height=119pt,keepaspectratio]{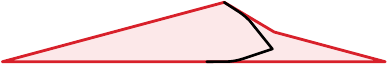}
\par{\scriptsize \mbox{$p=0.6,\ q=0.3,\ \alpha=15^\circ,\ \beta=15^\circ$}\par}
\end{subfigure}
&
\begin{subfigure}[b]{0.326\linewidth}\centering
\includegraphics[width=\linewidth,height=119pt,keepaspectratio]{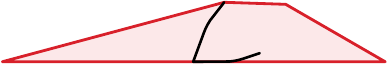}
\par{\scriptsize \mbox{$p=0.6,\ q=0.3,\ \alpha=15^\circ,\ \beta=30^\circ$}\par}
\end{subfigure}
&
\begin{subfigure}[b]{0.326\linewidth}\centering
\includegraphics[width=\linewidth,height=119pt,keepaspectratio]{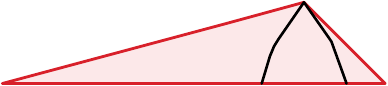}
\par{\scriptsize \mbox{$p=0.6,\ q=0.3,\ \alpha=15^\circ,\ \beta=45^\circ$}\par}
\end{subfigure}
\\
\begin{subfigure}[b]{0.326\linewidth}\centering
\includegraphics[width=\linewidth,height=119pt,keepaspectratio]{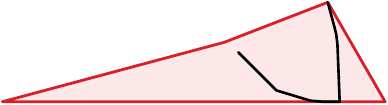}
\par{\scriptsize \mbox{$p=0.6,\ q=0.3,\ \alpha=15^\circ,\ \beta=60^\circ$}\par}
\end{subfigure}
&
\begin{subfigure}[b]{0.326\linewidth}\centering
\includegraphics[width=\linewidth,height=119pt,keepaspectratio]{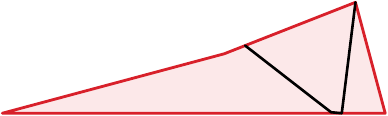}
\par{\scriptsize \mbox{$p=0.6,\ q=0.3,\ \alpha=15^\circ,\ \beta=75^\circ$}\par}
\end{subfigure}
&
\begin{subfigure}[b]{0.326\linewidth}\centering
\includegraphics[width=\linewidth,height=119pt,keepaspectratio]{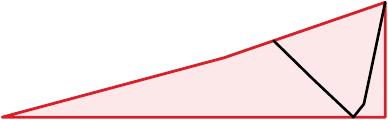}
\par{\scriptsize \mbox{$p=0.6,\ q=0.3,\ \alpha=15^\circ,\ \beta=90^\circ$}\par}
\end{subfigure}
\\
\begin{subfigure}[b]{0.326\linewidth}\centering
\includegraphics[width=\linewidth,height=119pt,keepaspectratio]{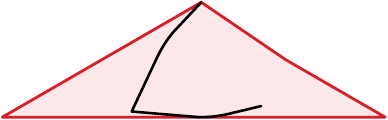}
\par{\scriptsize \mbox{$p=0.6,\ q=0.3,\ \alpha=30^\circ,\ \beta=30^\circ$}\par}
\end{subfigure}
&
\begin{subfigure}[b]{0.326\linewidth}\centering
\includegraphics[width=\linewidth,height=119pt,keepaspectratio]{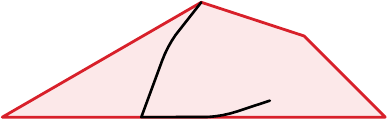}
\par{\scriptsize \mbox{$p=0.6,\ q=0.3,\ \alpha=30^\circ,\ \beta=45^\circ$}\par}
\end{subfigure}
&
\begin{subfigure}[b]{0.326\linewidth}\centering
\includegraphics[width=\linewidth,height=119pt,keepaspectratio]{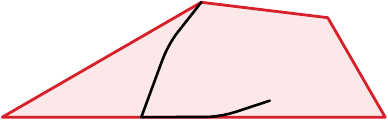}
\par{\scriptsize \mbox{$p=0.6,\ q=0.3,\ \alpha=30^\circ,\ \beta=60^\circ$}\par}
\end{subfigure}
\\
\begin{subfigure}[b]{0.326\linewidth}\centering
\includegraphics[width=\linewidth,height=119pt,keepaspectratio]{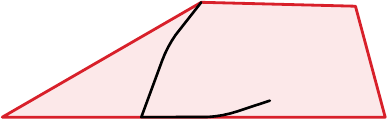}
\par{\scriptsize \mbox{$p=0.6,\ q=0.3,\ \alpha=30^\circ,\ \beta=75^\circ$}\par}
\end{subfigure}
&
\begin{subfigure}[b]{0.326\linewidth}\centering
\includegraphics[width=\linewidth,height=119pt,keepaspectratio]{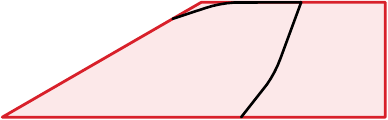}
\par{\scriptsize \mbox{$p=0.6,\ q=0.3,\ \alpha=30^\circ,\ \beta=90^\circ$}\par}
\end{subfigure}
&
\begin{subfigure}[b]{0.326\linewidth}\centering
\includegraphics[width=\linewidth,height=119pt,keepaspectratio]{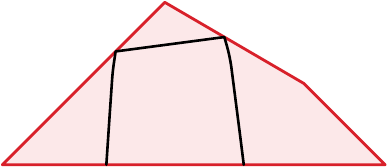}
\par{\scriptsize \mbox{$p=0.6,\ q=0.3,\ \alpha=45^\circ,\ \beta=45^\circ$}\par}
\end{subfigure}
\\
\begin{subfigure}[b]{0.326\linewidth}\centering
\includegraphics[width=\linewidth,height=119pt,keepaspectratio]{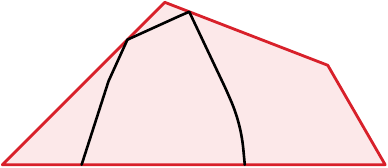}
\par{\scriptsize \mbox{$p=0.6,\ q=0.3,\ \alpha=45^\circ,\ \beta=60^\circ$}\par}
\end{subfigure}
&
\begin{subfigure}[b]{0.326\linewidth}\centering
\includegraphics[width=\linewidth,height=119pt,keepaspectratio]{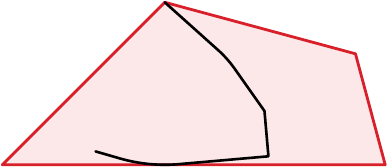}
\par{\scriptsize \mbox{$p=0.6,\ q=0.3,\ \alpha=45^\circ,\ \beta=75^\circ$}\par}
\end{subfigure}
&
\begin{subfigure}[b]{0.326\linewidth}\centering
\includegraphics[width=\linewidth,height=119pt,keepaspectratio]{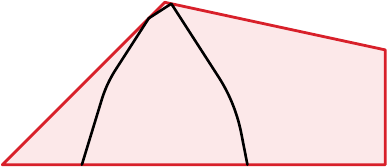}
\par{\scriptsize \mbox{$p=0.6,\ q=0.3,\ \alpha=45^\circ,\ \beta=90^\circ$}\par}
\end{subfigure}
\\
\begin{subfigure}[b]{0.326\linewidth}\centering
\includegraphics[width=\linewidth,height=119pt,keepaspectratio]{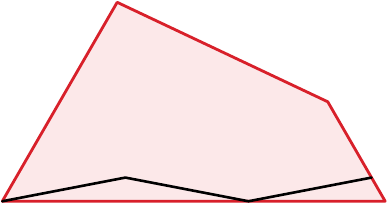}
\par{\scriptsize \mbox{$p=0.6,\ q=0.3,\ \alpha=60^\circ,\ \beta=60^\circ$}\par}
\end{subfigure}
&
\begin{subfigure}[b]{0.326\linewidth}\centering
\includegraphics[width=\linewidth,height=119pt,keepaspectratio]{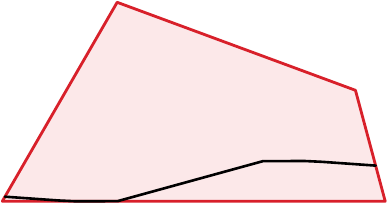}
\par{\scriptsize \mbox{$p=0.6,\ q=0.3,\ \alpha=60^\circ,\ \beta=75^\circ$}\par}
\end{subfigure}
&
\begin{subfigure}[b]{0.326\linewidth}\centering
\includegraphics[width=\linewidth,height=119pt,keepaspectratio]{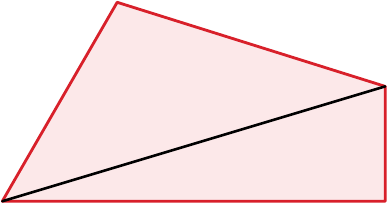}
\par{\scriptsize \mbox{$p=0.6,\ q=0.3,\ \alpha=60^\circ,\ \beta=90^\circ$}\par}
\end{subfigure}
\\
\begin{subfigure}[b]{0.326\linewidth}\centering
\includegraphics[width=\linewidth,height=119pt,keepaspectratio]{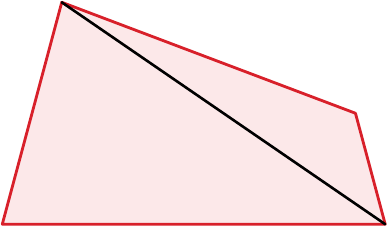}
\par{\scriptsize \mbox{$p=0.6,\ q=0.3,\ \alpha=75^\circ,\ \beta=75^\circ$}\par}
\end{subfigure}
&
\begin{subfigure}[b]{0.326\linewidth}\centering
\includegraphics[width=\linewidth,height=119pt,keepaspectratio]{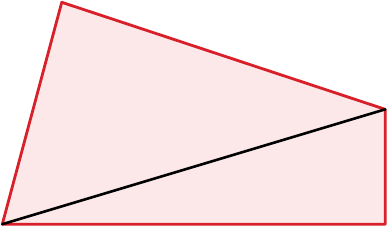}
\par{\scriptsize \mbox{$p=0.6,\ q=0.3,\ \alpha=75^\circ,\ \beta=90^\circ$}\par}
\end{subfigure}
&
\begin{subfigure}[b]{0.326\linewidth}\centering
\includegraphics[width=\linewidth,height=119pt,keepaspectratio]{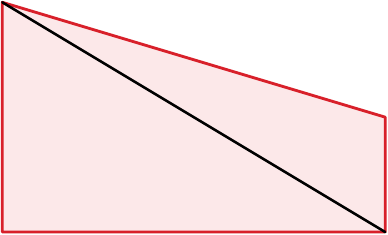}
\par{\scriptsize \mbox{$p=0.6,\ q=0.3,\ \alpha=90^\circ,\ \beta=90^\circ$}\par}
\end{subfigure}
\\
\begin{subfigure}[b]{0.326\linewidth}\centering
\includegraphics[width=\linewidth,height=119pt,keepaspectratio]{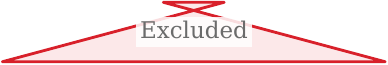}
\par{\scriptsize \mbox{$p=0.6,\ q=0.6,\ \alpha=15^\circ,\ \beta=15^\circ$}\par}
\end{subfigure}
&
\begin{subfigure}[b]{0.326\linewidth}\centering
\includegraphics[width=\linewidth,height=119pt,keepaspectratio]{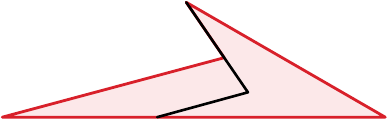}
\par{\scriptsize \mbox{$p=0.6,\ q=0.6,\ \alpha=15^\circ,\ \beta=30^\circ$}\par}
\end{subfigure}
&
\begin{subfigure}[b]{0.326\linewidth}\centering
\includegraphics[width=\linewidth,height=119pt,keepaspectratio]{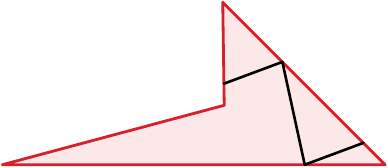}
\par{\scriptsize \mbox{$p=0.6,\ q=0.6,\ \alpha=15^\circ,\ \beta=45^\circ$}\par}
\end{subfigure}
\\
\begin{subfigure}[b]{0.326\linewidth}\centering
\includegraphics[width=\linewidth,height=119pt,keepaspectratio]{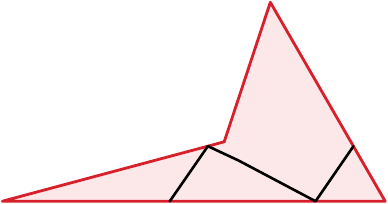}
\par{\scriptsize \mbox{$p=0.6,\ q=0.6,\ \alpha=15^\circ,\ \beta=60^\circ$}\par}
\end{subfigure}
&
\begin{subfigure}[b]{0.326\linewidth}\centering
\includegraphics[width=\linewidth,height=119pt,keepaspectratio]{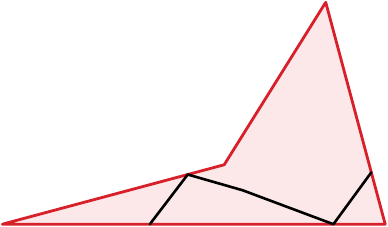}
\par{\scriptsize \mbox{$p=0.6,\ q=0.6,\ \alpha=15^\circ,\ \beta=75^\circ$}\par}
\end{subfigure}
&
\begin{subfigure}[b]{0.326\linewidth}\centering
\includegraphics[width=\linewidth,height=119pt,keepaspectratio]{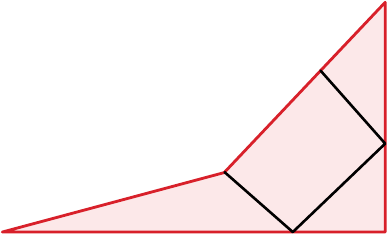}
\par{\scriptsize \mbox{$p=0.6,\ q=0.6,\ \alpha=15^\circ,\ \beta=90^\circ$}\par}
\end{subfigure}
\\
\begin{subfigure}[b]{0.326\linewidth}\centering
\includegraphics[width=\linewidth,height=119pt,keepaspectratio]{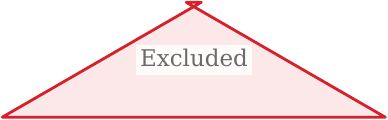}
\par{\scriptsize \mbox{$p=0.6,\ q=0.6,\ \alpha=30^\circ,\ \beta=30^\circ$}\par}
\end{subfigure}
&
\begin{subfigure}[b]{0.326\linewidth}\centering
\includegraphics[width=\linewidth,height=119pt,keepaspectratio]{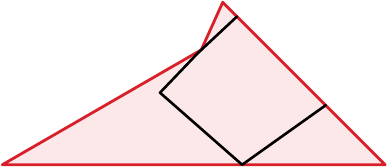}
\par{\scriptsize \mbox{$p=0.6,\ q=0.6,\ \alpha=30^\circ,\ \beta=45^\circ$}\par}
\end{subfigure}
&
\begin{subfigure}[b]{0.326\linewidth}\centering
\includegraphics[width=\linewidth,height=119pt,keepaspectratio]{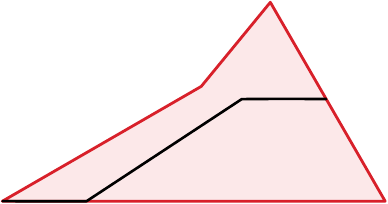}
\par{\scriptsize \mbox{$p=0.6,\ q=0.6,\ \alpha=30^\circ,\ \beta=60^\circ$}\par}
\end{subfigure}
\\
\begin{subfigure}[b]{0.326\linewidth}\centering
\includegraphics[width=\linewidth,height=119pt,keepaspectratio]{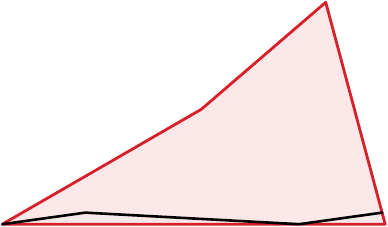}
\par{\scriptsize \mbox{$p=0.6,\ q=0.6,\ \alpha=30^\circ,\ \beta=75^\circ$}\par}
\end{subfigure}
&
\begin{subfigure}[b]{0.326\linewidth}\centering
\includegraphics[width=\linewidth,height=119pt,keepaspectratio]{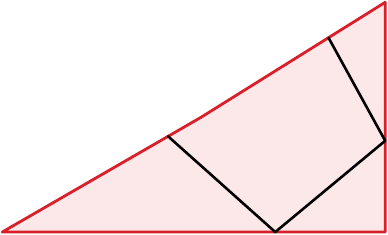}
\par{\scriptsize \mbox{$p=0.6,\ q=0.6,\ \alpha=30^\circ,\ \beta=90^\circ$}\par}
\end{subfigure}
&
\begin{subfigure}[b]{0.326\linewidth}\centering
\includegraphics[width=\linewidth,height=119pt,keepaspectratio]{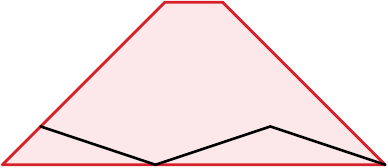}
\par{\scriptsize \mbox{$p=0.6,\ q=0.6,\ \alpha=45^\circ,\ \beta=45^\circ$}\par}
\end{subfigure}
\\
\begin{subfigure}[b]{0.326\linewidth}\centering
\includegraphics[width=\linewidth,height=119pt,keepaspectratio]{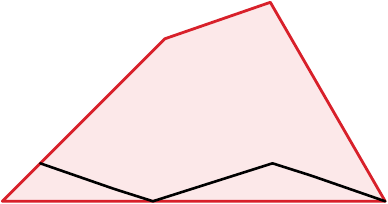}
\par{\scriptsize \mbox{$p=0.6,\ q=0.6,\ \alpha=45^\circ,\ \beta=60^\circ$}\par}
\end{subfigure}
&
\begin{subfigure}[b]{0.326\linewidth}\centering
\includegraphics[width=\linewidth,height=119pt,keepaspectratio]{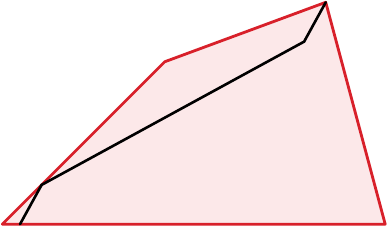}
\par{\scriptsize \mbox{$p=0.6,\ q=0.6,\ \alpha=45^\circ,\ \beta=75^\circ$}\par}
\end{subfigure}
&
\begin{subfigure}[b]{0.326\linewidth}\centering
\includegraphics[width=\linewidth,height=119pt,keepaspectratio]{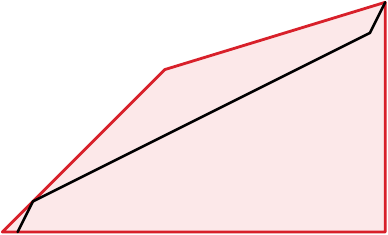}
\par{\scriptsize \mbox{$p=0.6,\ q=0.6,\ \alpha=45^\circ,\ \beta=90^\circ$}\par}
\end{subfigure}
\\
\begin{subfigure}[b]{0.326\linewidth}\centering
\includegraphics[width=\linewidth,height=119pt,keepaspectratio]{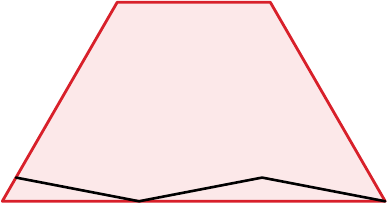}
\par{\scriptsize \mbox{$p=0.6,\ q=0.6,\ \alpha=60^\circ,\ \beta=60^\circ$}\par}
\end{subfigure}
&
\begin{subfigure}[b]{0.326\linewidth}\centering
\includegraphics[width=\linewidth,height=119pt,keepaspectratio]{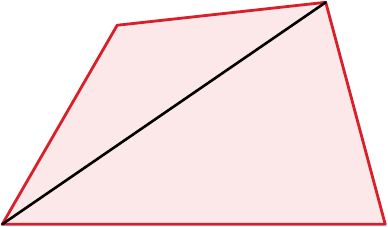}
\par{\scriptsize \mbox{$p=0.6,\ q=0.6,\ \alpha=60^\circ,\ \beta=75^\circ$}\par}
\end{subfigure}
&
\begin{subfigure}[b]{0.326\linewidth}\centering
\includegraphics[width=\linewidth,height=119pt,keepaspectratio]{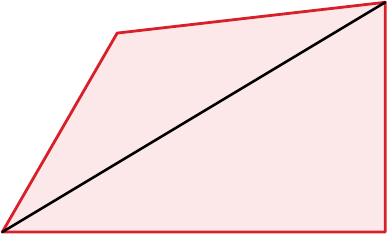}
\par{\scriptsize \mbox{$p=0.6,\ q=0.6,\ \alpha=60^\circ,\ \beta=90^\circ$}\par}
\end{subfigure}
\\
\begin{subfigure}[b]{0.326\linewidth}\centering
\includegraphics[width=\linewidth,height=119pt,keepaspectratio]{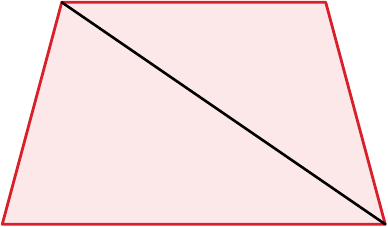}
\par{\scriptsize \mbox{$p=0.6,\ q=0.6,\ \alpha=75^\circ,\ \beta=75^\circ$}\par}
\end{subfigure}
&
\begin{subfigure}[b]{0.326\linewidth}\centering
\includegraphics[width=\linewidth,height=119pt,keepaspectratio]{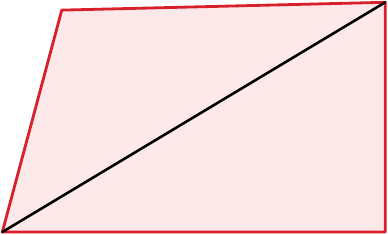}
\par{\scriptsize \mbox{$p=0.6,\ q=0.6,\ \alpha=75^\circ,\ \beta=90^\circ$}\par}
\end{subfigure}
&
\begin{subfigure}[b]{0.326\linewidth}\centering
\includegraphics[width=\linewidth,height=119pt,keepaspectratio]{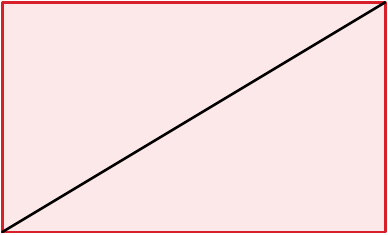}
\par{\scriptsize \mbox{$p=0.6,\ q=0.6,\ \alpha=90^\circ,\ \beta=90^\circ$}\par}
\end{subfigure}
\\
\begin{subfigure}[b]{0.326\linewidth}\centering
\includegraphics[width=\linewidth,height=119pt,keepaspectratio]{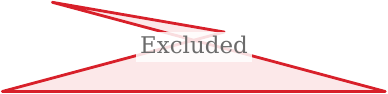}
\par{\scriptsize \mbox{$p=0.6,\ q=0.9,\ \alpha=15^\circ,\ \beta=15^\circ$}\par}
\end{subfigure}
&
\begin{subfigure}[b]{0.326\linewidth}\centering
\includegraphics[width=\linewidth,height=119pt,keepaspectratio]{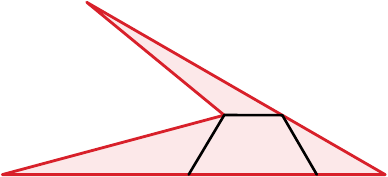}
\par{\scriptsize \mbox{$p=0.6,\ q=0.9,\ \alpha=15^\circ,\ \beta=30^\circ$}\par}
\end{subfigure}
&
\begin{subfigure}[b]{0.326\linewidth}\centering
\includegraphics[width=\linewidth,height=119pt,keepaspectratio]{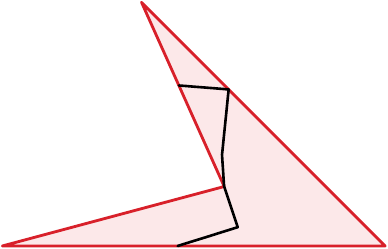}
\par{\scriptsize \mbox{$p=0.6,\ q=0.9,\ \alpha=15^\circ,\ \beta=45^\circ$}\par}
\end{subfigure}
\\
\begin{subfigure}[b]{0.326\linewidth}\centering
\includegraphics[width=\linewidth,height=119pt,keepaspectratio]{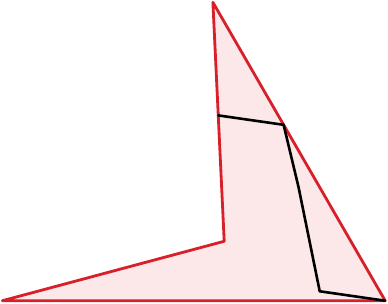}
\par{\scriptsize \mbox{$p=0.6,\ q=0.9,\ \alpha=15^\circ,\ \beta=60^\circ$}\par}
\end{subfigure}
&
\begin{subfigure}[b]{0.326\linewidth}\centering
\includegraphics[width=\linewidth,height=119pt,keepaspectratio]{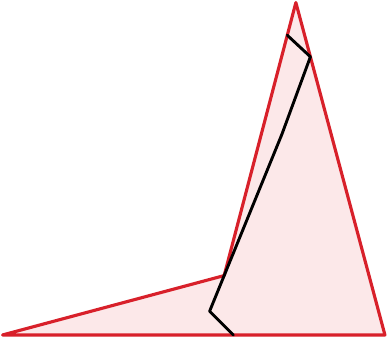}
\par{\scriptsize \mbox{$p=0.6,\ q=0.9,\ \alpha=15^\circ,\ \beta=75^\circ$}\par}
\end{subfigure}
&
\begin{subfigure}[b]{0.326\linewidth}\centering
\includegraphics[width=\linewidth,height=119pt,keepaspectratio]{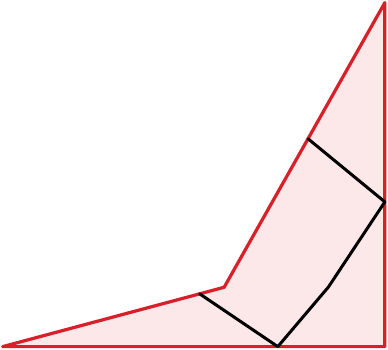}
\par{\scriptsize \mbox{$p=0.6,\ q=0.9,\ \alpha=15^\circ,\ \beta=90^\circ$}\par}
\end{subfigure}
\\
\begin{subfigure}[b]{0.326\linewidth}\centering
\includegraphics[width=\linewidth,height=119pt,keepaspectratio]{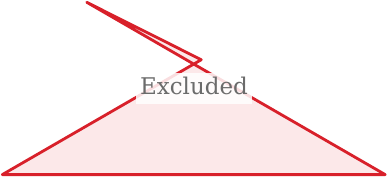}
\par{\scriptsize \mbox{$p=0.6,\ q=0.9,\ \alpha=30^\circ,\ \beta=30^\circ$}\par}
\end{subfigure}
&
\begin{subfigure}[b]{0.326\linewidth}\centering
\includegraphics[width=\linewidth,height=119pt,keepaspectratio]{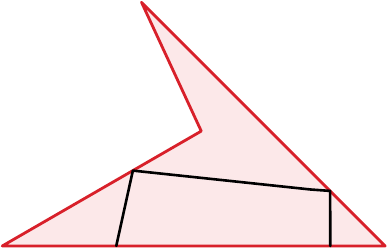}
\par{\scriptsize \mbox{$p=0.6,\ q=0.9,\ \alpha=30^\circ,\ \beta=45^\circ$}\par}
\end{subfigure}
&
\begin{subfigure}[b]{0.326\linewidth}\centering
\includegraphics[width=\linewidth,height=119pt,keepaspectratio]{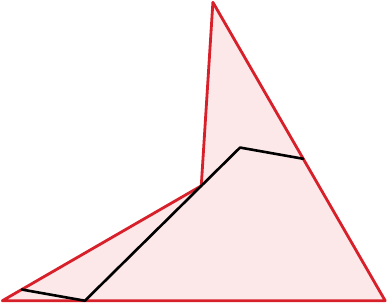}
\par{\scriptsize \mbox{$p=0.6,\ q=0.9,\ \alpha=30^\circ,\ \beta=60^\circ$}\par}
\end{subfigure}
\\
\begin{subfigure}[b]{0.326\linewidth}\centering
\includegraphics[width=\linewidth,height=119pt,keepaspectratio]{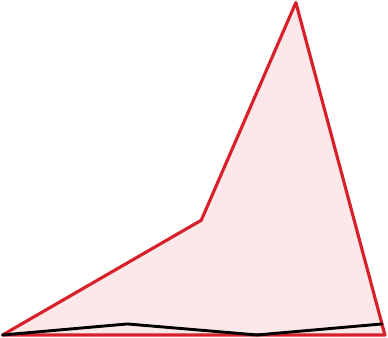}
\par{\scriptsize \mbox{$p=0.6,\ q=0.9,\ \alpha=30^\circ,\ \beta=75^\circ$}\par}
\end{subfigure}
&
\begin{subfigure}[b]{0.326\linewidth}\centering
\includegraphics[width=\linewidth,height=119pt,keepaspectratio]{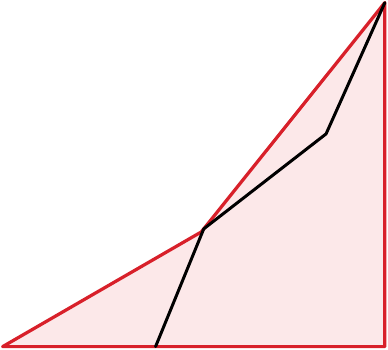}
\par{\scriptsize \mbox{$p=0.6,\ q=0.9,\ \alpha=30^\circ,\ \beta=90^\circ$}\par}
\end{subfigure}
&
\begin{subfigure}[b]{0.326\linewidth}\centering
\includegraphics[width=\linewidth,height=119pt,keepaspectratio]{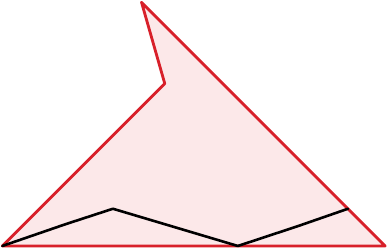}
\par{\scriptsize \mbox{$p=0.6,\ q=0.9,\ \alpha=45^\circ,\ \beta=45^\circ$}\par}
\end{subfigure}
\\
\begin{subfigure}[b]{0.326\linewidth}\centering
\includegraphics[width=\linewidth,height=119pt,keepaspectratio]{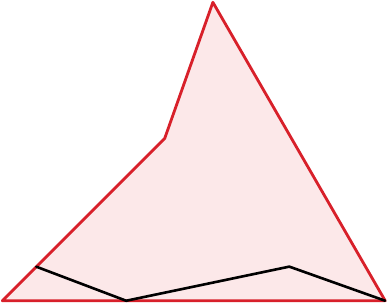}
\par{\scriptsize \mbox{$p=0.6,\ q=0.9,\ \alpha=45^\circ,\ \beta=60^\circ$}\par}
\end{subfigure}
&
\begin{subfigure}[b]{0.326\linewidth}\centering
\includegraphics[width=\linewidth,height=119pt,keepaspectratio]{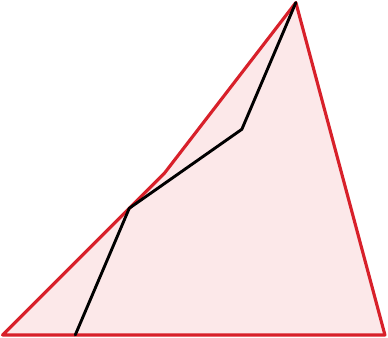}
\par{\scriptsize \mbox{$p=0.6,\ q=0.9,\ \alpha=45^\circ,\ \beta=75^\circ$}\par}
\end{subfigure}
&
\begin{subfigure}[b]{0.326\linewidth}\centering
\includegraphics[width=\linewidth,height=119pt,keepaspectratio]{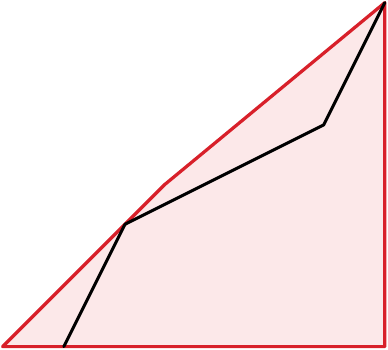}
\par{\scriptsize \mbox{$p=0.6,\ q=0.9,\ \alpha=45^\circ,\ \beta=90^\circ$}\par}
\end{subfigure}
\\
\begin{subfigure}[b]{0.326\linewidth}\centering
\includegraphics[width=\linewidth,height=119pt,keepaspectratio]{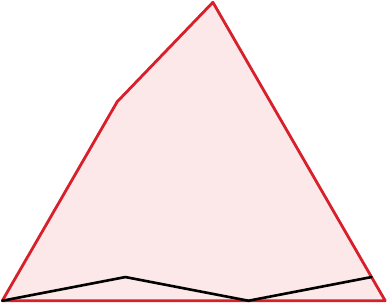}
\par{\scriptsize \mbox{$p=0.6,\ q=0.9,\ \alpha=60^\circ,\ \beta=60^\circ$}\par}
\end{subfigure}
&
\begin{subfigure}[b]{0.326\linewidth}\centering
\includegraphics[width=\linewidth,height=119pt,keepaspectratio]{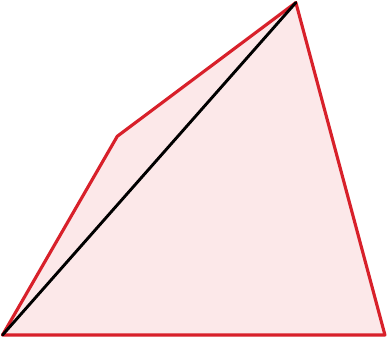}
\par{\scriptsize \mbox{$p=0.6,\ q=0.9,\ \alpha=60^\circ,\ \beta=75^\circ$}\par}
\end{subfigure}
&
\begin{subfigure}[b]{0.326\linewidth}\centering
\includegraphics[width=\linewidth,height=119pt,keepaspectratio]{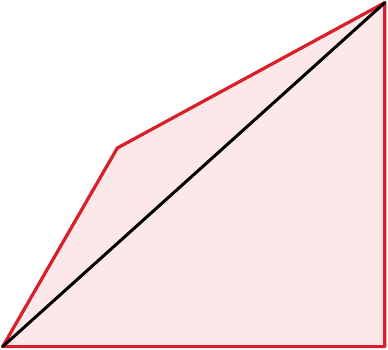}
\par{\scriptsize \mbox{$p=0.6,\ q=0.9,\ \alpha=60^\circ,\ \beta=90^\circ$}\par}
\end{subfigure}
\\
\begin{subfigure}[b]{0.326\linewidth}\centering
\includegraphics[width=\linewidth,height=119pt,keepaspectratio]{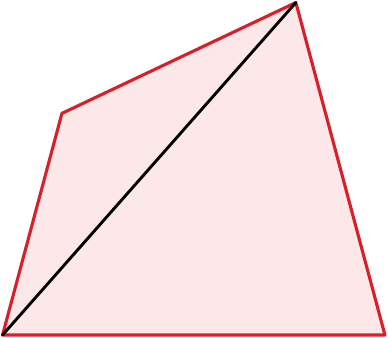}
\par{\scriptsize \mbox{$p=0.6,\ q=0.9,\ \alpha=75^\circ,\ \beta=75^\circ$}\par}
\end{subfigure}
&
\begin{subfigure}[b]{0.326\linewidth}\centering
\includegraphics[width=\linewidth,height=119pt,keepaspectratio]{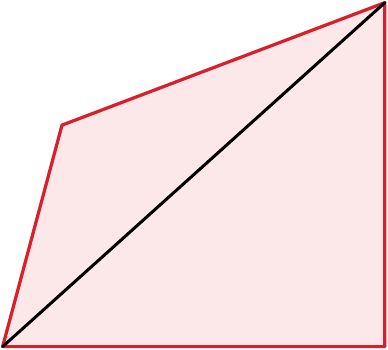}
\par{\scriptsize \mbox{$p=0.6,\ q=0.9,\ \alpha=75^\circ,\ \beta=90^\circ$}\par}
\end{subfigure}
&
\begin{subfigure}[b]{0.326\linewidth}\centering
\includegraphics[width=\linewidth,height=119pt,keepaspectratio]{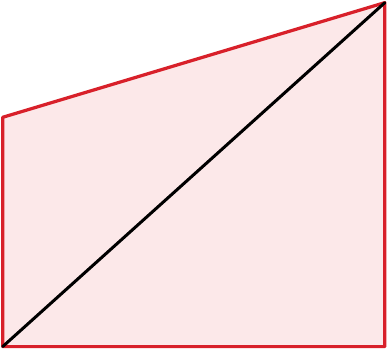}
\par{\scriptsize \mbox{$p=0.6,\ q=0.9,\ \alpha=90^\circ,\ \beta=90^\circ$}\par}
\end{subfigure}
\\
\begin{subfigure}[b]{0.326\linewidth}\centering
\includegraphics[width=\linewidth,height=119pt,keepaspectratio]{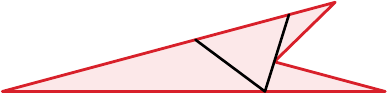}
\par{\scriptsize \mbox{$p=0.9,\ q=0.3,\ \alpha=15^\circ,\ \beta=15^\circ$}\par}
\end{subfigure}
&
\begin{subfigure}[b]{0.326\linewidth}\centering
\includegraphics[width=\linewidth,height=119pt,keepaspectratio]{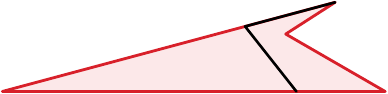}
\par{\scriptsize \mbox{$p=0.9,\ q=0.3,\ \alpha=15^\circ,\ \beta=30^\circ$}\par}
\end{subfigure}
&
\begin{subfigure}[b]{0.326\linewidth}\centering
\includegraphics[width=\linewidth,height=119pt,keepaspectratio]{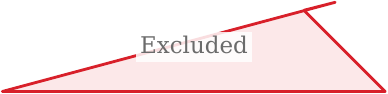}
\par{\scriptsize \mbox{$p=0.9,\ q=0.3,\ \alpha=15^\circ,\ \beta=45^\circ$}\par}
\end{subfigure}
\\
\begin{subfigure}[b]{0.326\linewidth}\centering
\includegraphics[width=\linewidth,height=119pt,keepaspectratio]{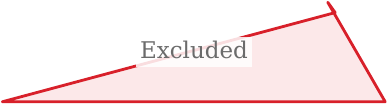}
\par{\scriptsize \mbox{$p=0.9,\ q=0.3,\ \alpha=15^\circ,\ \beta=60^\circ$}\par}
\end{subfigure}
&
\begin{subfigure}[b]{0.326\linewidth}\centering
\includegraphics[width=\linewidth,height=119pt,keepaspectratio]{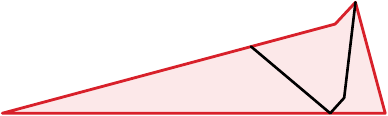}
\par{\scriptsize \mbox{$p=0.9,\ q=0.3,\ \alpha=15^\circ,\ \beta=75^\circ$}\par}
\end{subfigure}
&
\begin{subfigure}[b]{0.326\linewidth}\centering
\includegraphics[width=\linewidth,height=119pt,keepaspectratio]{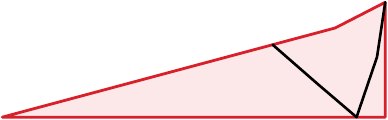}
\par{\scriptsize \mbox{$p=0.9,\ q=0.3,\ \alpha=15^\circ,\ \beta=90^\circ$}\par}
\end{subfigure}
\\
\begin{subfigure}[b]{0.326\linewidth}\centering
\includegraphics[width=\linewidth,height=119pt,keepaspectratio]{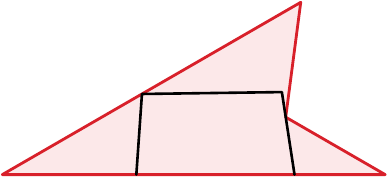}
\par{\scriptsize \mbox{$p=0.9,\ q=0.3,\ \alpha=30^\circ,\ \beta=30^\circ$}\par}
\end{subfigure}
&
\begin{subfigure}[b]{0.326\linewidth}\centering
\includegraphics[width=\linewidth,height=119pt,keepaspectratio]{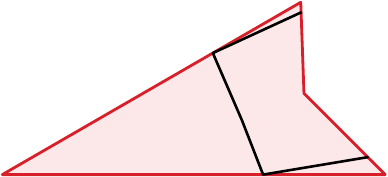}
\par{\scriptsize \mbox{$p=0.9,\ q=0.3,\ \alpha=30^\circ,\ \beta=45^\circ$}\par}
\end{subfigure}
&
\begin{subfigure}[b]{0.326\linewidth}\centering
\includegraphics[width=\linewidth,height=119pt,keepaspectratio]{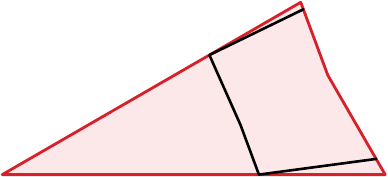}
\par{\scriptsize \mbox{$p=0.9,\ q=0.3,\ \alpha=30^\circ,\ \beta=60^\circ$}\par}
\end{subfigure}
\\
\begin{subfigure}[b]{0.326\linewidth}\centering
\includegraphics[width=\linewidth,height=119pt,keepaspectratio]{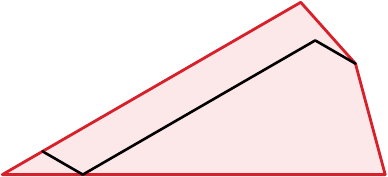}
\par{\scriptsize \mbox{$p=0.9,\ q=0.3,\ \alpha=30^\circ,\ \beta=75^\circ$}\par}
\end{subfigure}
&
\begin{subfigure}[b]{0.326\linewidth}\centering
\includegraphics[width=\linewidth,height=119pt,keepaspectratio]{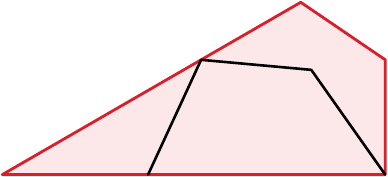}
\par{\scriptsize \mbox{$p=0.9,\ q=0.3,\ \alpha=30^\circ,\ \beta=90^\circ$}\par}
\end{subfigure}
&
\begin{subfigure}[b]{0.326\linewidth}\centering
\includegraphics[width=\linewidth,height=119pt,keepaspectratio]{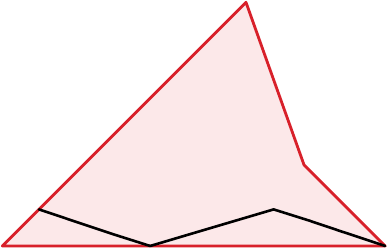}
\par{\scriptsize \mbox{$p=0.9,\ q=0.3,\ \alpha=45^\circ,\ \beta=45^\circ$}\par}
\end{subfigure}
\\
\begin{subfigure}[b]{0.326\linewidth}\centering
\includegraphics[width=\linewidth,height=119pt,keepaspectratio]{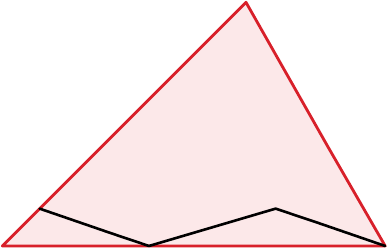}
\par{\scriptsize \mbox{$p=0.9,\ q=0.3,\ \alpha=45^\circ,\ \beta=60^\circ$}\par}
\end{subfigure}
&
\begin{subfigure}[b]{0.326\linewidth}\centering
\includegraphics[width=\linewidth,height=119pt,keepaspectratio]{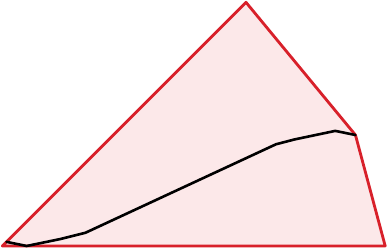}
\par{\scriptsize \mbox{$p=0.9,\ q=0.3,\ \alpha=45^\circ,\ \beta=75^\circ$}\par}
\end{subfigure}
&
\begin{subfigure}[b]{0.326\linewidth}\centering
\includegraphics[width=\linewidth,height=119pt,keepaspectratio]{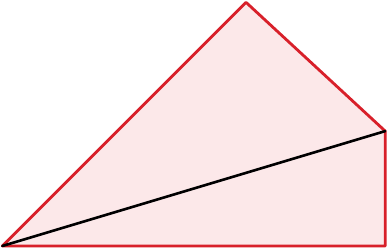}
\par{\scriptsize \mbox{$p=0.9,\ q=0.3,\ \alpha=45^\circ,\ \beta=90^\circ$}\par}
\end{subfigure}
\\
\begin{subfigure}[b]{0.326\linewidth}\centering
\includegraphics[width=\linewidth,height=119pt,keepaspectratio]{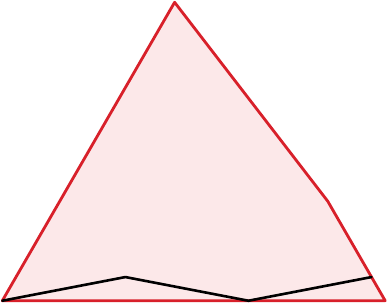}
\par{\scriptsize \mbox{$p=0.9,\ q=0.3,\ \alpha=60^\circ,\ \beta=60^\circ$}\par}
\end{subfigure}
&
\begin{subfigure}[b]{0.326\linewidth}\centering
\includegraphics[width=\linewidth,height=119pt,keepaspectratio]{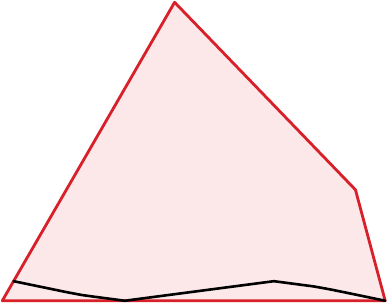}
\par{\scriptsize \mbox{$p=0.9,\ q=0.3,\ \alpha=60^\circ,\ \beta=75^\circ$}\par}
\end{subfigure}
&
\begin{subfigure}[b]{0.326\linewidth}\centering
\includegraphics[width=\linewidth,height=119pt,keepaspectratio]{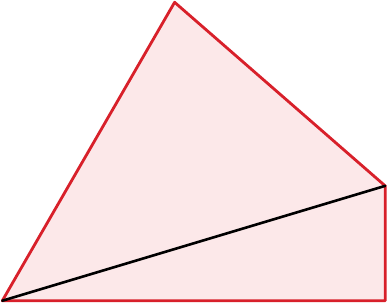}
\par{\scriptsize \mbox{$p=0.9,\ q=0.3,\ \alpha=60^\circ,\ \beta=90^\circ$}\par}
\end{subfigure}
\\
\begin{subfigure}[b]{0.326\linewidth}\centering
\includegraphics[width=\linewidth,height=119pt,keepaspectratio]{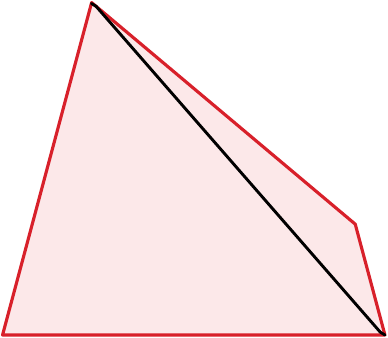}
\par{\scriptsize \mbox{$p=0.9,\ q=0.3,\ \alpha=75^\circ,\ \beta=75^\circ$}\par}
\end{subfigure}
&
\begin{subfigure}[b]{0.326\linewidth}\centering
\includegraphics[width=\linewidth,height=119pt,keepaspectratio]{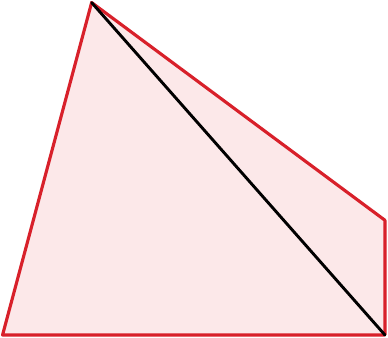}
\par{\scriptsize \mbox{$p=0.9,\ q=0.3,\ \alpha=75^\circ,\ \beta=90^\circ$}\par}
\end{subfigure}
&
\begin{subfigure}[b]{0.326\linewidth}\centering
\includegraphics[width=\linewidth,height=119pt,keepaspectratio]{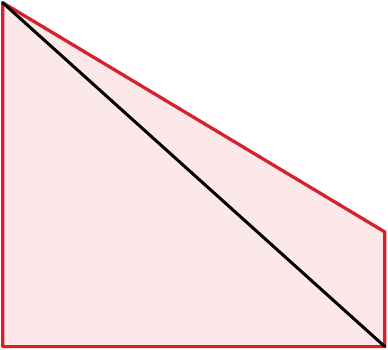}
\par{\scriptsize \mbox{$p=0.9,\ q=0.3,\ \alpha=90^\circ,\ \beta=90^\circ$}\par}
\end{subfigure}
\\
\begin{subfigure}[b]{0.326\linewidth}\centering
\includegraphics[width=\linewidth,height=119pt,keepaspectratio]{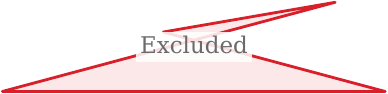}
\par{\scriptsize \mbox{$p=0.9,\ q=0.6,\ \alpha=15^\circ,\ \beta=15^\circ$}\par}
\end{subfigure}
&
\begin{subfigure}[b]{0.326\linewidth}\centering
\includegraphics[width=\linewidth,height=119pt,keepaspectratio]{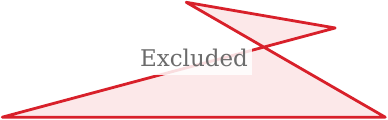}
\par{\scriptsize \mbox{$p=0.9,\ q=0.6,\ \alpha=15^\circ,\ \beta=30^\circ$}\par}
\end{subfigure}
&
\begin{subfigure}[b]{0.326\linewidth}\centering
\includegraphics[width=\linewidth,height=119pt,keepaspectratio]{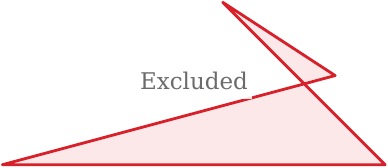}
\par{\scriptsize \mbox{$p=0.9,\ q=0.6,\ \alpha=15^\circ,\ \beta=45^\circ$}\par}
\end{subfigure}
\\
\begin{subfigure}[b]{0.326\linewidth}\centering
\includegraphics[width=\linewidth,height=119pt,keepaspectratio]{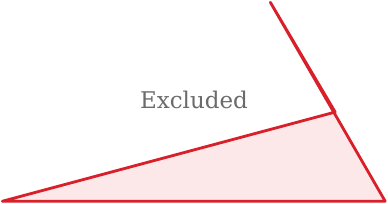}
\par{\scriptsize \mbox{$p=0.9,\ q=0.6,\ \alpha=15^\circ,\ \beta=60^\circ$}\par}
\end{subfigure}
&
\begin{subfigure}[b]{0.326\linewidth}\centering
\includegraphics[width=\linewidth,height=119pt,keepaspectratio]{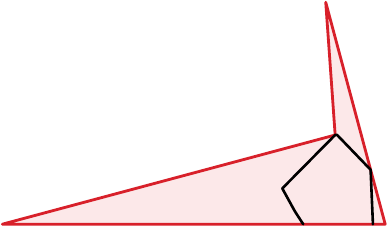}
\par{\scriptsize \mbox{$p=0.9,\ q=0.6,\ \alpha=15^\circ,\ \beta=75^\circ$}\par}
\end{subfigure}
&
\begin{subfigure}[b]{0.326\linewidth}\centering
\includegraphics[width=\linewidth,height=119pt,keepaspectratio]{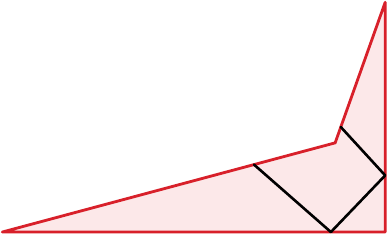}
\par{\scriptsize \mbox{$p=0.9,\ q=0.6,\ \alpha=15^\circ,\ \beta=90^\circ$}\par}
\end{subfigure}
\\
\begin{subfigure}[b]{0.326\linewidth}\centering
\includegraphics[width=\linewidth,height=119pt,keepaspectratio]{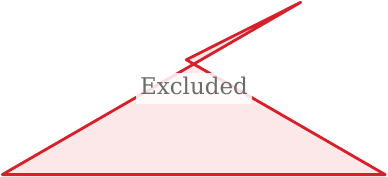}
\par{\scriptsize \mbox{$p=0.9,\ q=0.6,\ \alpha=30^\circ,\ \beta=30^\circ$}\par}
\end{subfigure}
&
\begin{subfigure}[b]{0.326\linewidth}\centering
\includegraphics[width=\linewidth,height=119pt,keepaspectratio]{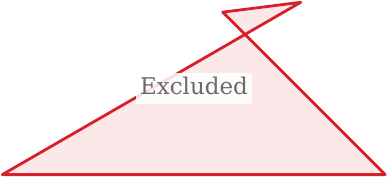}
\par{\scriptsize \mbox{$p=0.9,\ q=0.6,\ \alpha=30^\circ,\ \beta=45^\circ$}\par}
\end{subfigure}
&
\begin{subfigure}[b]{0.326\linewidth}\centering
\includegraphics[width=\linewidth,height=119pt,keepaspectratio]{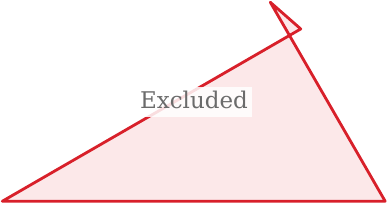}
\par{\scriptsize \mbox{$p=0.9,\ q=0.6,\ \alpha=30^\circ,\ \beta=60^\circ$}\par}
\end{subfigure}
\\
\begin{subfigure}[b]{0.326\linewidth}\centering
\includegraphics[width=\linewidth,height=119pt,keepaspectratio]{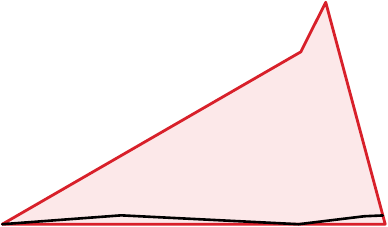}
\par{\scriptsize \mbox{$p=0.9,\ q=0.6,\ \alpha=30^\circ,\ \beta=75^\circ$}\par}
\end{subfigure}
&
\begin{subfigure}[b]{0.326\linewidth}\centering
\includegraphics[width=\linewidth,height=119pt,keepaspectratio]{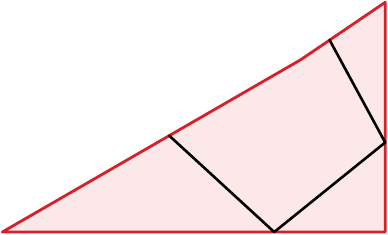}
\par{\scriptsize \mbox{$p=0.9,\ q=0.6,\ \alpha=30^\circ,\ \beta=90^\circ$}\par}
\end{subfigure}
&
\begin{subfigure}[b]{0.326\linewidth}\centering
\includegraphics[width=\linewidth,height=119pt,keepaspectratio]{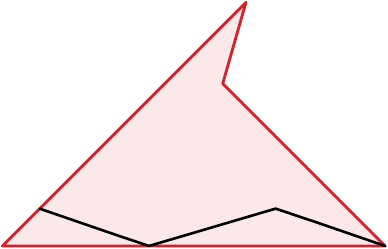}
\par{\scriptsize \mbox{$p=0.9,\ q=0.6,\ \alpha=45^\circ,\ \beta=45^\circ$}\par}
\end{subfigure}
\\
\begin{subfigure}[b]{0.326\linewidth}\centering
\includegraphics[width=\linewidth,height=119pt,keepaspectratio]{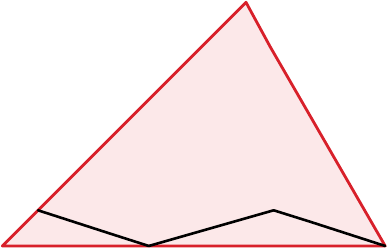}
\par{\scriptsize \mbox{$p=0.9,\ q=0.6,\ \alpha=45^\circ,\ \beta=60^\circ$}\par}
\end{subfigure}
&
\begin{subfigure}[b]{0.326\linewidth}\centering
\includegraphics[width=\linewidth,height=119pt,keepaspectratio]{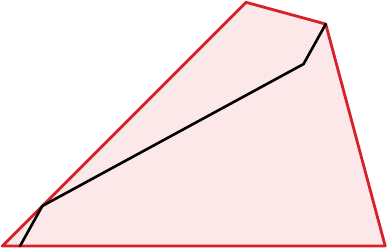}
\par{\scriptsize \mbox{$p=0.9,\ q=0.6,\ \alpha=45^\circ,\ \beta=75^\circ$}\par}
\end{subfigure}
&
\begin{subfigure}[b]{0.326\linewidth}\centering
\includegraphics[width=\linewidth,height=119pt,keepaspectratio]{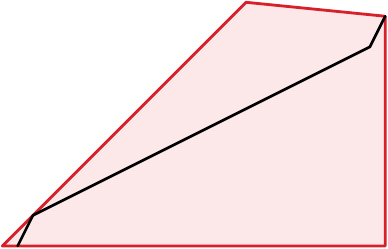}
\par{\scriptsize \mbox{$p=0.9,\ q=0.6,\ \alpha=45^\circ,\ \beta=90^\circ$}\par}
\end{subfigure}
\\
\begin{subfigure}[b]{0.326\linewidth}\centering
\includegraphics[width=\linewidth,height=119pt,keepaspectratio]{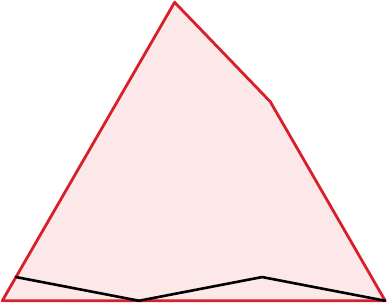}
\par{\scriptsize \mbox{$p=0.9,\ q=0.6,\ \alpha=60^\circ,\ \beta=60^\circ$}\par}
\end{subfigure}
&
\begin{subfigure}[b]{0.326\linewidth}\centering
\includegraphics[width=\linewidth,height=119pt,keepaspectratio]{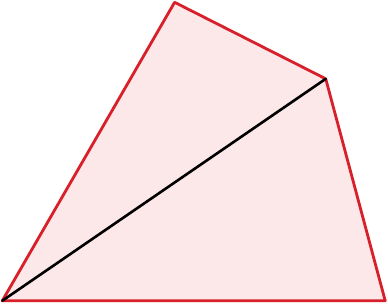}
\par{\scriptsize \mbox{$p=0.9,\ q=0.6,\ \alpha=60^\circ,\ \beta=75^\circ$}\par}
\end{subfigure}
&
\begin{subfigure}[b]{0.326\linewidth}\centering
\includegraphics[width=\linewidth,height=119pt,keepaspectratio]{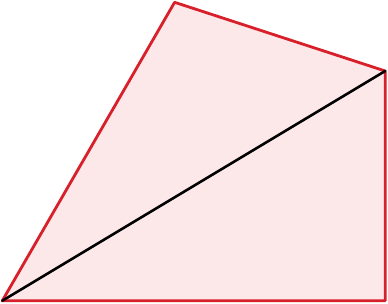}
\par{\scriptsize \mbox{$p=0.9,\ q=0.6,\ \alpha=60^\circ,\ \beta=90^\circ$}\par}
\end{subfigure}
\\
\begin{subfigure}[b]{0.326\linewidth}\centering
\includegraphics[width=\linewidth,height=119pt,keepaspectratio]{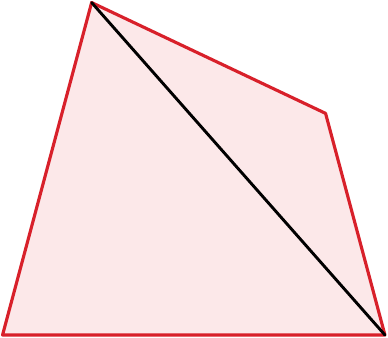}
\par{\scriptsize \mbox{$p=0.9,\ q=0.6,\ \alpha=75^\circ,\ \beta=75^\circ$}\par}
\end{subfigure}
&
\begin{subfigure}[b]{0.326\linewidth}\centering
\includegraphics[width=\linewidth,height=119pt,keepaspectratio]{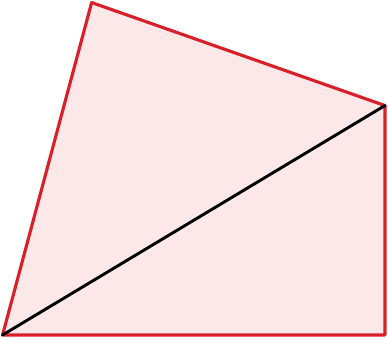}
\par{\scriptsize \mbox{$p=0.9,\ q=0.6,\ \alpha=75^\circ,\ \beta=90^\circ$}\par}
\end{subfigure}
&
\begin{subfigure}[b]{0.326\linewidth}\centering
\includegraphics[width=\linewidth,height=119pt,keepaspectratio]{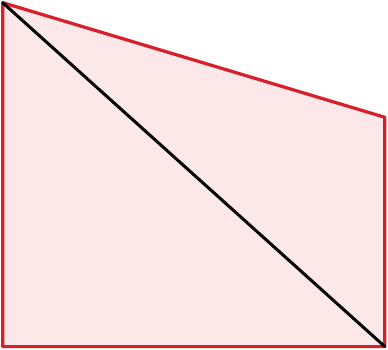}
\par{\scriptsize \mbox{$p=0.9,\ q=0.6,\ \alpha=90^\circ,\ \beta=90^\circ$}\par}
\end{subfigure}
\\
\begin{subfigure}[b]{0.326\linewidth}\centering
\includegraphics[width=\linewidth,height=119pt,keepaspectratio]{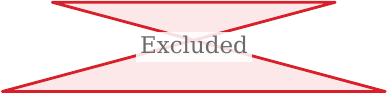}
\par{\scriptsize \mbox{$p=0.9,\ q=0.9,\ \alpha=15^\circ,\ \beta=15^\circ$}\par}
\end{subfigure}
&
\begin{subfigure}[b]{0.326\linewidth}\centering
\includegraphics[width=\linewidth,height=119pt,keepaspectratio]{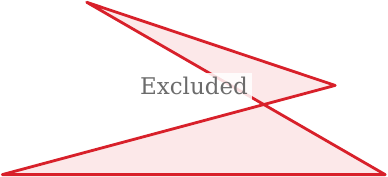}
\par{\scriptsize \mbox{$p=0.9,\ q=0.9,\ \alpha=15^\circ,\ \beta=30^\circ$}\par}
\end{subfigure}
&
\begin{subfigure}[b]{0.326\linewidth}\centering
\includegraphics[width=\linewidth,height=119pt,keepaspectratio]{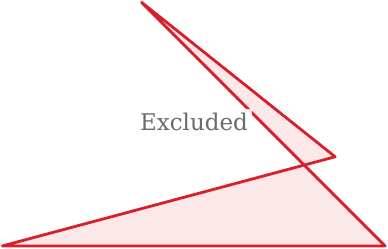}
\par{\scriptsize \mbox{$p=0.9,\ q=0.9,\ \alpha=15^\circ,\ \beta=45^\circ$}\par}
\end{subfigure}
\\
\begin{subfigure}[b]{0.326\linewidth}\centering
\includegraphics[width=\linewidth,height=119pt,keepaspectratio]{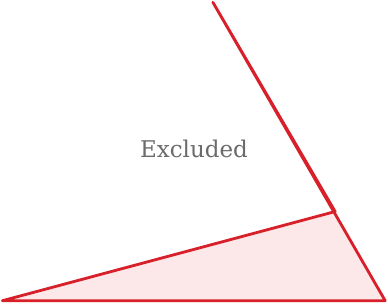}
\par{\scriptsize \mbox{$p=0.9,\ q=0.9,\ \alpha=15^\circ,\ \beta=60^\circ$}\par}
\end{subfigure}
&
\begin{subfigure}[b]{0.326\linewidth}\centering
\includegraphics[width=\linewidth,height=119pt,keepaspectratio]{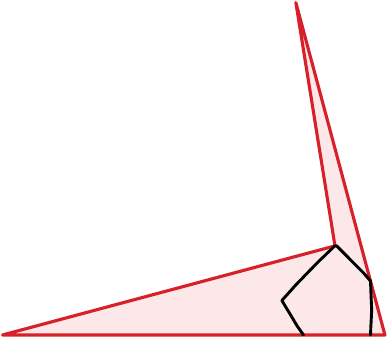}
\par{\scriptsize \mbox{$p=0.9,\ q=0.9,\ \alpha=15^\circ,\ \beta=75^\circ$}\par}
\end{subfigure}
&
\begin{subfigure}[b]{0.326\linewidth}\centering
\includegraphics[width=\linewidth,height=119pt,keepaspectratio]{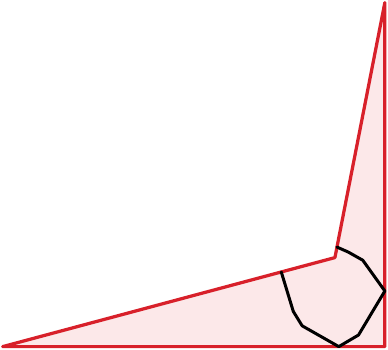}
\par{\scriptsize \mbox{$p=0.9,\ q=0.9,\ \alpha=15^\circ,\ \beta=90^\circ$}\par}
\end{subfigure}
\\
\begin{subfigure}[b]{0.326\linewidth}\centering
\includegraphics[width=\linewidth,height=119pt,keepaspectratio]{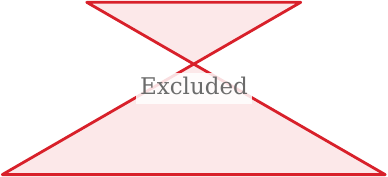}
\par{\scriptsize \mbox{$p=0.9,\ q=0.9,\ \alpha=30^\circ,\ \beta=30^\circ$}\par}
\end{subfigure}
&
\begin{subfigure}[b]{0.326\linewidth}\centering
\includegraphics[width=\linewidth,height=119pt,keepaspectratio]{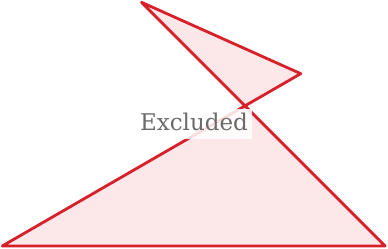}
\par{\scriptsize \mbox{$p=0.9,\ q=0.9,\ \alpha=30^\circ,\ \beta=45^\circ$}\par}
\end{subfigure}
&
\begin{subfigure}[b]{0.326\linewidth}\centering
\includegraphics[width=\linewidth,height=119pt,keepaspectratio]{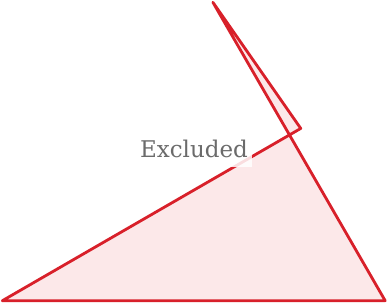}
\par{\scriptsize \mbox{$p=0.9,\ q=0.9,\ \alpha=30^\circ,\ \beta=60^\circ$}\par}
\end{subfigure}
\\
\begin{subfigure}[b]{0.326\linewidth}\centering
\includegraphics[width=\linewidth,height=119pt,keepaspectratio]{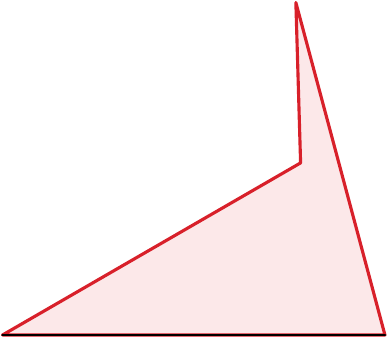}
\par{\scriptsize \mbox{$p=0.9,\ q=0.9,\ \alpha=30^\circ,\ \beta=75^\circ$}\par}
\end{subfigure}
&
\begin{subfigure}[b]{0.326\linewidth}\centering
\includegraphics[width=\linewidth,height=119pt,keepaspectratio]{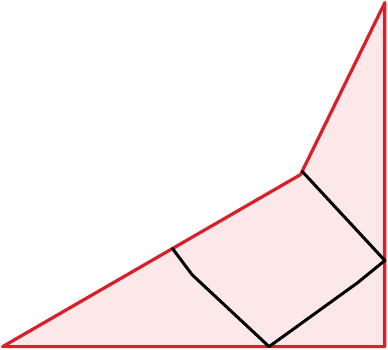}
\par{\scriptsize \mbox{$p=0.9,\ q=0.9,\ \alpha=30^\circ,\ \beta=90^\circ$}\par}
\end{subfigure}
&
\begin{subfigure}[b]{0.326\linewidth}\centering
\includegraphics[width=\linewidth,height=119pt,keepaspectratio]{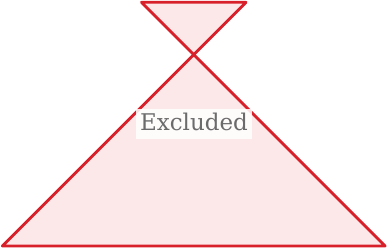}
\par{\scriptsize \mbox{$p=0.9,\ q=0.9,\ \alpha=45^\circ,\ \beta=45^\circ$}\par}
\end{subfigure}
\\
\begin{subfigure}[b]{0.326\linewidth}\centering
\includegraphics[width=\linewidth,height=119pt,keepaspectratio]{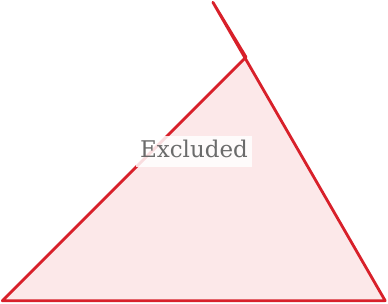}
\par{\scriptsize \mbox{$p=0.9,\ q=0.9,\ \alpha=45^\circ,\ \beta=60^\circ$}\par}
\end{subfigure}
&
\begin{subfigure}[b]{0.326\linewidth}\centering
\includegraphics[width=\linewidth,height=119pt,keepaspectratio]{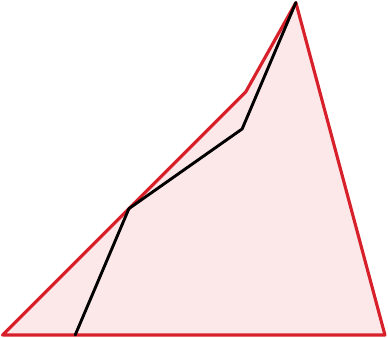}
\par{\scriptsize \mbox{$p=0.9,\ q=0.9,\ \alpha=45^\circ,\ \beta=75^\circ$}\par}
\end{subfigure}
&
\begin{subfigure}[b]{0.326\linewidth}\centering
\includegraphics[width=\linewidth,height=119pt,keepaspectratio]{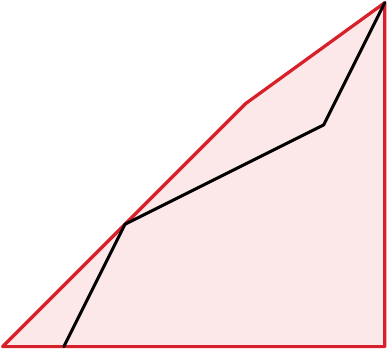}
\par{\scriptsize \mbox{$p=0.9,\ q=0.9,\ \alpha=45^\circ,\ \beta=90^\circ$}\par}
\end{subfigure}
\\
\begin{subfigure}[b]{0.326\linewidth}\centering
\includegraphics[width=\linewidth,height=119pt,keepaspectratio]{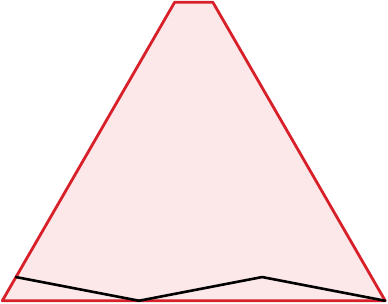}
\par{\scriptsize \mbox{$p=0.9,\ q=0.9,\ \alpha=60^\circ,\ \beta=60^\circ$}\par}
\end{subfigure}
&
\begin{subfigure}[b]{0.326\linewidth}\centering
\includegraphics[width=\linewidth,height=119pt,keepaspectratio]{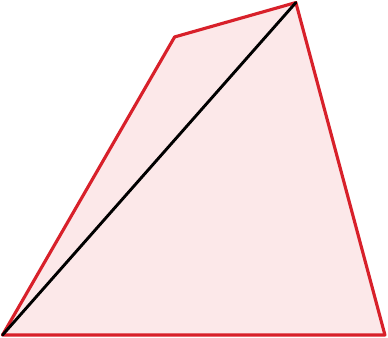}
\par{\scriptsize \mbox{$p=0.9,\ q=0.9,\ \alpha=60^\circ,\ \beta=75^\circ$}\par}
\end{subfigure}
&
\begin{subfigure}[b]{0.326\linewidth}\centering
\includegraphics[width=\linewidth,height=119pt,keepaspectratio]{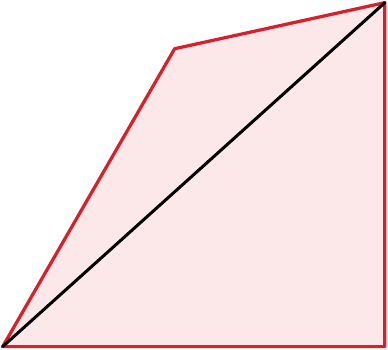}
\par{\scriptsize \mbox{$p=0.9,\ q=0.9,\ \alpha=60^\circ,\ \beta=90^\circ$}\par}
\end{subfigure}
\\
\begin{subfigure}[b]{0.326\linewidth}\centering
\includegraphics[width=\linewidth,height=119pt,keepaspectratio]{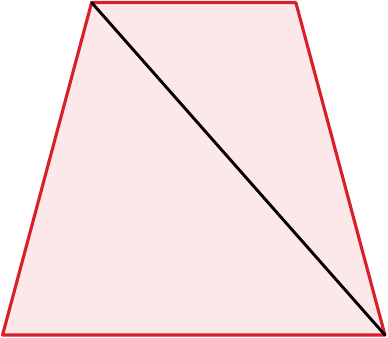}
\par{\scriptsize \mbox{$p=0.9,\ q=0.9,\ \alpha=75^\circ,\ \beta=75^\circ$}\par}
\end{subfigure}
&
\begin{subfigure}[b]{0.326\linewidth}\centering
\includegraphics[width=\linewidth,height=119pt,keepaspectratio]{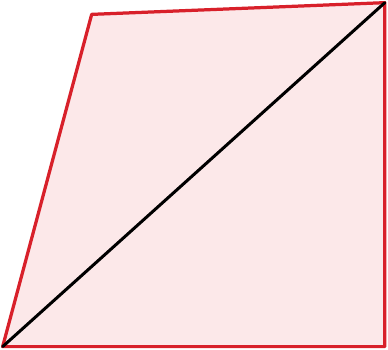}
\par{\scriptsize \mbox{$p=0.9,\ q=0.9,\ \alpha=75^\circ,\ \beta=90^\circ$}\par}
\end{subfigure}
&
\begin{subfigure}[b]{0.326\linewidth}\centering
\includegraphics[width=\linewidth,height=119pt,keepaspectratio]{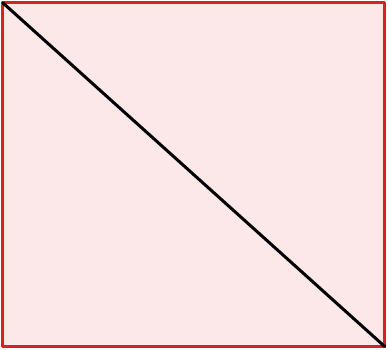}
\par{\scriptsize \mbox{$p=0.9,\ q=0.9,\ \alpha=90^\circ,\ \beta=90^\circ$}\par}
\end{subfigure}
\\
\end{longtable}
\endgroup

\section{Conclusion}
The numerical investigation covers five structured and general quadrilateral families. It reproduces the square and rectangle benchmarks, resolves competing bent and straight branches near the rhombus transition, and shows why a strip/diagonal comparison alone is insufficient for the broader parallelogram study. The trapezoid and arbitrary galleries retain the actual polygon geometry and consistently distinguish valid, limiting, and excluded configurations. 

These numerical results provide a reproducible basis for comparing quadrilateral escape paths and for further refinement of the most sensitive parameter regions. More refined numerical calculations and more choices for side lengths and angles can yield more accurate and detailed results of escape path for quadrilaterals. Future work includes a rigorous proof of the globally optimal shortest path for different types of quadrilateral and escape path shapes.

{\raggedright
}
\bigskip

AI usage disclosure: The author provided the methodological framework, while GPT-6 Astra handled numerical calculations and polished the language. 
\bigskip

\noindent College of Engineering and Computer Science, University of Central Florida, Orlando, FL, USA

\noindent Email: zhipeng.deng@ucf.edu
\end{document}